\documentclass{amsart}
\usepackage{moreverb,url}
\usepackage{amssymb,amsmath,amsfonts}
\usepackage{bm}
\usepackage[font=scriptsize, labelfont=bf]{caption}
\usepackage{float}
\usepackage{varwidth}
\usepackage{booktabs}
\usepackage{afterpage}
\usepackage{mathrsfs}
\usepackage{subcaption}
\usepackage[squaren]{SIunits}
\usepackage{verbatim}
\usepackage[square,numbers]{natbib}
\usepackage{stmaryrd}
\usepackage{fancyhdr}
\usepackage{syntax,etoolbox}
\usepackage{graphicx}
\usepackage{grffile}
\usepackage{hyperref}
\hypersetup{colorlinks=true,allcolors=blue}
\usepackage[square,numbers]{natbib}
\usepackage[left=1in, right=1in, bottom=1in,top=1in]{geometry}
\usepackage{lineno}
\usepackage{syntax,etoolbox}
\usepackage{algorithmicx}
\usepackage{algorithm}
\usepackage{algpseudocode}
\usepackage{tikz}
\algnewcommand{\Inputs}[1]{%
	\State \textbf{Input:}
	\Statex \hspace*{\algorithmicindent}\parbox[t]{.8\linewidth}{\raggedright #1}
}
\algnewcommand{\Initialize}[1]{%
	\State \textbf{Initialize:}
	\Statex \hspace*{\algorithmicindent}\parbox[t]{.8\linewidth}{\raggedright #1}
}
\algnewcommand{\Output}[1]{%
	\State \textbf{Output:}
	\Statex \hspace*{\algorithmicindent}\parbox[t]{.8\linewidth}{\raggedright #1}
}

\usepackage{hypcap}
\usepackage{leftidx}
\usepackage{xcolor}
\usepackage[nomain,symbols,nogroupskip]{glossaries}
\hypersetup{colorlinks=true,allcolors=blue}
\setglossarystyle{list}
\newglossary[mlg]{mathsymbols}{msy}{msn}{List of Mathematical Symbols}
\newglossary[plg]{physsymbols}{psy}{psn}{List of Physics Symbols}
\makeglossaries
\newglossaryentry{Omega}{
	name=\ensuremath{\Omega},
	description={open set or domain in $\mathbb{R}^3$},
	type=mathsymbols,
}
\newglossaryentry{DOmega}{
	name=\ensuremath{\mathcal D(\Omega)},
	description={space of infinitely differentiable functions with compact support in $\Omega$},
	type=mathsymbols,
}
\newglossaryentry{f+}{
	name=\ensuremath{\{f\}^{+}},
	description={positive part of the function $f$},
	type=mathsymbols,
}
\newglossaryentry{f-}{
	name=\ensuremath{\{f\}^{-}},
	description={negative part of the function $f$},
	type=mathsymbols,
}
\newglossaryentry{Gamma}{
	name=\ensuremath{\Gamma},
	description={topological boundary of a Lipschitz domain $\Omega\subset\mathbb{R}^3$},
	type=mathsymbols,
}
\newglossaryentry{dGamma}{
	name=\ensuremath{\mathrm{d}\Gamma},
	description={area element on $\Gamma$},
	type=mathsymbols,
}
\newglossaryentry{subsub}{
	name=\ensuremath{\subset\subset},
	description={set inclusion with positive distance between the boundaries},
	type=mathsymbols,
}
\newglossaryentry{VOmega}{
	name=\ensuremath{\bm{V}(\Omega)},
	description={space of admissible displacements in linearised elasticity},
	type=mathsymbols,
}
\newglossaryentry{A2}{
	name=\ensuremath{\mathcal{A}},
	description={linear elasticity operator},
	type=mathsymbols,
}
\newglossaryentry{pikappa}{
	name=\ensuremath{\pi_\kappa},
	description={tangential projection},
	type=mathsymbols,
}
\newglossaryentry{omegaproj}{
	name=\ensuremath{\omega},
	description={overlap set},
	type=mathsymbols,
}
\newglossaryentry{Falpha}{
	name=\ensuremath{F^{(\alpha)}(\xi)},
	description={open vertical fibre of body $\overline{\Omega^{(\alpha)}}$ above or below $\xi$},
	type=mathsymbols,
}
\newglossaryentry{U}{
	name=\ensuremath{\vec{\bm{U}}},
	description={set of admissible constrained displacements},
	type=mathsymbols,
}
\newglossaryentry{S}{
	name=\ensuremath{S},
	description={interface set},
	type=mathsymbols,
}
\newglossaryentry{Phi}{
	name=\ensuremath{\Phi},
	description={constraint qualification mapping},
	type=mathsymbols,
}
\newglossaryentry{KO1psi}{
	name=\ensuremath{\bm{K}(\Omega^{(1)};\Psi)},
	description={admissible displacements for the body $\overline{\Omega^{(1)}}$ under the virtual interface constraint},
	type=mathsymbols,
}
\newglossaryentry{KO2psi}{
	name=\ensuremath{\bm{K}(\Omega^{(2)};\Psi)},
	description={admissible displacements for the body $\overline{\Omega^{(2)}}$ under the virtual interface constraint},
	type=mathsymbols,
}
\newglossaryentry{VK}{
	name=\ensuremath{\mathcal{V}(K)},
	description={space of admissible controls},
	type=mathsymbols,
}
\newglossaryentry{Jbeta}{
	name=\ensuremath{\mathcal{J}_\beta},
	description={reduced cost functional characterising Problem~\ref{controlproblem}},
	type=mathsymbols,
}
\newglossaryentry{marginal}{
	name=\ensuremath{\tilde{\mathcal{F}}_\beta},
	description={marginal functional},
	type=mathsymbols,
}
\newglossaryentry{gammalimit}{
	name=\ensuremath{\tilde{\mathcal{F}}_0},
	description={$\Gamma$-limit of the family $\{\tilde{\mathcal{F}}_\beta\}_{\beta>0}$ in the strong topology of $\bm{L}^2(\Omega^{(1)})\times\bm{L}^2(\Omega^{(2)})$},
	type=mathsymbols,
}
\newglossaryentry{rho}{
	name=\ensuremath{\rho},
	description={positive penalty parameter meant to approach zero},
	type=mathsymbols,
}
\newglossaryentry{beta1}{
	name=\ensuremath{\bm{\beta}_1(\cdot;\Psi)},
	description={penalty operator associated with the deformation of the body $\overline{\Omega^{(1)}}$},
	type=mathsymbols,
}
\newglossaryentry{beta2}{
	name=\ensuremath{\bm{\beta}_2(\cdot;\Psi)},
	description={penalty operator associated with the deformation of the body $\overline{\Omega^{(2)}}$},
	type=mathsymbols,
}
\newglossaryentry{Jrho}{
	name=\ensuremath{\mathcal{J}_\rho},
	description={penalised reduced cost functional characterising Problem~\ref{cprho}},
	type=mathsymbols,
}
\newglossaryentry{c}{
	name=\ensuremath{c},
	description={positive penalty parameter in the Moreau--Yosida approximation meant to approach zero},
	type=mathsymbols,
}
\newglossaryentry{beta1c}{
	name=\ensuremath{\bm{\beta}_{1,c}(\cdot;\Psi)},
	description={smoothened penalty operator associated with the deformation of the body $\overline{\Omega^{(1)}}$},
	type=mathsymbols,
}
\newglossaryentry{beta2c}{
	name=\ensuremath{\bm{\beta}_{2,c}(\cdot;\Psi)},
	description={smoothened penalty operator associated with the deformation of the body $\overline{\Omega^{(2)}}$},
	type=mathsymbols,
}
\newglossaryentry{Jrhoc}{
	name=\ensuremath{\mathcal{J}_{\rho,c}},
	description={penalised reduced cost functional characterising Problem~\ref{cprhoc}},
	type=mathsymbols,
}
\newglossaryentry{A}{
	name=\ensuremath{A^{ijk\ell}},
	description={fourth order three-dimensional elasticity tensor},
	type=physsymbols,
}
\newglossaryentry{E3}{
	name=\ensuremath{\mathbb{E}^3},
	description={three-dimensional Euclidean space},
	type=physsymbols,
}
\newglossaryentry{ei}{
	name=\ensuremath{\bm{e}^i},
	description={reference frame in $\mathbb{R}^3$},
	type=physsymbols,
}
\newglossaryentry{lambda}{
	name=\ensuremath{\lambda},
	description={First Lam\'{e} constant},
	type=physsymbols,
}
\newglossaryentry{mu}{
	name=\ensuremath{\mu},
	description={Second Lam\'{e} constant (or shear modulus)},
	type=physsymbols,
}
\newglossaryentry{IOmega}{
	name=\ensuremath{\mathbb{I}_\Omega},
	description={identity function over the three-dimensional Lipschitz domain $\overline{\Omega}$},
	type=physsymbols,
}
\newglossaryentry{e}{
	name=\ensuremath{\bm{e}(\bm{v})=(e_{i\|j}(\bm{v}))},
	description={linearised strain tensor associated with the displacement field $\bm{v}$},
	type=physsymbols,
}
\newglossaryentry{Ce}{
	name=\ensuremath{C_e},
	description={positive constant associated with the uniform positive definiteness of the fourth order three-dimensional elasticity tensor $\{A^{ijk\ell}\}$},
	type=physsymbols,
}
\newglossaryentry{Omega1}{
	name=\ensuremath{\overline{\Omega^{(1)}}},
	description={reference configuration of the upper body},
	type=physsymbols,
}
\newglossaryentry{Omega2}{
	name=\ensuremath{\overline{\Omega^{(2)}}},
	description={reference configuration of the lower body},
	type=physsymbols,
}
\newglossaryentry{K}{
	name=\ensuremath{K},
	description={three-dimensional subset of $\mathbb{R}^3$ containing $\overline{\Omega^{(1)}}$ and $\overline{\Omega^{(2)}}$ and such that its boundary has positive distance from the boundaries of the reference configurations of the linearly elastic bodies $\overline{\Omega^{(1)}}$ and $\overline{\Omega^{(2)}}$},
	type=physsymbols,
}
\newglossaryentry{appforce}{
	name=\ensuremath{\bm{f}^{(\alpha)}},
	description={applied body force acting on $\Omega^{(\alpha)}$},
	type=physsymbols,
}
\newglossaryentry{virtualinterface}{
	name=\ensuremath{\psi},
	description={virtual interface separating the upper body from the lower body and primary variable of Problem~\ref{controlproblem}},
	type=physsymbols,
}
\newglossaryentry{3Dinterface}{
	name=\ensuremath{\Psi},
	description={three-dimensional analogue of $\psi$ that is constant along the transverse variable $x_3$},
	type=physsymbols,
}
\newglossaryentry{u1psi}{
	name=\ensuremath{\bm{u}^{(1)}(\Psi)},
	description={displacement field for the body $\overline{\Omega^{(1)}}$ that solves Problem~\ref{problem1}},
	type=physsymbols,
}
\newglossaryentry{u2psi}{
	name=\ensuremath{\bm{u}^{(2)}(\Psi)},
	description={displacement field for the body $\overline{\Omega^{(2)}}$ that solves Problem~\ref{problem2}},
	type=physsymbols,
}

\newtheorem{theorem}{Theorem}[section]
\newtheorem{lemma}{Lemma}[section]
\theoremstyle{definition}
\newtheorem{remark}[theorem]{Remark}

\newtheorem{innercustomgeneric}{\customgenericname}
\providecommand{\customgenericname}{}
\newcommand{\newcustomproblem}[2]{%
	\newenvironment{#1}[1]
	{%
		\renewcommand\customgenericname{#2}%
		\renewcommand\theinnercustomgeneric{##1}%
		\innercustomgeneric
	}
	{\endinnercustomgeneric}
}

\newcustomproblem{customprob}{Problem}

\newcommand*{\bqed}{\hfill\ensuremath{\blacksquare}}%

\newcommand{\vertiii}[1]{{\left\vert\kern-0.25ex\left\vert\kern-0.25ex\left\vert #1 
		\right\vert\kern-0.25ex\right\vert\kern-0.25ex\right\vert}}

\algblock{Input}{EndInput}

\algnotext{EndInput}

\algblock{Output}{EndOutput}

\algnotext{EndOutput}

\algblock{Iterate}{EndIterate}

\algnotext{EndIterate}

\allowdisplaybreaks[4]

\makeatletter
\newcommand{\addresseshere}{%
	\enddoc@text\let\enddoc@text\relax
}
\makeatother

\begin{document}
	
	
	\title[Optimal Interface Method]{The Optimal Interface Method: A Modelling Paradigm for Two-Body Contact Problems via Optimal Control and $\Gamma$-convergence}
	
	\author[Zhuopeng Liu]{Zhuopeng Liu}
	\address{Z.L. School of Science and Engineering, The Chinese University of Hong Kong (Shenzhen), 2001 Longxiang Blvd., Longgang District, Shenzhen, China}
	\email{zhuopengliu@link.cuhk.edu.cn}

	\author[Paolo Piersanti]{Paolo Piersanti}
	\address{P.P. School of Science and Engineering, The Chinese University of Hong Kong (Shenzhen), 2001 Longxiang Blvd., Longgang District, Shenzhen, China}
	\email[Corresponding author]{ppiersanti@cuhk.edu.cn}

\today

\begin{abstract}
This paper is concerned with the derivation of a model governing the equilibrium configuration of two linearly elastic bodies that may come into contact with each other. The contact zone, or coincidence set, is an \emph{a priori} unknown of the model, and the equilibrium configuration is governed by variational inequalities.

We introduce a new modelling paradigm in which the contact interface between the two linearly elastic bodies under consideration is now regarded as the \emph{primary variable} of an optimal control problem whose states are the unique solutions of two \emph{independent} unilateral contact problems for which this very contact interface represents the obstacle itself. As a result, the free boundary emerges as the solution of a low-dimensional optimisation problem rather than as the active set of coupled variational inequalities.

We establish the existence and uniqueness of optimal controls minimising the reduced cost functional arising from the modelling paradigm introduced beforehand, and we show that the family of marginal energy functionals associated with such model $\Gamma$-converges to the energy functional associated with the original variational inequalities as the positive \emph{cost parameter} $\beta$ in the reduced cost functional tends to zero. This $\Gamma$-convergence argument also yields a characterisation of the coincidence set of the limit problem in terms of the optimal interfaces.

We further show that, after a suitable regularisation, the reduced cost functional is G\^{a}teaux differentiable, and that its G\^{a}teaux derivative can be explicitly computed \emph{without solving any adjoint system}.

Finally, we corroborate the $\Gamma$-convergence result obtained beforehand with numerical experiments on a simplified and exactly solvable configuration.

\smallskip


\smallskip
\noindent \textbf{MSC~2020.} 49J40, 49J45, 49K21, 65K15, 74B05.
\end{abstract}

\maketitle
\tableofcontents

\section{Introduction}
The unilateral contact of a linearly elastic body with a rigid foundation was first formulated by Signorini~\cite{Sig33}. In its modern form, the problem consists in finding a displacement field satisfying the equilibrium equations together with the \emph{Signorini conditions} on the contact zone, i.e., non-penetration of the foundation and vanishing normal traction wherever contact is not active. Since these conditions are expressed as inequalities rather than equations, the problem admits no classical (namely, in terms of the Euler--Lagrange equations) formulation and is instead cast as \emph{variational inequalities} over the closed convex cone of kinematically admissible displacements. The existence and uniqueness of the solution for such contact problem were established by Fichera~\cite{Fic63,Fic64}, and the formulation was subsequently placed in the general framework of elliptic variational inequalities by Lions \& Stampacchia~\cite{DuvLio76,Lions1969}.

Obstacle problems in elasticity were recently studied by P.G. Ciarlet and his collaborators~\cite{CiaPie2018b,CiaPie2018bCR,CiaMarPie2018b,CiaMarPie2018,Pie22,Pie26-1,Pie26-2,Pie26-3}. The modelling and analysis of the deformation of linearly elastic virus capsids, modelled as a lattice, subject to lying over a flat rigid surface were carried out by means of the theory of obstacle problems by Temam and his group~\cite{HT17,PWDT2D,PWDT3D}.

The optimal control of variational inequalities is another interesting problem, and has long been studied in the mathematical literature. Adams, Lenhart and Yong studied an optimal control problem for \emph{scalar} variational inequalities with vanishing source term in~\cite{ALY98}. In this paper, the control was the obstacle determining the set of admissible solutions. In~\cite{AL03}, Adams and Lenhart extended the previous work to the case where the source term is non-vanishing, while maintaining the \emph{scalar} setting. In~\cite{AL02}, Adams and Lenhart studied an optimal control for \emph{scalar} parabolic obstacle problems where, again, the control is the obstacle identifying the set of admissible solutions. It is worth mentioning that in~\cite{AL02,CCT07} the authors chose to model the energy paid by the system to implement non-flat control by means of a biharmonic operator to be able to leverage compactness that is not, at least \emph{a priori}, available for the standard modelling by means of an energy of Dirichlet type employed in~\cite{AL03,ALY98,IK07}.

The problem we are here concerned with amounts to finding the equilibrium position of two linearly elastic bodies that, when deformed, may come to contact with each other. The central difficulties in the study of such problems lie not in determining the existence and uniqueness of displacement fields associated with the minimal energy configuration, but in the analysis of the \emph{coincidence set}, defined as the subset of the contact zone on which the two bodies actually touch, and in determining a suitable relaxation of the original variational inequalities.
The difficulty in the study of the coincidence set is due to the fact that it is \emph{a priori} unknown and, generally, not smooth. For what concerns the search for an approximation of the original variational inequalities modelling the contact problem under consideration, the classical approximation by penalty method fails to be a stepping stone for the numerical treatment of the problem, since the penalty term couples the two displacement fields through a four-dimensional fibre product.

\subsection*{Novelties and main results}
In this paper we propose a different point of view, which is the main \textbf{modelling contribution} of this work. Our approach consists in \emph{promoting the separating interface from an unknown to a primary variable}. More precisely, the separating interface is regarded as the \emph{control} variable of an optimal control problem, that we denote by Problem~\ref{controlproblem}, and the mechanical system is naturally driven towards choosing the interface that is least energetically expensive. Specifically, we look for a surface $\psi$ of $H^2_0$ regularity, identified with its three-dimensional extension $\Psi$ that is constant along the transverse variable, and that separates the two deformed bodies. We then minimise the total elastic energy of the two bodies \emph{augmented by a cost term} penalising the curvature of the interface. The cost parameter, that we denote by $\beta>0$, weights the interface curvature against the elastic energy. For large values of $\beta$, the model drives the mechanical system towards selecting interfaces that are as flat as possible, while for small values of $\beta$ the main energetic contribution is due to the elastic energy of the two linearly elastic bodies under consideration, and the interface is thus free to adapt to the contact geometry at a negligible cost.

For a prescribed interface, we are then able to obtain two independent single-body unilateral problems in which the interface plays the role of the obstacle. The non-interpenetration constraint is enforced by the interface itself, thus lowering the complexity of the problem under consideration. As a result, the free boundary is no longer an unknown active set to be tracked, but the solution of a less complex optimisation problem posed over the space of admissible interfaces.
We now discuss the two new results presented here:
\begin{itemize}
	\item[$(1)$]Our \textbf{first main result} is a characterisation of the coincidence set of Problem~\ref{problem0}, the variational problem formulated in terms of variational inequalities governing the deformation of two linearly elastic bodies that are not allowed to interpenetrate. We obtain this characterisation by approximating Problem~\ref{problem0} via an optimal control problem, Problem~\ref{controlproblem}, in which the admissible controls are \emph{virtual interfaces} separating the two bodies and acting as obstacles. The connection between the coincidence set of the \emph{original} obstacle problem and the optimal control formulation is investigated via a $\Gamma$-convergence argument which, to our knowledge, has not been previously resorted to in the literature. More precisely, we show that the optimal states for Problem~\ref{controlproblem} converge to the solution of Problem~\ref{problem0} in $L^2$ as $\beta\to 0$ (Theorem~\ref{gamma-4}).
	\item[$(2)$] Our \textbf{second main result} amounts to proposing a formulation for the regularised cost functional for which its G\^{a}teaux derivative is obtained \emph{without} solving any adjoint system (Theorem~\ref{th:Gateaux:Jrhoc}). More precisely, we will show that an application of the envelope theorem (cf., e.g., \cite{MS02} and Section~6.1 in~\cite{Car01}) shows that the G\^{a}teaux derivative of the regularised reduced cost functional reduces to a direct post-processing of the state, rather than requiring the solution of an adjoint problem. As a result, the algorithm we propose for implementing the method presented here does not require the implementation of the loop associated with the adjoint problem, as it is instead customary for other models~\cite{AL03,ALY98,AHL09,CCT07,IK07}.
\end{itemize}

It is worth mentioning that the scalar nature of the models treated in the literature discussed beforehand~\cite{AL02,AL03,ALY98,AHL09,CCT07,IK07} allows the integration domain of the optimal control problem to be the same as the integration domain of the state problems, and leads to significant simplifications in the mathematical analysis. These simplifications are lost in our formulation, as the state problems are integrated over two separate three-dimensional sets and the confinement condition is defined along vertical fibres.

\textbf{The main goal of this paper is the development of a modelling paradigm for two-body unilateral contact in which the separating interface is promoted from an unknown to the primary variable of an optimal control problem, so that the free boundary is recovered as the solution of a low-dimensional energy minimisation. This paradigm is well-posed, is consistent with the classical formulation based on variational inequalities, and gives rise to a gradient descent algorithm that does not require solving any adjoint system}.

\subsection*{Structure of the paper}
The paper is divided into seven sections, including this one. In Section~\ref{Sec:2} we present the notation and the geometrical framework. In Section~\ref{Sec:3} we formulate the contact problem between two linearly elastic bodies as a set of variational inequalities posed over a non-empty, closed and convex subset of the product space of admissible displacements, and we prove that this problem admits a unique solution (Theorem~\ref{th-ex-un-0}). In Section~\ref{Sec:3-1} we introduce the optimal control approximation, in which a virtual interface acts both as control and as obstacle, we prove existence and uniqueness of the optimal control and of the optimal states (Theorem~\ref{ex-un-op-con}), and we show that the family of marginal energy functionals $\Gamma$-converges to the energy functional associated with the variational problem presented in Section~\ref{Sec:3} (Theorem~\ref{gamma-4}), thereby justifying the optimal control formulation. In Section~\ref{Sec:3bis} we penalise the optimal control by means of a non-differentiable term and prove convergence of the penalised solutions and controls (Theorem~\ref{th:optimal-1}). In Section~\ref{Sec:4} we consider a Moreau--Yosida regularisation of this penalty term, show that the control-to-state mapping is G\^{a}teaux differentiable, and that the reduced cost functional is G\^{a}teaux differentiable as well. This differentiability is used to implement a gradient descent method, that does not require solving any adjoint system, in the numerical experiments reported in Section~\ref{Sec:5}. A list of symbols is appended to facilitate the reading.

\section{Geometrical and Analytical Preliminaries}\label{Sec:2}
For a review of the standard concepts in differential geometry utilized here, one may consult references such as~\cite{Ciarlet2000} or~\cite{Ciarlet2005}.
Unless employed for sequence indexing, Latin indices are constrained to the set $\{1,2,3\}$, Greek indices are constrained to the set $\{1,2\}$, and the Einstein summation convention is systematically applied in conjunction with these two rules.

Given an open set $\gls{Omega}\subset \mathbb{R}^3$, the symbols $L^2(\Omega)$ and $H^1(\Omega)$ represent the conventional Lebesgue and Sobolev spaces. The space $\gls{DOmega}$ consists of all infinitely differentiable functions defined on $\Omega$ whose support is compact in $\Omega$. The norm in a vector space $X$ is denoted by $\|\cdot\|_X$. Bold typeface, as in $\bm{L}^2(\Omega)$ or $\bm{H}^1(\Omega)$, is used to designate spaces consisting of vector-valued functions. Given a measurable function $f:\Omega\to\mathbb{R}$, we denote by $\gls{f+}$ its positive part, whose definition is given by (cf., e.g., \cite{ABM14,EG15}) $\{f\}^{+}(x):=\max\{f(x),0\}$ for a.a. $x\in\Omega$, and we denote by $\{f\}^{-}$ its negative part, whose definition is given by (cf., e.g., \cite{ABM14,EG15}) $\{f\}^-(x):=-\min\{f(x),0\}$ for a.a. $x\in\Omega$.

A \emph{Lipschitz domain} $\Omega\subset\mathbb{R}^3$ is characterised as a non-empty, bounded, connected, open subset whose boundary $\gls{Gamma}$ is Lipschitz continuous, with $\Omega$ situated locally on a single side of $\Gamma$. The volume element in $\Omega$ is written as $\mathrm{d}x$, while $\gls{dGamma}$ denotes the area element on $\Gamma$. A more comprehensive treatment is available in sources like Section~8.2 of~\cite{Cia25}. Furthermore, let $\Gamma=\Gamma_0\cup\Gamma_1$ be a subset of the boundary that is measurable with respect to $\mathrm{d}\Gamma$, satisfying $\Gamma_0\cap\Gamma_1 =\emptyset$ and $\textup{area }\Gamma_0>0$. We denote by $\bm{\nu}$ the outer unit normal vector along $\Gamma$. We denote by $\overline{\Omega}$ the closure of $\Omega$ with respect to the Euclidean metric of the ambiance space. The special notation $\Omega_1\gls{subsub}\Omega$ means that $\overline{\Omega_1}\subset\Omega$ and $\textup{dist}(\partial\Omega,\partial\Omega_1):=\min\{|x-y|;x\in\partial\Omega\textup{ and } y\in\partial\Omega_1\}>0$.

The three-dimensional ``physical space'' $\mathbb{R}^3$ is modelled as a \emph{real three-dimensional affine Euclidean space}. This is constructed by choosing an origin $O$ and associating it with a \emph{real three-dimensional Euclidean space} $\gls{E3}$. An \emph{orthonormal basis} of vectors $\gls{ei}$ is defined for $\mathbb{E}^3$. The Euclidean inner product of $\bm{a}$ and $\bm{b}$ in $\mathbb{E}^3$ is expressed as $\bm{a} \cdot \bm{b}$; the Euclidean norm of an element $\bm{a} \in \mathbb{E}^3$ is $|\bm{a}|$; the Kronecker symbol is $\delta^{ij}$.

Without any loss of rigour or generality, we systematically identify $\mathbb{R}^3$ with the affine Euclidean space $\mathbb{E}^3$, meaning that every point $x\in\mathbb{R}^3$ is associated with a unique geometrical vector $\bm{Ox}\in\mathbb{E}^3$. The combination of the origin $O\in\mathbb{R}^3$ and the orthonormal basis vectors $\bm{e}^i\in\mathbb{E}^3$ constitutes a \emph{Cartesian frame} for $\mathbb{R}^3$. The components $x_i$ of the vector $\bm{Ox}$ with respect to the basis $\bm{e}^i$ are the \emph{Cartesian coordinates} of the point $x$, or the \emph{Cartesian components} of the vector $\bm{Ox}$. After fixing a Cartesian frame, any point $x$ is \emph{identified} with the vector $\bm{Ox}=x_i \bm{e}^i \in \mathbb{E}^3$. The partial derivative $\partial/\partial x_i$ is abbreviated as $\partial_i$.

The closure $\overline{\Omega}\subset\mathbb{E}^3$ symbolises the \emph{reference configuration} of a \emph{linearly elastic body} free from external loads. In line with \cite{Ciarlet1988}, the material is taken to be \emph{isotropic}, \emph{homogeneous}, and \emph{linearly elastic}. The mechanical response of this material is thus determined by its two \emph{Lam\'{e} constants} $\gls{lambda}\ge 0$ and $\gls{mu}>0$.

Herein, the abbreviation ``a.e.'' means ``almost everywhere'' and ``a.a.'' means ``almost all''. 
Denote by $\gls{IOmega}$ the vector field $\mathbb{I}_\Omega:\overline{\Omega}\to\mathbb{E}^3$ defined by $\mathbb{I}_\Omega(x):=x$ for all $x\in\overline{\Omega}$. Deformations of the reference configuration $\overline{\Omega}$ are modelled via displacement vector fields $\bm{v}:\overline{\Omega}\to\mathbb{R}^3$. The deformation of the linearly elastic body $\overline{\Omega}$ associated with the displacement $\bm{v}$ takes the form $\mathbb{I}_\Omega(x)+v_i(x)\bm{e}^i$ for all or, possibly, a.a. $x\in\overline{\Omega}$.
The function space of admissible displacements for the linearly elastic models under consideration is defined as
\begin{equation*}
	\gls{VOmega}:=\{\bm{v}=(v_i)\in\bm{H}^1(\Omega);\bm{v}={\bf{0}}\,\mathrm{d}\Gamma\textup{-a.e. on }\Gamma_0\},
\end{equation*}
and is equipped with the norm
\begin{equation*}
	\|\bm{v}\|_{\bm{V}(\Omega)}:=\left(\sum_{i}\|v_i\|_{H^1(\Omega)}^2\right)^{1/2}.
\end{equation*}

The fourth order three-dimensional elasticity tensor is expressed in Cartesian coordinates, with its components written as \gls{A}.
Its contravariant components are given by (see, e.g., \cite{Ciarlet1988})
\begin{equation*}
	A^{ijk\ell}:=\lambda\delta^{ij}\delta^{k\ell}+\mu(\delta^{ik}\delta^{j\ell}+\delta^{i\ell}\delta^{jk}),
\end{equation*}
and it holds that $A^{ijk\ell}=A^{jik\ell}=A^{k\ell ij}$.

For a vector field $\bm{v}\in\bm{H}^1(\Omega)$, the \emph{linearised change of metric tensor} $\gls{e}$ is introduced, with components $e_{i\|j}(\bm{v})$ defined by:
\begin{equation*}
	e_{i\|j}(\bm{v}):=\dfrac{1}{2}(\partial_j v_i +\partial_i v_j) \in L^2(\Omega).
\end{equation*}

This tensor is symmetric, meaning that $e_{i\|j}(\bm{v})=e_{j\|i}(\bm{v})$ for every $\bm{v}\in\bm{H}^1(\Omega)$.
We now present \emph{Korn's inequality} formulated in Cartesian coordinates (refer to, for example, Theorem~6.3-4 in~\cite{Ciarlet1988}).
\begin{theorem}
	\label{KornCart}
	Consider a Lipschitz domain $\Omega$ in $\mathbb{R}^3$ where $\Gamma_0$ is a subset of the boundary $\Gamma$ with positive area.
	Then, there exists a constant $c_0>0$ for which
	$$
	c_0^{-1}\|\bm{v}\|_{\bm{H}^1(\Omega)}\le\|\bm{e}(\bm{v})\|_{\mathbb{L}^2(\Omega)}\le c_0\|\bm{v}\|_{\bm{H}^1(\Omega)},
	$$
	for all $\bm{v}\in\bm{V}(\Omega)$.
	\qed
\end{theorem}

This profound inequality has been established through multiple proofs. Significant contributions include~\cite{Friedrichs1947,Gobert1962,HlaNec1970a,HlaNec1970b}, page~110 of \cite{DuvLio76}, and Section~6.3 of \cite{NecHla1981}. A proof within the space $W^{1,p}(\Omega)$ for $1 < p \le \infty$ is provided in~\cite{Temam1983}. A straightforward demonstration is available in \cite{Nit1981}, with supplementary information found in Appendix (A) of \cite{Miyoshi1985}.

Let $\gls{A2}:\bm{V}(\Omega) \to \bm{V}^\ast(\Omega)$ be the linear and continuous operator defined by the relation:
\begin{equation}
	\label{A}
	\langle\mathcal{A}\bm{w},\bm{v}\rangle_{\bm{V}^\ast(\Omega),\bm{V}(\Omega)}:=\int_{\Omega}A^{ijk\ell} e_{k\|\ell}(\bm{w}) e_{i\|j}(\bm{v})\,\mathrm{d}x,\quad\textup{ for all }\bm{w}, \bm{v}\in\bm{V}(\Omega).
\end{equation}

Owing to the characteristics of the elasticity tensor $\{A^{ijk\ell}\}$ and Korn's inequality in Cartesian coordinates (Theorem~\ref{KornCart}), the operator $\mathcal{A}$ introduced in~\eqref{A} is symmetric and uniformly positive-definite. This means a constant $\gls{Ce}>0$ exists such that
\begin{equation}
	\label{eq:pos-def}
	\sum_{i,j=1}^{3}|t_{ij}|^2 \le C_e A^{ijk\ell}(x) t_{k\ell}t_{ij},
\end{equation}
for every point $x\in\overline{\Omega}$ and every symmetric matrix $(t_{ij})$.

\section{An Optimal Interface Method for Modelling the Contact between two Linearly Elastic Bodies}\label{Sec:3}

The purpose of this section is to propose a mathematical model that describes the deformation of two linearly elastic bodies, whose states at rest are respectively denoted by $\gls{Omega1}$ and $\gls{Omega2}$. We assume that $\Omega^{(\alpha)}$ is a Lipschitz domain in $\mathbb{R}^3$. To fix ideas, assume without loss of generality that the linearly elastic body $\overline{\Omega^{(1)}}$ lies above the linearly elastic body $\overline{\Omega^{(2)}}$. For the sake of physical admissibility, we require that
\begin{equation}
	\label{eq:psi}
	\max_{\overline{\Omega^{(2)}}}(\mathbb{I}_{\Omega^{(2)}}\cdot\bm{e}_3)<0<\min_{\overline{\Omega^{(1)}}}(\mathbb{I}_{\Omega^{(1)}}\cdot\bm{e}_3),
\end{equation}
which means that the two bodies do not make contact if no forces are applied. Since the linearly elastic bodies $\overline{\Omega^{(\alpha)}}$ are compact in $\mathbb{E}^3$ and since the admissible deformations for these bodies are \emph{linearly elastic}, it is licit to \emph{assume} that there exists a compact set $\gls{K}\subset\mathbb{R}^3$ of the form $K=\overline{\kappa}\times[-a,a]$, where $\kappa$ is a Lipschitz domain $\mathbb{R}^2$ and $a>0$, containing all the possible admissible deformations of the linearly elastic bodies $\overline{\Omega^{(\alpha)}}$. We denote by $\Gamma_0^{(\alpha)}$ the $\mathrm{d}\mathcal{H}^2$-measurable non-zero-area portion of $\partial\Omega^{(\alpha)}$ where homogeneous Dirichlet boundary conditions of place are imposed.

Let $\gls{pikappa}:K\to\overline{\kappa}$ denote the \emph{tangential projection} $\pi_\kappa(x_1,x_2,x_3):=(x_1,x_2)$, and define the overlap set $\gls{omegaproj}$ by:
\begin{equation*}
	\omega:=\pi_\kappa(\Omega^{(1)})\cap\pi_\kappa(\Omega^{(2)}).
\end{equation*}

Assume that $\Omega^{(1)}$ and $\Omega^{(2)}$ are such that:
\begin{itemize}
	\item[$(i)$] $\omega$ is bounded and measurable;
	\item[$(ii)$] $\mathcal{L}^2(\omega)>0$.
\end{itemize}

Assumption $(ii)$ guarantees that all $\mathrm{d}\mathcal{L}^2$-a.e. statements on $\omega$ and all the integrals over $\omega$ in the sequel are well defined while assumption $(ii)$ makes the contact region non-degenerate.

For each tangential coordinate $\xi=(\xi_1,\xi_2)\in\pi_\kappa(\Omega^{(\alpha)})$ and each $\alpha\in\{1,2\}$, the \emph{open vertical fibre} of body $\overline{\Omega^{(\alpha)}}$ corresponding to $\xi$ is given by
\begin{equation*}
	\gls{Falpha}:=\left\{x_3\in\mathbb{R}; (\xi_1,\xi_2,x_3)\in\Omega^{(\alpha)}\right\},
\end{equation*}
and we extend $F^{(\alpha)}$ by $F^{(\alpha)}(\xi)=\emptyset$ for $\xi\in\kappa\setminus\pi_\kappa(\Omega^{(\alpha)})$. We assume that the two bodies are \emph{vertical graphs over their projections}, namely, that for $\alpha\in\{1,2\}$ there exist continuous functions $a^{(\alpha)}:\overline{\pi_\kappa(\Omega^{(\alpha)})}\to\mathbb{R}$ and $b^{(\alpha)}:\overline{\pi_\kappa(\Omega^{(\alpha)})}\to\mathbb{R}$ such that $a^{(\alpha)}\le b^{(\alpha)}$, and such that:
\begin{equation}
	\label{eq:graph}
	\begin{aligned}
		\gls{Omega1}&=\left\{(\xi,t); \xi\in\overline{\pi_\kappa(\Omega^{(1)})} \textup{ and } a^{(1)}(\xi)\le t\le b^{(1)}(\xi)\right\},\\
		\gls{Omega2}&=\left\{(\xi,t); \xi\in\overline{\pi_\kappa(\Omega^{(2)})} \textup{ and } a^{(2)}(\xi)\le t\le b^{(2)}(\xi)\right\}.
	\end{aligned}
\end{equation}

Consequently, it results that $\overline{F^{(\alpha)}(\xi)}$ is a non-empty compact interval, which possibly reduces to a single point, and the mapping $\xi\in\overline{\pi_\kappa(\Omega^{(\alpha)})}\mapsto \overline{F^{(\alpha)}(\xi)}$ has non-empty compact values and is continuous in the Hausdorff metric (cf., e.g., \cite{Manetti2023}).
Define the set
\begin{equation}
	\label{eq:U-step2}
	\begin{aligned}
		\gls{U}:=\Big\{&(\bm{v}^{(1)},\bm{v}^{(2)})\in \bm{V}(\Omega^{(1)})\times\bm{V}(\Omega^{(2)});\\
		&\quad\textup{ for } \mathrm{d}\mathcal{L}^2\textup{-a.e. } \xi=(\xi_1,\xi_2)\in\omega,\textup{ for } \mathrm{d}\mathcal{L}^1\textup{-a.e. } x_3^{(1)}\in F^{(1)}(\xi)\textup{ and } \mathrm{d}\mathcal{L}^1\textup{-a.e. } x_3^{(2)}\in F^{(2)}(\xi),\\
		&\quad (\mathbb{I}_{\Omega^{(1)}}(\xi_1,\xi_2,x_3^{(1)})+\bm{v}^{(1)}(\xi_1,\xi_2,x_3^{(1)}))\cdot\bm{e}_3
		\ge(\mathbb{I}_{\Omega^{(2)}}(\xi_1,\xi_2,x_3^{(2)})+\bm{v}^{(2)}(\xi_1,\xi_2,x_3^{(2)}))\cdot\bm{e}_3
		\Big\},
	\end{aligned}
\end{equation}
and note that a convenient reformulation of the constraint in~\eqref{eq:U-step2} can be given by defining the \emph{interface set}
\begin{equation*}
	\gls{S}:=\left\{((\xi_1,\xi_2),t_1,t_2)\in\omega\times\mathbb{R}\times\mathbb{R}; t_1\in F^{(1)}(\xi_1,\xi_2) \textup{ and } t_2\in F^{(2)}(\xi_1,\xi_2)\right\},
\end{equation*}
equipped with the product measure $\mu:=\mathcal{L}^2\otimes\mathcal{L}^1\otimes\mathcal{L}^1$ restricted to $S$. Define the \emph{constraint qualification mapping} \gls{Phi} as follows:
\begin{equation}
	\label{eq:Phidef}
	\begin{aligned}
		\bm{\Phi}(\bm{v}^{(1)},\bm{v}^{(2)})(\xi_1,\xi_2,t_1,t_2):=&
		\left(\mathbb{I}_{\Omega^{(1)}}(\xi_1,\xi_2,t_1)+\bm{v}^{(1)}(\xi_1,\xi_2,t_1)\right)\cdot\bm{e}_3\\
		&-\left(\mathbb{I}_{\Omega^{(2)}}(\xi_1,\xi_2,t_2)+\bm{v}^{(2)}(\xi_1,\xi_2,t_2)\right)\cdot\bm{e}_3 .
	\end{aligned}
\end{equation}

Let $p_1(\xi_1,\xi_2,t_1,t_2):=(\xi_1,\xi_2,t_1)$ and $p_2(\xi_1,\xi_2,t_1,t_2):=(\xi_1,\xi_2,t_2)$ denote the coordinate projections of $S$ onto $\overline{\Omega^{(1)}}$ and $\overline{\Omega^{(2)}}$, respectively.
Note that $\Phi$ is affine and continuous as a map from $\bm{V}(\Omega^{(1)})\times\bm{V}(\Omega^{(2)})$ into $L^1(S)$. Indeed, since every fibre lies in $[-a,a]$, we have $|t_1-t_2|\le 2a$ on $S$, and Fubini's theorem (cf., e.g., Theorem~1.15-6$(b)$ in~\cite{Cia25}) gives:
\begin{equation*}
	\begin{aligned}
		\|v^{(1)}_3\circ p_1\|_{L^1(S)}=\int_{\omega}|F^{(2)}(\xi)|\int_{F^{(1)}(\xi)}|v^{(1)}_3(\xi,t_1)|\,\mathrm{d}t_1\,\mathrm{d}\xi\le 2a\|v^{(1)}_3\|_{L^1(\Omega^{(1)})},\\
		\|v^{(2)}_3\circ p_2\|_{L^1(S)}=\int_{\omega}|F^{(1)}(\xi)|\int_{F^{(2)}(\xi)}|v^{(2)}_3(\xi,t_2)|\,\mathrm{d}t_2\,\mathrm{d}\xi\le 2a\|v^{(2)}_3\|_{L^1(\Omega^{(2)})}.
	\end{aligned}
\end{equation*}

The fact that $\omega\subset\kappa$ is bounded gives $\mu(S)\le(2a)^2\mathcal{L}^2(\kappa)<\infty$ so that it is possible to infer:
\begin{equation*}
\begin{aligned}
	&\|\Phi(\bm{v}^{(1)},\bm{v}^{(2)})\|_{L^1(S)}\le 2a\mu(S)+\|v^{(1)}_3\circ p_1\|_{L^1(S)}+\|v^{(2)}_3\circ p_2\|_{L^1(S)}\\
	&\le (2a)^3\mathcal{L}^2(\kappa)+2a\|v^{(1)}_3\|_{L^1(\Omega^{(1)})}+2a\|v^{(2)}_3\|_{L^1(\Omega^{(2)})}.
\end{aligned}
\end{equation*}

By Fubini's theorem (cf., e.g., Theorem~1.15-6$(b)$ in~\cite{Cia25}), we can re-cast the constraints in~\eqref{eq:U-step2} in terms of the following single condition:
\begin{equation}
	\label{eq:Phi}
	\Phi(\bm{v}^{(1)},\bm{v}^{(2)})\ge 0\quad\mathrm{d}\mu\textup{-a.e. in }S.
\end{equation}

We now establish that $\vec{\bm{U}}$ is non-empty, closed and convex.
\begin{lemma}
	\label{lem:U-step2}
	Let $\vec{\bm{U}}$ be defined by \eqref{eq:U-step2}, and assume that $\omega\neq\emptyset$ and that $\overline{\Omega^{(\alpha)}}\subset\subset K=\overline{\kappa}\times[-a,a]$ for $\alpha\in\{1,2\}$. Then $\vec{\bm{U}}$ is non-empty, convex and closed in $\bm{V}(\Omega^{(1)})\times\bm{V}(\Omega^{(2)})$.
\end{lemma}
\begin{proof}
	The proof is divided into three steps, numbered $(i)$--$(iii)$.
	
	$(i)$\emph{The set $\vec{\bm{U}}$ is non-empty.} Since both $\overline{\Omega^{(1)}}$ and $\overline{\Omega^{(2)}}$ are compact, and the gap function
	$x\in\overline{\Omega^{(\alpha)}}\mapsto\left(\mathbb{I}_{\Omega^{(\alpha)}}(x)\cdot\bm{e}_3\right)$ is continuous, so the maximum and the minimum for one such gap function in \eqref{eq:psi} are attained. Moreover, there exists a strictly positive gap between the compact sets $\overline{\Omega^{(1)}}$ and $\overline{\Omega^{(2)}}$, that we denote and define by:
	\begin{equation*}
		\delta:=\min_{\overline{\Omega^{(1)}}}(\mathbb{I}_{\Omega^{(1)}}\cdot\bm{e}_3)-\max_{\overline{\Omega^{(2)}}}(\mathbb{I}_{\Omega^{(2)}}\cdot\bm{e}_3)>0.
	\end{equation*}
	
	Take $\bm{v}^{(1)}=\bm{0}$ and $\bm{v}^{(2)}=\bm{0}$. For every $\xi\in\omega$, every $x_3^{(1)}\in F^{(1)}(\xi)$ and every $x_3^{(2)}\in F^{(2)}(\xi)$, it results $(\xi_1,\xi_2,x_3^{(\alpha)})\in\Omega^{(\alpha)}\subset\overline{\Omega^{(\alpha)}}$, whence:
	\begin{equation*}
	\begin{aligned}
		&\left(\mathbb{I}_{\Omega^{(1)}}(\xi_1,\xi_2,x_3^{(1)})+\bm{v}^{(1)}(\xi_1,\xi_2,x_3^{(1)})\right)\cdot\bm{e}_3-\left(\mathbb{I}_{\Omega^{(2)}}(\xi_1,\xi_2,x_3^{(2)})+\bm{v}^{(2)}(\xi_1,\xi_2,x_3^{(2)})\right)\cdot\bm{e}_3\\
		&\qquad=\left(\mathbb{I}_{\Omega^{(1)}}(\xi_1,\xi_2,x_3^{(1)})-\mathbb{I}_{\Omega^{(2)}}(\xi_1,\xi_2,x_3^{(2)})\right)\cdot\bm{e}_3\\
		&\qquad\ge \min_{\overline{\Omega^{(1)}}}(\mathbb{I}_{\Omega^{(1)}}\cdot\bm{e}_3)-\max_{\overline{\Omega^{(2)}}}(\mathbb{I}_{\Omega^{(2)}}\cdot\bm{e}_3)=\delta>0.
	\end{aligned}
	\end{equation*}
	
	Thus the inequality defining $\vec{\bm{U}}$ in \eqref{eq:U-step2} holds \emph{strictly}, with the uniform gap $\delta$, for \emph{every} admissible point $(\xi,x_3^{(1)},x_3^{(2)})$. In particular, it holds for $\mathrm{d}\mathcal{L}^1$-a.e. $x_3^{(1)}\in F^{(1)}(\xi)$, $\mathrm{d}\mathcal{L}^1$-a.e. $x_3^{(2)}\in F^{(2)}(\xi)$ and $\mathrm{d}\mathcal{L}^2$-a.e. $\xi\in\omega$. This shows that $(\bm{0},\bm{0})\in\vec{\bm{U}}$ and, thus, that the set $\vec{\bm{U}}$ is non-empty, as it was to be proved.
	
	$(ii)$ \emph{The set $\vec{\bm{U}}$ is closed in $\bm{V}(\Omega^{(1)})\times\bm{V}(\Omega^{(2)})$.} Let $(\bm{v}^{(1)}_n,\bm{v}^{(2)}_n)\to(\bm{v}^{(1)},\bm{v}^{(2)})$ strongly in $\bm{V}(\Omega^{(1)})\times\bm{V}(\Omega^{(2)})$ as $n\to\infty$. By the continuity of $\Phi$ defined in~\eqref{eq:Phidef}, it results that the functions $\Phi_n:=\Phi(\bm{v}^{(1)}_n,\bm{v}^{(2)}_n)$ converge to $\Phi:=\Phi(\bm{v}^{(1)},\bm{v}^{(2)})$ in $L^1(S)$ as $n\to\infty$. Since $(\bm{v}^{(1)}_n,\bm{v}^{(2)}_n)\in\vec{\bm{U}}$ for each integer $n\ge 1$, it results that $\Phi_n\ge 0$ $\mathrm{d}\mu$-a.e. in $S$ for each integer $n\ge 1$. Since the cone $C:=\{f\in L^1(S); f\ge 0\,\mathrm{d}\mu\textup{-a.e. in }S\}$ is closed in $L^1(S)$, we obtain that $\Phi\ge 0$ $\mathrm{d}\mu$-a.e. in $S$, that is equivalent to stating that $(\bm{v}^{(1)},\bm{v}^{(2)})\in\vec{\bm{U}}$.
	
	$(iii)$ \emph{The set $\vec{\bm{U}}$ is convex.} In order to establish the sought property, it suffices to observe that, since $C$ is convex and $\Phi$ is affine, the set $\vec{\bm{U}}=\Phi^{-1}(C)$ is convex. This completes the proof.
\end{proof}

For the sake of simplicity, we do not consider the action of applied surface forces, even though it is licit to consider the action of applied surface forces on the \emph{a priori} unknown portion of boundary of the deformation of $\overline{\Omega^{(1)}}$ that does not engage contact with the obstacle and, \emph{a fortiori} with the other linearly elastic body $\overline{\Omega^{(2)}}$. We denote by $\gls{appforce}\in\bm{L}^2(\Omega^{(\alpha)})$ the applied body force acting on the linearly elastic body $\overline{\Omega^{(\alpha)}}$.
An intuitive attempt in formulating the variational problem governing the deformation of two linearly elastic bodies that could make contact in the fashion described by the confinement condition characterising the set $\vec{\bm{U}}$ is the following.

\begin{customprob}{$\mathcal{P}(\Omega^{(1)};\Omega^{(2)})$}\label{problem0}
	Find $(\bm{u}^{(1)},\bm{u}^{(2)})\in\vec{\bm{U}}$ satisfying the variational inequalities
	\begin{equation*}
		\begin{aligned}
			&\int_{\Omega^{(1)}} A^{ijk\ell,1} e_{k\|\ell}(\bm{u}^{(1)}) e_{i\|j}(\bm{v}_1-\bm{u}^{(1)}) \,\mathrm{d}x^{(1)}+\int_{\Omega^{(2)}} A^{ijk\ell,2} e_{k\|\ell}(\bm{u}^{(2)}) e_{i\|j}(\bm{v}_2-\bm{u}^{(2)}) \,\mathrm{d}x^{(2)}\\
			&\ge \int_{\Omega^{(1)}} \bm{f}^{(1)} \cdot (\bm{v}_1-\bm{u}^{(1)}) \,\mathrm{d}x^{(1)}+\int_{\Omega^{(2)}} \bm{f}^{(2)} \cdot (\bm{v}_2-\bm{u}^{(2)}) \,\mathrm{d}x^{(2)},
		\end{aligned}
	\end{equation*}
	for all $(\bm{v}_1,\bm{v}_2)\in\vec{\bm{U}}$.
	\bqed	
\end{customprob}

The following existence and uniqueness result is straightforward to establish.
\begin{theorem}
	\label{th-ex-un-0}
	Problem~\ref{problem0} admits a unique solution $(\bm{u}^{(1)},\bm{u}^{(2)})\in\vec{\bm{U}}$.
\end{theorem}
\begin{proof}
	Equip $\bm{V}(\Omega^{(1)})\times\bm{V}(\Omega^{(2)})$ with the norm $\|\cdot\|_{1,2}$ defined by:
	\begin{equation*}
		\|(\bm{v}_1,\bm{v}_2)\|_{1,2}:=\left(\|\bm{v}_1\|_{\bm{V}(\Omega^{(1)})}^2+\|\bm{v}_2\|_{\bm{V}(\Omega^{(2)})}^2\right)^{1/2},\quad\textup{ for all }(\bm{v}_1,\bm{v}_2)\in\bm{V}(\Omega^{(1)})\times\bm{V}(\Omega^{(2)}).
	\end{equation*}
	
	Thanks to the symmetry of the fourth order three-dimensional elasticity tensors $\{A^{ijk\ell,\alpha}\}$ and the symmetry of the linearised strain tensor, the bilinear form
	\begin{equation*}
		(\bm{v}_1,\bm{v}_2)\in\bm{V}(\Omega^{(1)})\times\bm{V}(\Omega^{(2)})\mapsto\sum_{\alpha=1}^{2}\left\{\int_{\Omega^{(\alpha)}}A^{ijk\ell,\alpha}e_{k\|\ell}(\bm{v}_\alpha) e_{i\|j}(\bm{v}_\alpha) \,\mathrm{d}x^{(\alpha)}\right\}
	\end{equation*}
	is symmetric. Boundedness for this bilinear form follows from H\"older's inequality (cf., e.g., \cite{Brez11}).
	Besides, the uniform positive definiteness of the fourth order three-dimensional elasticity tensors $\{A^{ijk\ell,\alpha}\}$ (cf., \eqref{eq:pos-def}) and Korn's inequality in Cartesian coordinates (Theorem~\ref{KornCart}) give:
	\begin{equation*}
		\sum_{\alpha=1}^{2}\left\{\int_{\Omega^{(\alpha)}}A^{ijk\ell,\alpha}e_{k\|\ell}(\bm{v}_\alpha) e_{i\|j}(\bm{v}_\alpha) \,\mathrm{d}x^{(\alpha)}\right\}
		\ge\sum_{\alpha=1}^{2}\left\{\dfrac{\|\bm{v}_\alpha\|_{\bm{V}(\Omega^{(\alpha)})}^2}{\left(c_0^{(\alpha)}\right)^2 C_e^{(\alpha)}}\right\},\quad\textup{ for all }(\bm{v}_1,\bm{v}_2)\in\bm{V}(\Omega^{(1)})\times\bm{V}(\Omega^{(2)}).
	\end{equation*}
	
	The previous estimate in turn implies that the bilinear form above is $(\bm{V}(\Omega^{(1)})\times\bm{V}(\Omega^{(2)}))$-elliptic.
	Finally, the linear form
	\begin{equation*}
		(\bm{v}_1,\bm{v}_2)\in\bm{V}(\Omega^{(1)})\times\bm{V}(\Omega^{(2)})\mapsto\sum_{\alpha=1}^{2}\int_{\Omega^{(\alpha)}}\bm{f}^{(\alpha)}\cdot\bm{v}_\alpha\,\mathrm{d}x^{(\alpha)}
	\end{equation*}
	is continuous. Since the set $\vec{\bm{U}}$ is non-empty, closed and convex (Lemma~\ref{lem:U-step2}), an application the Lions--Stampacchia theorem (Theorem~4.8-1 in~\cite{Cia25}) allows us to complete the proof.
\end{proof}

While the deformation of such linearly elastic bodies can be clearly described by means of the classical three-dimensional linearly elastic energy (cf., e.g., \cite{Ciarlet1988}), what is less clear is how to track the contact area which, at least \emph{a priori}, is a hypersurface of codimension one. The purpose of the next section is to provide a modelling methodology that allows to keep track of the contact area while gaining insight into the way the two linearly elastic bodies under consideration deform.

\section{An approximation of Problem~\ref{problem0} via optimal control theory}\label{Sec:3-1}
We now propose a method to approximate the free boundary of Problem~\ref{problem0} by the graph of a \emph{sufficiently smooth} function $\gls{virtualinterface}:\mathbb{R}^2\to\mathbb{R}$ such that:
\begin{itemize}
	\item[$(i)$] The function $\psi$ is, up to a change in the coordinate system, \emph{as close as possible to being zero};
	\item[$(ii)$] The function $\psi$ \emph{separates} the two deformed elastic bodies, in the sense that deformations of $\overline{\Omega^{(1)}}$ lie in the epigraph of $\psi$, denoted by $\textup{epi}\psi:=\{(x_1,x_2,x_3)\in\mathbb{R}^3; x_3\ge\psi(x_1,x_2)\}$, and deformations of $\overline{\Omega^{(2)}}$ lie in the closed hypograph of $\psi$, that is, in $\{(x_1,x_2,x_3)\in\mathbb{R}^3; x_3\le\psi(x_1,x_2)\}$;
	\item[$(iii)$] The function $\psi$ satisfies:
	\begin{equation}
		\label{bdry:feas}
		\begin{aligned}
			&\mathbb{I}_{\Omega^{(1)}}(x^{(1)})\cdot\bm{e}_3\ge\psi(x_1,x_2),\quad\textup{ for all }x^{(1)}=(x_1,x_2,x_3^{(1)})\in\Gamma_0^{(1)},\\
			&\mathbb{I}_{\Omega^{(2)}}(x^{(2)})\cdot\bm{e}_3\le\psi(x_1,x_2),\quad\textup{ for all }x^{(2)}=(x_1,x_2,x_3^{(2)})\in\Gamma_0^{(2)}.
		\end{aligned}
	\end{equation}
\end{itemize}

Condition~$(i)$ suggests that the function $\psi$ can be interpreted as an \emph{optimal control}~\cite{Lio71} ensuring that the deformations of the linearly elastic bodies $\overline{\Omega^{(1)}}$ and $\overline{\Omega^{(2)}}$ do not interpenetrate each other.
In condition $(ii)$, we deliberately use non-strict inequalities, so that each deformed body may \emph{touch} the graph of the interface $\psi$.
In condition $(iii)$ the constraint is imposed at \emph{every} point of $\Gamma_0^{(\alpha)}$, as the map $x\mapsto\mathbb{I}_{\Omega^{(\alpha)}}(x)\cdot\bm{e}_3$ is continuous, and $\psi$ will be chosen in $\mathcal V(K)$ (defined later in this section), whose elements are continuous. By contrast, in~\eqref{eq:conf-cond} the displacement $\bm{v}^{(\alpha)}$ is only defined a.e., so the confinement is stated in the a.e.-sense.

The optimality in the search of the function $\psi$, which acts as separating interface between the deformations of the two linearly elastic bodies $\overline{\Omega^{(1)}}$ and $\overline{\Omega^{(2)}}$, is due to the fact that we want to find the surface which \emph{complies the most} with the contact geometry of Problem~\ref{problem0}.
Admissible displacements are thus vector fields $\bm{v}^{(\alpha)}:\overline{\Omega^{(\alpha)}}\to\mathbb{R}^3$ satisfying the following confinement conditions:
\begin{equation*}
	\begin{aligned}
		(\mathbb{I}_{\Omega^{(1)}}(x)+\bm{v}^{(1)}(x))\cdot\bm{e}_3 &\ge \psi((\mathbb{I}_{\Omega^{(1)}}(x)+\bm{v}^{(1)}(x))\cdot\bm{e}_1,(\mathbb{I}_{\Omega^{(1)}}(x)+\bm{v}^{(1)}(x))\cdot\bm{e}_2),\quad\textup{ for $\mathrm{d}\Gamma^{(1)}$-a.a. }x\in\Gamma^{(1)},\\
		(\mathbb{I}_{\Omega^{(2)}}(x)+\bm{v}^{(2)}(x))\cdot\bm{e}_3 &\le \psi((\mathbb{I}_{\Omega^{(2)}}(x)+\bm{v}^{(2)}(x))\cdot\bm{e}_1,(\mathbb{I}_{\Omega^{(2)}}(x)+\bm{v}^{(2)}(x))\cdot\bm{e}_2),\quad\textup{ for $\mathrm{d}\Gamma^{(2)}$-a.a. }x\in\Gamma^{(2)}.
	\end{aligned}
\end{equation*}

The first confinement condition amounts to requiring that $\mathrm{d}\Gamma^{(1)}$-a.a. the points of the boundary of the linearly elastic body $\overline{\Omega^{(1)}}$ lie above the graph of the given function $\psi$. The second inequality is the symmetric counterpart for $\overline{\Omega^{(2)}}$, whose deformed points must remain below the graph of $\psi$.

Since $\psi$ will be sought among functions with small gradient (condition $(i)$) and the displacements are small (this is the main feature characterising linearised elasticity~\cite{Ciarlet1988}), we can drop, in the same fashion as page~24 in~\cite{KikuchiOden1988}, the tangential displacement in the argument of $\psi$, which renders the constraint affine in $\bm{v}^{(\alpha)}$ for each fixed $\psi$, obtaining:
\begin{equation}
	\label{eq:conf-cond}
	\begin{aligned}
		(\mathbb{I}_{\Omega^{(1)}}(x^{(1)})+\bm{v}^{(1)}(x^{(1)}))\cdot\bm{e}_3 &\ge \psi(\mathbb{I}_{\Omega^{(1)}}(x^{(1)})\cdot\bm{e}_1,\mathbb{I}_{\Omega^{(1)}}(x^{(1)})\cdot\bm{e}_2),\quad\textup{ for a.a. }x^{(1)}\in\overline{\Omega^{(1)}},\\
		(\mathbb{I}_{\Omega^{(2)}}(x^{(2)})+\bm{v}^{(2)}(x^{(2)}))\cdot\bm{e}_3 &\le \psi(\mathbb{I}_{\Omega^{(2)}}(x^{(2)})\cdot\bm{e}_1,\mathbb{I}_{\Omega^{(2)}}(x^{(2)})\cdot\bm{e}_2),\quad\textup{ for a.a. }x^{(2)}\in\overline{\Omega^{(2)}}.
	\end{aligned}
\end{equation}

\begin{remark}
	We do \emph{not} infer the interior confinement from the boundary one. Indeed, a body whose boundary lies in $\textup{epi}\psi$ may have interior points below the graph of $\psi$. The confinement is therefore imposed \emph{directly} at a.e. point of $\Omega^{(\alpha)}$, as in \eqref{eq:conf-cond}. Besides we will see that defining the constraint in the interior will be essential for the calculation of the gradient of the cost functional for the optimal control problem formulated in this section, which we will need for numerically implementing the model (Section~\ref{Sec:5}).
	\bqed
\end{remark}

To begin with, we state the variational problem governing the displacement state associated with a \emph{given} function $\psi$, which acts as a \emph{control}.
We assume that the competitors to the role of optimal interface are functions $\psi\in H^2_0(\kappa)$. With a function $\psi$ of one such kind, we associate the function $\gls{3Dinterface}:K\to\mathbb{R}$ defined by
\begin{equation}
	\label{Psi}
	\Psi(x):=\psi(x_1,x_2),\quad\textup{ for all }x=(x_i)\in K,
\end{equation}
meaning that $\Psi$ is the three-dimensional analogue of $\psi$ and $\Psi$ is constant along the transverse variable $x_3$.

It is immediate to observe, however, that the formulation proposed in~\eqref{eq:conf-cond} is not always satisfied for linearly elastic bodies at rest. Indeed, it may happen that the graph of $\psi$ intersects the reference configuration of the linearly elastic body $\overline{\Omega^{(1)}}$ when no external forces are applied. This situation is reasonable to be considered if, for instance, the linearly elastic body $\overline{\Omega^{(1)}}$ undergoes the action of no applied forces while the linearly elastic body $\overline{\Omega^{(2)}}$ pushes so as to \emph{penetrate} the linearly elastic body $\overline{\Omega^{(1)}}$.
To address this issue, observe that $(iii)$ and Stampacchia's theorem (cf., e.g., Theorem~4.4 in Section~4.2.2 of~\cite{EG15}) give:
\begin{equation}
	\label{eq:bc-neg-pos-parts}
	\begin{aligned}
		\left\{(\mathbb{I}_{\Omega^{(1)}}-\Psi\bm{e}_3)\cdot\bm{e}_3\right\}^{-}\bm{e}_3\in\bm{V}(\Omega^{(1)}),\\
		-\left\{(\mathbb{I}_{\Omega^{(2)}}-\Psi\bm{e}_3)\cdot\bm{e}_3\right\}^{+}\bm{e}_3\in\bm{V}(\Omega^{(2)}).
	\end{aligned}
\end{equation}

In light of this considerations, we introduce the set $\gls{KO1psi}$ of admissible displacements for the linearly elastic body $\overline{\Omega^{(1)}}$ when subject to an applied body force $\bm{f}^{(1)}\in \bm{L}^2(\Omega^{(1)})$
\begin{equation}
	\label{eq:set-comp}
	\bm{K}(\Omega^{(1)};\Psi):=\left\{\bm{v}\in\bm{V}(\Omega^{(1)});\left(\mathbb{I}_{\Omega^{(1)}}+\bm{v}\right)\cdot\bm{e}_3 \ge \Psi, \textup{ a.e. in } \Omega^{(1)}\right\},
\end{equation}
where we observe that, thanks to~\eqref{eq:bc-neg-pos-parts}, it results $\left\{(\mathbb{I}_{\Omega^{(1)}}-\Psi\bm{e}_3)\cdot\bm{e}_3\right\}^{-}\bm{e}_3\in\bm{K}(\Omega^{(1)};\Psi)$.

The variational problem modelling the displacement of the linearly elastic body $\overline{\Omega^{(1)}}$ subject to remaining confined above the \emph{given} function $\psi\in H^2_0(\kappa)$, that naturally defines a function $\Psi\in H^2(K)$ that is constant with respect to the transverse variable, takes the following form.
\begin{customprob}{$\mathcal{P}(\Omega^{(1)};\Psi)$}\label{problem1}
	Find $\bm{u}^{(1)}=(u_i^{(1)})\in\bm{K}(\Omega^{(1)};\Psi)$ that satisfies the variational inequalities
	\begin{equation*}
		\int_{\Omega^{(1)}} A^{ijk\ell,1} e_{k\|\ell}(\bm{u}^{(1)}) e_{i\|j}(\bm{v}-\bm{u}^{(1)}) \,\mathrm{d}x^{(1)} \ge \int_{\Omega^{(1)}} \bm{f}^{(1)} \cdot (\bm{v}-\bm{u}^{(1)}) \,\mathrm{d}x^{(1)},
	\end{equation*}
	for all $\bm{v}=(v_i) \in \bm{K}(\Omega^{(1)};\Psi)$.
	\bqed
\end{customprob}

We begin by establishing the existence and uniqueness of solutions for Problem~\ref{problem1}.
\begin{theorem}
	\label{th:ex1}
	Problem~\ref{problem1} admits a unique solution.
\end{theorem}
\begin{proof}
	To establish the existence and uniqueness of solutions for Problem~\ref{problem1}, we verify the assumptions of Theorem~4.8-2 in~\cite{Cia25}. Thanks to the symmetry and uniform positive definiteness of the fourth order three-dimensional elasticity tensor $\{A^{ijk\ell,1}\}$, the bilinear form
	\begin{equation*}
		(\bm{u},\bm{v})\in\bm{V}(\Omega^{(1)})\times\bm{V}(\Omega^{(1)}) \mapsto \int_{\Omega^{(1)}} A^{ijk\ell,1} e_{k\|\ell}(\bm{u})e_{i\|j}(\bm{v})\,\mathrm{d}x^{(1)},
	\end{equation*}
	is symmetric and continuous. Thanks to the Korn inequality in Cartesian coordinates (Theorem~\ref{KornCart}), the bilinear form considered above is $\bm{V}(\Omega^{(1)})$-elliptic.
	The linear form
	\begin{equation*}
		\bm{v}\in\bm{V}(\Omega^{(1)})\mapsto\int_{\Omega^{(1)}}\bm{f}^{(1)}\cdot\bm{v}\,\mathrm{d}x^{(1)},
	\end{equation*}
	is continuous thanks to H\"older's inequality (cf., e.g., \cite{Brez11}). The set $\bm{K}(\Omega^{(1)};\Psi)$ introduced in~\eqref{eq:set-comp} is clearly non-empty, as we have already observed that $\left\{(\mathbb{I}_{\Omega^{(1)}}-\Psi\bm{e}_3)\cdot\bm{e}_3\right\}^{-}\bm{e}_3\in\bm{K}(\Omega^{(1)};\Psi)$, closed in $\bm{V}(\Omega^{(1)})$, and convex. Therefore, an application of Theorem~4.8-2 in~\cite{Cia25} gives the sought conclusion.
\end{proof}

In a completely similar fashion, solutions for the variational problem modelling the displacement of the linearly elastic body $\overline{\Omega^{(2)}}$ subject to a given applied body force $\bm{f}^{(2)}\in\bm{L}^2(\Omega^{(2)})$ remaining confined below the \emph{given} function $\psi\in H^1_0(\kappa)$ are to be sought in the set
\begin{equation*}
	\gls{KO2psi}:=\left\{\bm{v}\in\bm{V}(\Omega^{(2)});\left(\mathbb{I}_{\Omega^{(2)}}+\bm{v}\right)\cdot\bm{e}_3 \le \Psi, \textup{ a.e. in } \Omega^{(2)}\right\},
\end{equation*}
which is non-empty (as~\eqref{eq:bc-neg-pos-parts} implies $\left\{(\mathbb{I}_{\Omega^{(2)}}-\Psi\bm{e}_3)\cdot\bm{e}_3\right\}^{+}\bm{e}_3\in\bm{K}(\Omega^{(2)};\Psi)$), closed in $\bm{V}(\Omega^{(2)})$, and convex. One such variational problem can be formulated as follows and, in the same fashion as Theorem~\ref{th:ex1}, admits a unique solution.

\begin{customprob}{$\mathcal{P}(\Omega^{(2)};\Psi)$}\label{problem2}
	Find $\bm{u}^{(2)}=(u_i^{(2)})\in\bm{K}(\Omega^{(2)};\Psi)$ that satisfies the variational inequalities
	\begin{equation*}
		\int_{\Omega^{(2)}} A^{ijk\ell,2} e_{k\|\ell}(\bm{u}^{(2)}) e_{i\|j}(\bm{v}-\bm{u}^{(2)})\,\mathrm{d}x^{(2)}\ge\int_{\Omega^{(2)}}\bm{f}^{(2)}\cdot(\bm{v}-\bm{u}^{(2)})\,\mathrm{d}x^{(2)},
	\end{equation*}
	for all $\bm{v}=(v_i)\in\bm{K}(\Omega^{(2)};\Psi)$.
	\bqed	
\end{customprob}

Given a number $\beta>0$, the optimal control problem we formulate next aims to find the least energetically expensive interfaces separating the deformations of the linearly elastic bodies $\overline{\Omega^{(1)}}$ and $\overline{\Omega^{(2)}}$ solving Problem~\ref{problem1} and Problem~\ref{problem2}, respectively.
We define the set of admissible controls as:
\begin{equation*}
	\gls{VK}:=\left\{\Phi\in H^2(K); \Phi(x_1,x_2,x_3)=\psi(x_1,x_2)\textup{ for some }\psi\in H^2_0(\kappa)\textup{ satisfying~\eqref{bdry:feas}}\right\}.
\end{equation*}

Equivalently, $\Phi\in\mathcal{V}(K)$ if and only if $\Phi\in H^2(K)$, $\partial_3\Phi=0$ a.e. in $K$, $\Phi=\partial_{\bm\nu}\Phi=0$ $\mathrm{d}\mathcal{H}^2$-a.e. on $\partial\kappa\times(-a,a)$, and \eqref{bdry:feas} holds with $\psi:=\Phi(\cdot,\cdot,\cdot)$. Observe that $\Phi\equiv 0\in\mathcal{V}(K)$, which renders $\mathcal{V}(K)$ non-empty. It can also be observed that $\mathcal{V}(K)$ is convex and closed.
In line with Section~\ref{Sec:3}, we assume that
\begin{equation}
	\label{kappa:compact}
	\overline{\pi_\kappa(\Omega^{(1)})}\cup\overline{\pi_\kappa(\Omega^{(2)})}\subset\subset\kappa,
\end{equation}
that is, the projections of the reference configurations of the linearly elastic bodies under consideration are compactly contained in $\kappa$. The reason why we assume~\eqref{kappa:compact} is that we want the projections of the bodies to be compactly contained in $\kappa$. This guarantees that all admissible deformations remain in the compact container $K=\overline\kappa\times[-a,a]$ and the contact region lies away from $\partial\kappa$, on which the admissible controls are hypothesised to vanish.

In addition to \eqref{kappa:compact}, we assume that $\Gamma_0^{(\alpha)}$ is \emph{relatively closed} in $\Gamma^{(\alpha)}$ and that the following \emph{separation property}, employed in the construction of the recovery sequence for the $\Gamma$-convergence argument (Lemma~\ref{gamma-3}), holds:
\begin{equation}
	\label{hyp:sep}
	\textup{dist}\left(\overline{\omega},\pi_\kappa(\Gamma_0^{(1)})\cup\pi_\kappa(\Gamma_0^{(2)})\right)>0.
\end{equation}

The existence of separating controls that can be cut off near $\partial\kappa$ follows instead from the geometric hypotheses~\eqref{kappa:compact} and~\eqref{hyp:sep}, which state that $\overline\omega\subset\subset\kappa$ and that $\overline\omega$ lies at positive distance from $\pi_\kappa(\Gamma_0^{(1)})\cup\pi_\kappa(\Gamma_0^{(2)})$.
The announced optimal control problem takes the following form.
\begin{customprob}{$\mathscr{C}_\beta(K)$}\label{controlproblem}
	Find $\Psi_\beta\in\mathcal{V}(K)$ that minimises the reduced cost functional
	\begin{equation*}
		\begin{aligned}
		&\gls{Jbeta}(\Psi):=\dfrac{\beta}{2}\int_{K}|\Delta\Psi|^2\,\mathrm{d}x
		+\dfrac{1}{2}\int_{\Omega^{(1)}}A^{ijk\ell,1}e_{k\|\ell}(\bm{u}^{(1)}(\Psi))e_{i\|j}(\bm{u}^{(1)}(\Psi))\,\mathrm{d}x-\int_{\Omega^{(1)}}\bm{f}^{(1)}\cdot\bm{u}^{(1)}(\Psi)\,\mathrm{d}x\\
		&\quad+\dfrac{1}{2}\int_{\Omega^{(2)}}A^{ijk\ell,2}e_{k\|\ell}(\bm{u}^{(2)}(\Psi))e_{i\|j}(\bm{u}^{(2)}(\Psi))\,\mathrm{d}x-\int_{\Omega^{(2)}}\bm{f}^{(2)}\cdot\bm{u}^{(2)}(\Psi)\,\mathrm{d}x,
		\end{aligned}
	\end{equation*}
	subject to $\gls{u1psi}$ being the solution of Problem~\ref{problem1} and to $\gls{u2psi}$ being the solution of Problem~\ref{problem2}.
	\bqed
\end{customprob}

\begin{remark}[Physical motivation of the proposed modelling approach]
The reformulation of unilateral contact between two linearly elastic bodies as an optimal control problem, where the obstacle's action is regarded as a primary variable minimising a composite cost of deformation energy and control cost, offers an advantageous perspective. This approach is theoretically sound and practically beneficial for the following reasons:
\begin{itemize}
	\item[$(1)$] It provides a \emph{unified variational framework} that seamlessly integrates the non-penetration constraint into the energy minimization principle. By augmenting the classical potential energy with a Tikhonov regularization term $\frac{\beta}{2}\|\Delta\Psi\|_{L^2(K)}^2$, the inherently non-smooth, inequality-constrained problem is transformed into a smooth, unconstrained optimization in the control variable;
	\item[$(2)$] This formulation admits a compelling \emph{physical interpretation}. The fact that the parameter $\beta>0$ is small implies that there is \emph{a very small cost} to pay to choose a non-zero obstacle $\Psi$ over $\Psi\equiv 0$. Indeed, if we take one linearly elastic body into account, we know by the classical theory of minimisation of quadratic functionals (cf., e.g., Chapter~8 in~\cite{Cia25}) that there exists one and only one minimiser for one such quadratic energy functional. Therefore, any deformed reference configuration constrained not to cross the obstacle $\Psi$ is not energetically more convenient than the solution of the corresponding obstacle-free problem. As a result, despite the presence of the obstacle $\Psi$, the Problem~\ref{controlproblem} will produce a control and a state which are as close as possible to producing the solution of the corresponding obstacle-free problem;
	\item[$(3)$] From a computational standpoint, the approach offers an alternative to the classical penalty method which is in general not applicable to this kind of contact problems (viz. Remark~\ref{rem:comp} below).
\end{itemize}
\bqed
\end{remark}

While the existence and uniqueness of solutions for Problem~\ref{controlproblem} can be determined by resorting to the Direct Method of the Calculus of Variations (viz., e.g., Chapter~11 in~\cite{Cia25}), the continuity of the control-to-state operator is more challenging to establish as the non-linearity associated with the geometrical constraint makes a non-smoothness arise. For remedying this situation, we relax Problem~\ref{problem1} and Problem~\ref{problem2} so as to be able to make terms associated with the confinement condition defined in terms of the given function $\Psi$ appear in the governing equations. Before stating this relaxation, let us establish the existence and uniqueness of solutions for Problem~\ref{controlproblem}.
\begin{theorem}
	\label{ex-un-op-con}
	Let $\bm{f}^{(\alpha)}\in\bm{L}^2(\Omega^{(\alpha)})$ and let $\beta>0$ be given. Then, Problem~\ref{controlproblem} admits one and only one solution $\Psi_\beta\in\mathcal{V}(K)$.
\end{theorem}
\begin{proof}
	The cost functional $\mathcal{J}_\beta(\Psi)$ is a reduced version of the full cost functional
	\begin{equation*}
		\begin{aligned}
			&\mathcal{J}_\beta(\Psi,\bm{u}^{(1)},\bm{u}^{(2)}):=\dfrac{\beta}{2}\int_{K}|\Delta\Psi|^2\,\mathrm{d}x
			+\dfrac{1}{2}\int_{\Omega^{(1)}}A^{ijk\ell,1}e_{k\|\ell}(\bm{u}^{(1)})e_{i\|j}(\bm{u}^{(1)})\,\mathrm{d}x-\int_{\Omega^{(1)}}\bm{f}^{(1)}\cdot\bm{u}^{(1)}\,\mathrm{d}x\\
			&\quad+\dfrac{1}{2}\int_{\Omega^{(2)}}A^{ijk\ell,2}e_{k\|\ell}(\bm{u}^{(2)})e_{i\|j}(\bm{u}^{(2)})\,\mathrm{d}x-\int_{\Omega^{(2)}}\bm{f}^{(2)}\cdot\bm{u}^{(2)}\,\mathrm{d}x,
		\end{aligned}
	\end{equation*}
	which is strictly convex in all variables, being $\Psi\in \mathcal{V}(K)\mapsto\|\Delta\Psi\|_{L^2(K)}^2$ strictly convex in $\mathcal{V}(K)$ and being
	\begin{equation*}
		\bm{u}^{(\alpha)}\in\bm{V}(\Omega^{(\alpha)})\mapsto E_\alpha(\bm{u}^{(\alpha)}):=\int_{\Omega^{(\alpha)}}A^{ijk\ell,\alpha}e_{k\|\ell}(\bm{u}^{(\alpha)})e_{i\|j}(\bm{u}^{(\alpha)})\,\mathrm{d}x^{(\alpha)}-\int_{\Omega^{(\alpha)}}\bm{f}^{(\alpha)}\cdot\bm{u}^{(\alpha)}\,\mathrm{d}x^{(\alpha)}
	\end{equation*}
	strictly convex in $\bm{V}(\Omega^{(\alpha)})$. Additionally, an application of the Korn inequality in Cartesian coordinates (Theorem~\ref{KornCart}) and Theorem~8.7-1 in~\cite{Cia25} give that the functional $\mathcal{J}(\Psi,\bm{u}^{(1)},\bm{u}^{(2)})$ is coercive over $\mathcal{V}(K)\times\bm{V}(\Omega^{(1)})\times\bm{V}(\Omega^{(2)})$. Finally, the feasible set
	\begin{equation*}
		\begin{aligned}
			\bm{\mathcal{U}}:=&\Big\{(\Psi,\bm{u}^{(1)},\bm{u}^{(2)})\in\mathcal{V}(K)\times\bm{V}(\Omega^{(1)})\times\bm{V}(\Omega^{(2)});\\
			&\quad\left(\mathbb{I}_{\Omega^{(1)}}+\bm{u}^{(1)}\right)\cdot\bm{e}_3 \ge \Psi, \textup{ a.e. in } \Omega^{(1)}\textup{ and }\\
			&\quad\left(\mathbb{I}_{\Omega^{(2)}}+\bm{u}^{(2)}\right)\cdot\bm{e}_3 \le \Psi, \textup{ a.e. in } \Omega^{(2)}\Big\}
		\end{aligned}
	\end{equation*}
	is non-empty thanks to~\eqref{eq:bc-neg-pos-parts}, closed and convex. Therefore, we are in position to apply Theorem~11.7-1 and Theorem~9.13-3$(b)$ in~\cite{Cia25} so as to infer that the \emph{full} cost functional has at most one minimiser.
	
	For each fixed $\Psi\in\mathcal V(K)$, the unique solution $\bm{u}^{(\alpha)}(\Psi)$ of Problem~\ref{problem1} (resp. Problem~\ref{problem2}) is, by Stampacchia's theorem, the unique minimiser of $E_\alpha$ over $\bm{K}(\Omega^{(\alpha)};\Psi)$. Hence the minimiser $(\Psi_\beta,\bm{u}^{(1)}_\beta,\bm{u}^{(2)}_\beta)$ of the full cost functional satisfies $\bm{u}^{(\alpha)}_\beta=\bm{u}^{(\alpha)}(\Psi_\beta)$, so that
	\begin{equation*}
	\min_{\Psi\in\mathcal V(K)}\mathcal{J}_\beta(\Psi)=\min_{(\Psi,\bm{u}^{(1)},\bm{u}^{(2)})\in\bm{\mathcal{U}}}\mathcal{J}_\beta(\Psi,\bm u^{(1)},\bm u^{(2)}),
	\end{equation*}
	and $\mathcal{J}_\beta$ has exactly one minimiser $\Psi_\beta$.
\end{proof}

To complete this section, we establish a connection between Problem~\ref{controlproblem} and Problem~\ref{problem0}. Consider the energy functional $\gls{gammalimit}:\bm{L}^2(\Omega^{(1)})\times\bm{L}^2(\Omega^{(2)})\to\mathbb{R}\cup\{\infty\}$ defined by:
\begin{equation*}
	\tilde{\mathcal{F}}_0(\bm{u}^{(1)},\bm{u}^{(2)}):=
	\begin{cases}
		\sum_{\alpha=1}^{2}E_\alpha(\bm{u}^{(\alpha)})&,\textup{ if }(\bm{u}^{(1)},\bm{u}^{(2)})\in\vec{\bm{U}},\\
		\infty&,\textup{ otherwise}.
	\end{cases}
\end{equation*}

To begin with, observe that $\tilde{\mathcal{F}}_0$ admits a unique minimiser in light of Theorem~\ref{th-ex-un-0}.
For each $\beta>0$, define the marginal functional $\tilde{\mathcal{F}}_\beta:\bm{V}(\Omega^{(1)})\times\bm{V}(\Omega^{(2)})\to\mathbb{R}$ as follows
\begin{equation}
	\label{marginal}
	\gls{marginal}(\bm{u}^{(1)},\bm{u}^{(2)}):=\sum_{\alpha=1}^{2}E_\alpha(\bm{u}^{(\alpha)})+\inf_{\substack{\Psi\in\mathcal{V}(K);\\\bm{u}^{(1)}\in\bm{K}(\Omega^{(1)};\Psi)\\\bm{u}^{(2)}\in\bm{K}(\Omega^{(2)};\Psi)}}\left(\dfrac{\beta}{2}\|\Delta\Psi\|_{L^2(K)}^2\right),
\end{equation}
for all $(\bm{u}^{(1)},\bm{u}^{(2)})\in\bm{V}(\Omega^{(1)})\times\bm{V}(\Omega^{(2)})$. Observe that, given $(\bm{u}^{(1)},\bm{u}^{(2)})\in\bm{V}(\Omega^{(1)})\times\bm{V}(\Omega^{(2)})$, if there exists no element $\Psi\in\mathcal{V}(K)$ such that $\bm{u}^{(1)}\in\bm{K}(\Omega^{(1)};\Psi)$ and $\bm{u}^{(2)}\in\bm{K}(\Omega^{(2)};\Psi)$ then the infimum on the right-hand side of~\eqref{marginal} equals infinity. It is then clear that
\begin{equation}
	\label{id-min}
	\inf_{\substack{\bm{u}^{(1)}\in\bm{V}(\Omega^{(1)})\\\bm{u}^{(2)}\in\bm{V}(\Omega^{(2)})}}\tilde{\mathcal{F}}_\beta(\bm{u}^{(1)},\bm{u}^{(2)})=\min_{\Psi\in\mathcal{V}(K)}\mathcal{J}_\beta(\Psi),
\end{equation}
and the infimum on the left-hand side of~\eqref{id-min} is attained \emph{only} at the optimal states $(\bm{u}^{(1)}(\Psi_\beta),\bm{u}^{(2)}(\Psi_\beta))$ associated with the optimal control $\Psi_\beta$ that uniquely minimises the right-hand side of~\eqref{id-min}. Therefore, instead of carrying out a rigorous $\Gamma$-convergence analysis as $\beta\to0$ of the family $\{\mathcal{J}_\beta\}_{\beta>0}$, which would be challenging due to the dependence of the states on the controls in the reduced cost functional, we conduct one such $\Gamma$-convergence analysis over the family $\{\mathcal{F}_\beta\}_{\beta>0}$ of marginal functionals, where the analogues of the state variables are independent of the controls. Despite changing the departure point of the $\Gamma$-convergence analysis, we will still be able to show in the end that the family $\{(\bm{u}^{(1)}(\Psi_\beta),\bm{u}^{(2)}(\Psi_\beta))\}_{\beta>0}$ of minimisers for $\tilde{\mathcal{F}}_\beta$ weakly converges in $\bm{V}(\Omega^{(1)})\times\bm{V}(\Omega^{(2)})$ to the unique minimiser of $\tilde{\mathcal{F}}_0$ (cf., e.g., \cite{DalMaso93} and Theorem~11.18-1 in~\cite{Cia25}).
We now establish the $\Gamma$-convergence of the family of functionals $\{\tilde{\mathcal{F}}_\beta\}_{\beta>0}$ to the functional $\tilde{\mathcal{F}}_0$ as $\beta\to0$. We break the proof in three lemmas.

\begin{lemma}
	\label{gamma-1}
	The family of functionals $\{\tilde{\mathcal{F}}_\beta\}_{\beta>0}$ is equi-coercive (viz., e.g., Definition~7.6 in~\cite{DalMaso93}) in the weak topology of $\bm{V}(\Omega^{(1)})\times\bm{V}(\Omega^{(2)})$ and in the strong topology of $\bm{L}^2(\Omega^{(1)})\times\bm{L}^2(\Omega^{(2)})$.
\end{lemma}
\begin{proof}
	Let $C>0$ be given and let $\{\beta_n\}_{n=1}^\infty$ be a sequence such that $\beta_n\to0$ as $n\to\infty$. Let $(\bm{u}^{(1)},\bm{u}^{(2)})\in\bm{V}(\Omega^{(1)})\times\bm{V}(\Omega^{(2)})$ be such that $\tilde{\mathcal{F}}_{\beta_n}(\bm{u}^{(1)},\bm{u}^{(2)})\le C$ for all integers $n\ge 1$.
	An application of the uniform positive definiteness of the fourth order three-dimensional elasticity tensors $\{A^{ijk\ell,\alpha}\}$ (viz.~\eqref{eq:pos-def}), Korn's inequality in Cartesian coordinates (Theorem~\ref{KornCart}) and Young's inequality~\cite{Young1912} gives:
	\begin{equation}
		\label{bound1}
		\begin{aligned}
			&\tilde{\mathcal{F}}_{\beta_n}(\bm{u}^{(1)},\bm{u}^{(2)})\ge\sum_{\alpha=1}^2 E_\alpha(\bm{u}^{(\alpha)})\ge\dfrac{1}{2C_e^{(1)}}\sum_{i,j=1}^{3}\|e_{i\|j}(\bm{u}^{(1)})\|_{L^2(\Omega^{(1)})}^2+\dfrac{1}{2C_e^{(2)}}\sum_{i,j=1}^{3}\|e_{i\|j}(\bm{u}^{(2)})\|_{L^2(\Omega^{(2)})}^2\\
			&\quad-c_0^2 C_e^{(1)}\|\bm{f}^{(1)}\|_{\bm{L}^2(\Omega^{(1)})}^2-\dfrac{1}{4c_0^2C_e^{(1)}}\|\bm{u}^{(1)}\|_{\bm{V}(\Omega^{(1)})}^2\\
			&\quad-|c_0^{(1)}|^2 C_e^{(2)}\|\bm{f}^{(2)}\|_{\bm{L}^2(\Omega^{(2)})}^2-\dfrac{1}{4|c_0^{(1)}|^2C_e^{(2)}}\|\bm{u}^{(2)}\|_{\bm{V}(\Omega^{(2)})}^2\\
			&\ge\dfrac{1}{4|c_0^{(1)}|^2 C_e^{(1)}}\|\bm{u}^{(1)}\|_{\bm{V}(\Omega^{(1)})}^2+\dfrac{1}{4|c_0^{(2)}|^2 C_e^{(2)}}\|\bm{u}^{(2)}\|_{\bm{V}(\Omega^{(2)})}^2\\
			&\quad-|c_0^{(1)}|^2 C_e^{(1)}\|\bm{f}^{(1)}\|_{\bm{L}^2(\Omega^{(1)})}^2-|c_0^{(2)}|^2 C_e^{(2)}\|\bm{f}^{(2)}\|_{\bm{L}^2(\Omega^{(2)})}^2.
		\end{aligned}
	\end{equation}
	
	Thanks to~\eqref{bound1}, we obtain that:
	\begin{equation*}
		\begin{aligned}
			\|\bm{u}^{(1)}\|_{\bm{V}(\Omega^{(1)})}^2&\le 4C |c_0^{(1)}|^2 C_e^{(1)}+4 |c_0^{(1)}|^4 |C_e^{(1)}|^2\|\bm{f}^{(1)}\|_{\bm{L}^2(\Omega^{(1)})}^2+4 |c_0^{(1)}|^2|c_0^{(2)}|^2 C_e^{(1)}C_e^{(2)}\|\bm{f}^{(2)}\|_{\bm{L}^2(\Omega^{(2)})}^2,\\
			\|\bm{u}^{(2)}\|_{\bm{V}(\Omega^{(2)})}^2&\le 4C |c_0^{(2)}|^2 C_e^{(2)}+4 |c_0^{(1)}|^2 |c_0^{(2)}|^2 C_e^{(1)} C_e^{(2)}\|\bm{f}^{(1)}\|_{\bm{L}^2(\Omega^{(1)})}^2+4 |c_0^{(2)}|^4 |C_e^{(2)}|^2\|\bm{f}^{(2)}\|_{\bm{L}^2(\Omega^{(2)})}^2.
		\end{aligned}
	\end{equation*}
	
	Define the number
	\begin{equation*}
		\begin{aligned}
			R:=&4 |c_0^{(1)}|^4|C_e^{(1)}|^2\|\bm{f}^{(1)}\|_{\bm{L}^2(\Omega^{(1)})}^2+4 |c_0^{(1)}|^2 |c_0^{(2)}|^2C_e^{(1)}C_e^{(2)}\left(\|\bm{f}^{(1)}\|_{\bm{L}^2(\Omega^{(1)})}^2+\|\bm{f}^{(2)}\|_{\bm{L}^2(\Omega^{(2)})}^2\right)\\
			&\quad+4 |c_0^{(2)}|^4 |C_e^{(2)}|^2\|\bm{f}^{(2)}\|_{\bm{L}^2(\Omega^{(2)})}^2+4C |c_0^{(1)}|^2 C_e^{(1)}+4C |c_0^{(2)}|^2 C_e^{(2)},
		\end{aligned}
	\end{equation*}
	and observe that $\{(\bm{u}^{(1)},\bm{u}^{(2)})\in\bm{V}(\Omega^{(1)})\times\bm{V}(\Omega^{(2)});\tilde{\mathcal{F}}_\beta(\bm{u}^{(1)},\bm{u}^{(2)})\le C\}\subset B(0_{\bm{V}(\Omega^{(1)})\times\bm{V}(\Omega^{(2)})};R)$. Since the ball $B(0_{\bm{V}(\Omega^{(1)})\times\bm{V}(\Omega^{(2)})};R)$ is relatively compact in in the weak topology of $\bm{V}(\Omega^{(1)})\times\bm{V}(\Omega^{(2)})$ and in the strong topology of $\bm{L}^2(\Omega^{(1)})\times\bm{L}^2(\Omega^{(2)})$, the conclusion follows.
\end{proof}

The next intermediary result is the so-called $\liminf$ inequality.
\begin{lemma}
	\label{gamma-2}
	Let $\{(\bm{u}^{(1)}_n,\bm{u}^{(2)}_n)\}_{n=1}^\infty$ be a sequence that weakly converges to $(\bm{u}^{(1)},\bm{u}^{(2)})$ in $\bm{V}(\Omega^{(1)})\times\bm{V}(\Omega^{(2)})$ as $n\to\infty$. Then, for every sequence $\{\beta_n\}_{n=1}^\infty$ such that $\beta_n\to 0$ as $n\to\infty$, we have:
	\begin{equation*}
		\tilde{\mathcal{F}}_0(\bm{u}^{(1)},\bm{u}^{(2)})\le\liminf_{n\to\infty}\tilde{\mathcal{F}}_{\beta_n}(\bm{u}^{(1)}_n,\bm{u}^{(2)}_n).
	\end{equation*}
\end{lemma}
\begin{proof}
	If $\liminf_{n\to\infty}\tilde{\mathcal{F}}_{\beta_n}(\bm{u}^{(1)}_n,\bm{u}^{(2)}_n)=+\infty$, there is nothing to prove. Otherwise, extract a subsequence, still denoted by
	$\{(\bm{u}^{(1)}_n,\bm{u}^{(2)}_n)\}_{n=1}^\infty$, for which the $\liminf$ is attained. For one such subsequence, it results that there exists $C>0$ such that:
	\begin{equation*}
		\tilde{\mathcal{F}}_{\beta_n}(\bm{u}^{(1)}_n,\bm{u}^{(2)}_n)\le C, \quad\textup{ for all }n\ge1.
	\end{equation*}
	
	Since each $E_\alpha$ is bounded below on $\bm{V}(\Omega^{(\alpha)})$ (Lemma~\ref{gamma-1}), the finiteness of each $\tilde{\mathcal{F}}_{\beta_n}(\bm{u}^{(1)}_n,\bm{u}^{(2)}_n)$ implies that the infimum on the right-hand side of \eqref{marginal} is finite. Therefore, the definition of infimum ensures that, for each element $(\bm{u}^{(1)}_n,\bm{u}^{(2)}_n)$ of the sequence under consideration, there exists $\Psi_n\in\mathcal{V}(K)$ which depends on the choice of $(\bm{u}^{(1)}_n,\bm{u}^{(2)}_n)$ and is such that:
	\begin{equation*}
		\begin{aligned}
			(\mathbb{I}_{\Omega^{(1)}}+\bm{u}^{(1)}_n)\cdot\bm{e}_3\ge\Psi_n&,\quad\textup{ a.e. in }\Omega^{(1)},\\
			(\mathbb{I}_{\Omega^{(2)}}+\bm{u}^{(2)}_n)\cdot\bm{e}_3\le\Psi_n&,\quad\textup{ a.e. in }\Omega^{(2)}.
		\end{aligned}
	\end{equation*}
	
	Writing $\Psi_n=\psi_n\circ\pi$ with $\psi_n\in H^2_0(\kappa)$, these two inequalities imply that $(\bm{u}^{(1)}_n,\bm{u}^{(2)}_n)\in\vec{\bm{U}}$ for every $n\ge1$.
	As $\vec{\bm{U}}$ is non-empty, closed and convex (Lemma~\ref{lem:U-step2}), it is weakly closed in $\bm{V}(\Omega^{(1)})\times\bm{V}(\Omega^{(2)})$, whence the weak limit $(\bm{u}^{(1)},\bm{u}^{(2)})$ belongs to $\vec{\bm{U}}$. Consequently, we obtain:
	\begin{equation*}
		\tilde{\mathcal{F}}_0(\bm{u}^{(1)},\bm{u}^{(2)})=\sum_{\alpha=1}^{2}E_\alpha(\bm{u}^{(\alpha)}).
	\end{equation*}
	
	Each $E_\alpha$ is convex and continuous on $\bm{V}(\Omega^{(\alpha)})$ and, hence, sequentially weakly lower semicontinuous, so that it results:
	\begin{equation*}
		\sum_{\alpha=1}^{2}E_\alpha(\bm{u}^{(\alpha)})\le\sum_{\alpha=1}^{2}\liminf_{n\to\infty}E_\alpha(\bm{u}^{(\alpha)}_n).
	\end{equation*}
	
	Finally, since the infimum in~\eqref{marginal} is non-negative,
	\begin{equation*}
		\tilde{\mathcal{F}}_{\beta_n}(\bm{u}^{(1)}_n,\bm{u}^{(2)}_n)\ge E_1(\bm{u}^{(1)}_n)+E_2(\bm{u}^{(2)}_n),
	\end{equation*}
	whence:
	\begin{equation*}
		\sum_{\alpha=1}^{2}\liminf_{n\to\infty}E_\alpha(\bm{u}^{(\alpha)}_n)\le\liminf_{n\to\infty}\left(E_1(\bm{u}^{(1)}_n)+E_2(\bm{u}^{(2)}_n)\right)\le\liminf_{n\to\infty}\tilde{\mathcal{F}}_{\beta_n}(\bm{u}^{(1)}_n,\bm{u}^{(2)}_n).
	\end{equation*}
	
	Combining the three inequalities above yields the claimed $\liminf$ inequality.
	\end{proof}

Finally, we show that there exists a \emph{recovery sequence} that verifies the $\limsup$ inequality.
\begin{lemma}
	\label{gamma-3}
	Let $(\bm{u}^{(1)},\bm{u}^{(2)})$ in $\bm{V}(\Omega^{(1)})\times\bm{V}(\Omega^{(2)})$ be such that $\tilde{\mathcal{F}}_0(\bm{u}^{(1)},\bm{u}^{(2)})<\infty$. Then, for every sequence $\{\beta_n\}_{n=1}^\infty$ such that $\beta_n\to 0$ as $n\to\infty$, there exists a sequence $\{(\bm{u}^{(1)}_n,\bm{u}^{(2)}_n)\}_{n=1}^\infty$ that strongly converges to $(\bm{u}^{(1)},\bm{u}^{(2)})$ in $\bm{V}(\Omega^{(1)})\times\bm{V}(\Omega^{(2)})$ as $n\to\infty$ for which:
	\begin{equation*}
		\tilde{\mathcal{F}}_0(\bm{u}^{(1)},\bm{u}^{(2)})=\lim_{n\to\infty}\tilde{\mathcal{F}}_{\beta_n}(\bm{u}^{(1)}_n,\bm{u}^{(2)}_n).
	\end{equation*}
\end{lemma}
\begin{proof}
	The proof is divided into four steps, numbered $(i)$--$(iv)$. In the context of this proof, we will work with functions $\psi\in H^2_0(\kappa)$ rather than their three-dimensional analogues $\Psi\in\mathcal{V}(K)$ defined in~\eqref{Psi}.
	
	$(i)$ \emph{The set
		\begin{equation*}
			\begin{aligned}
				\vec{\bm{U}}^{(1)}:=\Big\{&(\bm{u}^{(1)},\bm{u}^{(2)})\in\vec{\bm{U}};
				\textup{ there exists }\delta>0 \textup{ such that}\\
				&\quad(\mathbb{I}_{\Omega^{(1)}}(\xi_1,\xi_2,x_3^{(1)})+\bm{u}^{(1)}(\xi_1,\xi_2,x_3^{(1)}))\cdot\bm{e}_3\ge(\mathbb{I}_{\Omega^{(2)}}(\xi_1,\xi_2,x_3^{(2)})+\bm{u}^{(2)}(\xi_1,\xi_2,x_3^{(2)}))\cdot\bm{e}_3+\delta\\
				&\quad\textup{ for }\mathrm{d}\mathcal{L}^2\textup{-a.e. }\xi=(\xi_1,\xi_2)\in\omega,\textup{ for }\mathrm{d}\mathcal{L}^1\textup{-a.a. }x_3^{(1)}\in F^{(1)}(\xi)\textup{  and }\mathrm{d}\mathcal{L}^1\textup{-a.a. }x_3^{(2)}\in F^{(2)}(\xi)\Big\},
			\end{aligned}
		\end{equation*}
		is dense in $\vec{\bm{U}}$.}
		Let $(\bm{u}^{(1)},\bm{u}^{(2)})\in\vec{\bm{U}}$. Thanks to the assumption~\eqref{hyp:sep} we made beforehand, there exists a cut-off function $\psi\in\mathcal{D}(\kappa)$ such that:
		\begin{equation*}
			\begin{aligned}
				0\le\psi\le1&,\quad\textup{ in }\overline{\kappa},\\
				\psi\equiv 1&,\quad\textup{ in a neighbourhood of }\overline{\omega},\\
				\psi\equiv 0&,\quad\textup{ in a neighbourhood of }\pi_\kappa(\Gamma_0^{(1)})\cup\pi_\kappa(\Gamma_0^{(2)}).
			\end{aligned}
		\end{equation*}
		
		With the customary abuse of notation, we still denote by $\psi$ the function $x=(x_1,x_2,x_3)\mapsto\psi(x_1,x_2)$. For each integer $n\ge1$, define:
		\begin{equation*}
		\begin{aligned}
			\bm{u}^{(1)}_n:=\bm{u}^{(1)}+\dfrac{\psi}{n}\bm{e}_3&, \quad\textup{ in }\Omega^{(1)},\\
			\bm{u}^{(2)}_n:=\bm{u}^{(2)}-\dfrac{\psi}{n}\bm{e}_3&, \quad\textup{ in }\Omega^{(2)}.
		\end{aligned}
		\end{equation*}
		
		Since $\psi$ vanishes on a neighbourhood of $\pi_\kappa(\Gamma_0^{(\alpha)})$ and $\bm{u}^{(\alpha)}\in\bm{V}(\Omega^{(\alpha)})$, it results that $\bm{u}^{(\alpha)}_n\in\bm{V}(\Omega^{(\alpha)})$ for each integer $n\ge 1$ and for each $\alpha\in\{1,2\}$. Additionally, it results
		\begin{equation*}
			\|\bm{u}^{(\alpha)}_n-\bm{u}^{(\alpha)}\|_{\bm{V}(\Omega^{(\alpha)})}=\dfrac{1}{n}\|\psi\|_{H^1(\Omega^{(\alpha)})}\to0,\quad\textup{ as }n\to\infty,
		\end{equation*}
		whence the sequence $\{(\bm{u}^{(1)}_n,\bm{u}^{(2)}_n)\}_{n=1}^\infty$ strongly converges to $(\bm{u}^{(1)},\bm{u}^{(2)})$ in $\bm{V}(\Omega^{(1)})\times\bm{V}(\Omega^{(2)})$ as $n\to\infty$.
		
		For $\mathrm{d}\mathcal{L}^2$-a.e. $\xi\in\omega$, for $\mathrm{d}\mathcal{L}^1$-a.e. $x_3^{(1)}\in F^{(1)}(\xi)$ and $\mathrm{d}\mathcal{L}^1$-a.e. $x_3^{(2)}\in F^{(2)}(\xi)$, the assumed admissibility of $(\bm{u}^{(1)},\bm{u}^{(2)})\in\vec{\bm{U}}$ gives:
		\begin{equation*}
		(\mathbb{I}_{\Omega^{(1)}}(\xi,x_3^{(1)})+\bm{u}^{(1)}(\xi,x_3^{(1)}))\cdot\bm{e}_3\ge(\mathbb{I}_{\Omega^{(2)}}(\xi,x_3^{(2)})+\bm{u}^{(2)}(\xi,x_3^{(2)}))\cdot\bm{e}_3.
		\end{equation*}
		
		Since $\psi\equiv1$ in $\overline{\omega}$, we obtain
		\begin{equation*}
		\begin{aligned}
			&(\mathbb{I}_{\Omega^{(1)}}(\xi,x_3^{(1)})+\bm{u}^{(1)}_n(\xi,x_3^{(1)}))\cdot\bm{e}_3-(\mathbb{I}_{\Omega^{(2)}}(\xi,x_3^{(2)})+\bm{u}^{(2)}_n(\xi,x_3^{(2)}))\cdot\bm{e}_3\\
			&\quad=\left[(\mathbb{I}_{\Omega^{(1)}}(\xi,x_3^{(1)})+\bm{u}^{(1)}(\xi,x_3^{(1)}))\cdot\bm{e}_3-(\mathbb{I}_{\Omega^{(2)}}(\xi,x_3^{(2)})+\bm{u}^{(2)}(\xi,x_3^{(2)}))\cdot\bm{e}_3\right]+\dfrac{2\psi}{n}\ge\dfrac{2\psi}{n}=\dfrac{2}{n}>0,
		\end{aligned}
		\end{equation*}
		for $\mathrm{d}\mathcal{L}^2$-a.e. $\xi\in\omega$, for $\mathrm{d}\mathcal{L}^1$-a.e. $x_3^{(1)}\in F^{(1)}(\xi)$ and $\mathrm{d}\mathcal{L}^1$-a.e. $x_3^{(2)}\in F^{(2)}(\xi)$. Therefore $(\bm{u}^{(1)}_n,\bm{u}^{(2)}_n)\in\vec{\bm{U}}^{(1)}$ with $\delta=\frac2n$ for every integer $n\ge1$, and $\vec{\bm{U}}^{(1)}$ is dense in $\vec{\bm{U}}$. This proves $(i)$.
	
	$(ii)$ \emph{The set
	\begin{equation*}
		\begin{aligned}
			\vec{\bm{U}}^{(2)}:=\Big\{&(\bm{u}^{(1)},\bm{u}^{(2)})\in\vec{\bm{U}}^{(1)}
			\cap\left(\bm{\mathcal{C}}^\infty(\overline{\Omega^{(1)}})\times\bm{\mathcal{C}}^\infty(\overline{\Omega^{(2)}})\right);\\
			&\textup{there exists }\psi=\psi(\bm{u}^{(1)},\bm{u}^{(2)})\in\mathcal{D}(\kappa)
			\textup{ such that}\\
			&(\mathbb{I}_{\Omega^{(2)}}(x^{(2)})+\bm{u}^{(2)}(x^{(2)}))\cdot\bm{e}_3<\psi(x^{(2)}_1,x^{(2)}_2),\quad\textup{for }\mathrm{d}\mathcal{L}^3\textup{-a.a. }x^{(2)}\in\Omega^{(2)},\\
			&\psi(x^{(1)}_1,x^{(1)}_2)<(\mathbb{I}_{\Omega^{(1)}}(x^{(1)})+\bm{u}^{(1)}(x^{(1)}))\cdot\bm{e}_3,\quad\textup{for }\mathrm{d}\mathcal{L}^3\textup{-a.a. }x^{(1)}\in\Omega^{(1)}\Big\}
		\end{aligned}
	\end{equation*}
	is dense in $\vec{\bm{U}}^{(1)}$.}
	
	Let $(\bm{u}^{(1)},\bm{u}^{(2)})\in\vec{\bm{U}}^{(1)}$, and let $\delta>0$ be such
	that
	\begin{equation*}
	(\mathbb{I}_{\Omega^{(1)}}(\xi_1,\xi_2,x_3^{(1)})+\bm{u}^{(1)}(\xi_1,\xi_2,x_3^{(1)}))\cdot\bm{e}_3\ge(\mathbb{I}_{\Omega^{(2)}}(\xi_1,\xi_2,x_3^{(2)})+\bm{u}^{(2)}(\xi_1,\xi_2,x_3^{(2)}))\cdot\bm{e}_3+\delta,
	\end{equation*}
	for $\mathrm{d}\mathcal{L}^2$-a.e. $\xi=(\xi_1,\xi_2)\in\omega$, for $\mathrm{d}\mathcal{L}^1$-a.e. $x_3^{(1)}\in F^{(1)}(\xi)$ and $\mathrm{d}\mathcal{L}^1$-a.e. $x_3^{(2)}\in F^{(2)}(\xi)$. The proof has two sub-steps: first we regularise the pair, then we construct a separating control.
	
	\underline{\emph{Sub-step 1: Regularisation.}} Since $\Omega^{(\alpha)}$ is a Lipschitz domain, there exists a larger Lipschitz domain $\Omega^{(\alpha)}\subset\subset\hat{\Omega}^{(\alpha)}$. Choose $0<\eta<\mathrm{dist}(\partial\Omega^{(\alpha)},\partial\hat{\Omega}^{(\alpha)})$ and set:
	\begin{equation*}
		\omega^\eta:=\{\xi\in\kappa:\ \mathrm{dist}(\xi,\omega)<\eta\}.
	\end{equation*}
	
	By~\eqref{hyp:sep}, after shrinking $\eta$ if necessary, we obtain that the set $\pi_\kappa(\Gamma_0^{(1)})\cup\pi_\kappa(\Gamma_0^{(2)})$ lies at positive distance from $\omega^\eta$.
	We extend each displacement to $\hat{\Omega}^{(\alpha)}$ by two independent even reflections across the boundary. Tangentially, we reflect across the lateral part of $\partial\Omega^{(\alpha)}$ in the horizontal direction $(x_1,x_2)$. This keeps $x_3$ fixed, and the fields become defined at every tangential coordinate of $\omega^\eta$. For $\xi\in\omega$, define:
	\begin{equation*}
		\begin{aligned}
			\ell(\xi)&:=\textup{essinf}_{x_3\in F^{(1)}(\xi)}\left\{\left(\mathbb{I}_{\Omega^{(1)}}(\xi_1,\xi_2,x_3)+\bm{u}^{(1)}(\xi_1,\xi_2,x_3)\right)\cdot\bm{e}_3\right\},\\
			m(\xi)&:=\textup{esssup}_{x_3\in F^{(2)}(\xi)}\left\{\left(\mathbb{I}_{\Omega^{(2)}}(\xi_1,\xi_2,x_3)+\bm{u}^{(2)}(\xi_1,\xi_2,x_3)\right)\cdot\bm{e}_3\right\}.
		\end{aligned}
	\end{equation*}
	
	The gap condition satisfied by $(\bm{u}^{(1)},\bm{u}^{(2)})\in\vec{\bm{U}}^{(1)}$ reads:
	\begin{equation*}
		\ell(\xi)\ge m(\xi)+\delta,\qquad\textup{ for }\mathrm{d}\mathcal{L}^2\text{-a.a. }\xi\in\omega.
	\end{equation*}
	
	The tangential reflection copies the fields at fixed height $x_3$, so the envelopes of the tangentially extended fields, still denoted $\ell$ and $m$, are the reflection of the above two functions across $\partial\omega$. Since this reflection is one-to-one, we obtain:
	\begin{equation*}
		\ell(\xi)\ge m(\xi)+\delta,\qquad\textup{ for }\mathrm{d}\mathcal{L}^2\textup{-a.a. }\xi\in\omega^\eta.
	\end{equation*}
	
	For $\xi\in\omega^\eta$ denote by $\tilde{F}^{(1)}(\xi)$ and $\tilde{F}^{(2)}(\xi)$ the fibres of the tangentially extended bodies. Vertically, we reflect the tangentially extended third component of the displacement across $\inf \tilde{F}^{(1)}(\xi)$ for $\Omega^{(1)}$ and across $\sup \tilde{F}^{(2)}(\xi)$ for $\Omega^{(2)}$.
	Confining this vertical reflection to a layer adjacent to $\partial\Omega^{(\alpha)}$ thin enough, the deformed heights are perturbed by at most $\delta/4$. On the remaining parts of $\partial\Omega^{(\alpha)}$ we extend by an arbitrary bounded $H^1$ extension given by the standard extension theorem (cf., e.g., \cite{Brez11}). The previous procedure yields an extension $\hat{\bm{u}}^{(\alpha)}\in\bm{H}^1(\hat{\Omega}^{(\alpha)})$ of $\bm{u}^{(\alpha)}$ such that
	\begin{equation}
		\label{eq:ext-props}
		\left(\mathbb{I}_{\hat{\Omega}^{(1)}}(\xi_1,\xi_2,x_3^{(1)})+\hat{\bm{u}}^{(1)}(\xi_1,\xi_2,x_3^{(1)})\right)\cdot\bm{e}_3-\left(\mathbb{I}_{\hat{\Omega}^{(2)}}(\xi_1,\xi_2,x_3^{(2)})+\hat{\bm{u}}^{(2)}(\xi_1,\xi_2,x_3^{(2)})\right)\cdot\bm{e}_3\ge\dfrac{\delta}{2},
	\end{equation}
	for $\mathrm{d}\mathcal{L}^2$-a.a. $\xi\in\omega^\eta$, for $\mathrm{d}\mathcal{L}^1$-a.a. $x_3^{(1)}\in\hat{F}^{(1)}(\xi)$ and $\mathrm{d}\mathcal{L}^1$-a.a. $x_3^{(2)}\in\hat{F}^{(2)}(\xi)$. Indeed, the tangential reflection preserves the inequality $\ell\ge m+\delta$ on $\omega^\eta$, and the vertical reflection loses at most $\delta/4$ on each side, so the remaining gap is $\delta-\delta/4-\delta/4=\delta/2$.
	
	We now construct two auxiliary sequences. First, by the density theorem of~\cite{Dok73,DZ06} (see also Theorem~6.2 in~\cite{Pie26-3}) there exists a sequence $\bm{g}^{(\alpha)}_n\in\bm{\mathcal{C}}^\infty(\overline{\Omega^{(\alpha)}})$ such that:
	\begin{equation*}
		\begin{aligned}
			\bm{g}^{(\alpha)}_n=\bm{0}&,\quad\text{ on }\Gamma_0^{(\alpha)},\\
			\bm{g}^{(\alpha)}_n\to\bm{u}^{(\alpha)}&,\quad\textup{ in }\bm{H}^1(\Omega^{(\alpha)})\textup{ as }n\to\infty.
		\end{aligned}
	\end{equation*}
	
	Second, let $\rho\in\mathcal{D}(\mathbb{R}^3)$ be a mollifier (cf., e.g., \cite{Brez11}), and define
	\begin{equation}
		\label{eq:h}
		\bm{h}^{(\alpha)}_n:=\left(\rho_{1/n}\star\hat{\bm{u}}^{(\alpha)}\right)\big|_{\Omega^{(\alpha)}},
	\end{equation}
	where $\hat{\bm{u}}^{(\alpha)}$ is extended by zero outside $\hat{\Omega}^{(\alpha)}$. It results that $\bm{h}^{(\alpha)}_n\in\bm{\mathcal{C}}^\infty(\overline{\Omega^{(\alpha)}})$ and:
	\begin{equation*}
		\bm{h}^{(\alpha)}_n\to\bm{u}^{(\alpha)}\quad\textup{ in }\bm{H}^1(\Omega^{(\alpha)})\textup{ as }n\to\infty.
	\end{equation*}
	
	Since $\omega^\eta$ and $\pi_\kappa(\Gamma_0^{(1)})\cup\pi_\kappa(\Gamma_0^{(2)})$ are at positive distance, there is a fixed cut-off $\varphi\in\mathcal{C}^\infty(\kappa)$ with $0\le\varphi\le1$, $\varphi=1$ on $\omega^\eta$ and $\varphi=0$ on a neighbourhood of $\pi_\kappa(\Gamma_0^{(1)})\cup\pi_\kappa(\Gamma_0^{(2)})$. With the abuse of notation $\varphi(x):=\varphi(\pi_\kappa(x))$, define:
	\begin{equation*}
		\bm{u}^{(\alpha)}_n:=\varphi\,\bm{h}^{(\alpha)}_n+(1-\varphi)\,\bm{g}^{(\alpha)}_n,\quad\textup{ in }\overline{\Omega^{(\alpha)}}.
	\end{equation*}
	
	Both $\bm{h}^{(\alpha)}_n$ and $\bm{g}^{(\alpha)}_n$ are smooth, and the function $\varphi\circ\pi_\kappa$ is smooth because it is constant in $x_3$. Hence, it results that $\bm{u}^{(\alpha)}_n\in\bm{\mathcal{C}}^\infty(\overline{\Omega^{(\alpha)}})$. On $\Gamma_0^{(\alpha)}$ we have $\varphi=0$ and $\bm{g}^{(\alpha)}_n=\bm{0}$, so $\bm{u}^{(\alpha)}_n=\bm{0}$ on $\Gamma_0^{(\alpha)}$, that is, $\bm{u}^{(\alpha)}_n\in\bm{V}(\Omega^{(\alpha)})$.
	
	For $\mathrm{d}\mathcal{L}^2$-a.a. $\xi\in\omega$, for $\mathrm{d}\mathcal{L}^1$-a.a. $x_3^{(1)}\in F^{(1)}(\xi)$ and $\mathrm{d}\mathcal{L}^1$-a.a. $x_3^{(2)}\in F^{(2)}(\xi)$, the tangential projection $\xi$ lies in $\omega\subset\omega^\eta$, so $\varphi=1$ there and $\bm{u}^{(\alpha)}_n=\bm{h}^{(\alpha)}_n$. Therefore, by \eqref{eq:h}, we obtain that:
	\begin{equation}
		\label{eq:gap-conv}
		\begin{aligned}
			&\left(\mathbb{I}_{\Omega^{(1)}}(\xi_1,\xi_2,x_3^{(1)})+\bm{u}^{(1)}_n(\xi_1,\xi_2,x_3^{(1)})\right)\cdot\bm{e}_3-\left(\mathbb{I}_{\Omega^{(2)}}(\xi_1,\xi_2,x_3^{(2)})+\bm{u}^{(2)}_n(\xi_1,\xi_2,x_3^{(2)})\right)\cdot\bm{e}_3\\
			&=\int_{B(0;1/n)}\rho_{1/n}(z)\Big[\left(\mathbb{I}_{\hat{\Omega}^{(1)}}(\xi_1-z_1,\xi_2-z_2,x_3^{(1)}-z_3)+\hat{\bm{u}}^{(1)}(\xi_1-z_1,\xi_2-z_2,x_3^{(1)}-z_3)\right)\cdot\bm{e}_3\\
			&\qquad-\left(\mathbb{I}_{\hat{\Omega}^{(2)}}(\xi_1-z_1,\xi_2-z_2,x_3^{(2)}-z_3)+\hat{\bm{u}}^{(2)}(\xi_1-z_1,\xi_2-z_2,x_3^{(2)}-z_3)\right)\cdot\bm{e}_3\Big]\,\mathrm{d}z\ge\dfrac{\delta}{2}.
		\end{aligned}
	\end{equation}
	
	Note in passing that the identity in~\eqref{eq:gap-conv} carries no correction term because $(\mathbb{I}_{\hat{\Omega}^{(\alpha)}}\cdot\bm{e}_3)$ is affine with respect to $x_3$, so the two shifts by $z_3$ cancel, and $\int_{B(0;1/n)}\rho_{1/n}\,\mathrm{d}z=1$. The inequality follows from~\eqref{eq:ext-props} as the choice $0<\eta<\mathrm{dist}(\partial\Omega^{(\alpha)},\partial\hat{\Omega}^{(\alpha)})$ guarantees that $(\xi_1-z_1,\xi_2-z_2,x_3^{(\alpha)}-z_3)\in\hat{\Omega}^{(\alpha)}$ whenever $|z|\le1/n\le\eta$, and by Fubini's theorem (cf., e.g., Theorem~1.15-6$(b)$ in~\cite{Cia25}) the integrand inequality holds for a.a. $z$ at a.e. base point. Thus $(\bm{u}^{(1)}_n,\bm{u}^{(2)}_n)\in\vec{\bm{U}}^{(1)}$ and the gap extent is $\delta/2$.
	
	Since both $\bm{h}^{(\alpha)}_n$ and $\bm{g}^{(\alpha)}_n$ converge strongly to $\bm{u}^{(\alpha)}$ in $\bm{H}^1(\Omega^{(\alpha)})$, and $\varphi$ is a fixed smooth function, we obtain that
	\begin{equation*}
	\varphi\bm{h}^{(\alpha)}_n+(1-\varphi)\bm{g}^{(\alpha)}_n\to\varphi\bm{u}^{(\alpha)}+(1-\varphi)\bm{u}^{(\alpha)}=\bm{u}^{(\alpha)},\quad\textup{ in }\bm{H}^1(\Omega^{(\alpha)})\textup{ as }n\to\infty,
	\end{equation*}
	that is, $\bm{u}^{(\alpha)}_n\to\bm{u}^{(\alpha)}$ in $\bm{V}(\Omega^{(\alpha)})$ as $n\to\infty$.
	
	\underline{\emph{Sub-step 2: Construction of a separating control.}} Since $\bm{u}^{(\alpha)}_n$ is continuous on the compact set $\overline{\Omega^{(\alpha)}}$ and the gap on $\omega$ is at least $\delta/2$, continuity upgrades the almost-everywhere gap to every vertical coordinate. We use here that each body is a graph over $\omega$ in the $\bm{e}_3$ direction, so that every fibre $\overline{F^{(\alpha)}(\xi)}$ is an interval and the envelopes below are continuous. For each $\xi\in\overline{\omega}$, define the envelopes:
	\begin{equation*}
	\begin{aligned}
		L_n(\xi)&:=\min_{x_3\in\overline{F^{(1)}(\xi)}}\left\{\left(\mathbb{I}_{\Omega^{(1)}}(\xi_1,\xi_2,x_3)+\bm{u}^{(1)}_n(\xi_1,\xi_2,x_3)\right)\cdot\bm{e}_3\right\},\\
		U_n(\xi)&:=\max_{x_3\in\overline{F^{(2)}(\xi)}}\left\{\left(\mathbb{I}_{\Omega^{(2)}}(\xi_1,\xi_2,x_3)+\bm{u}^{(2)}_n(\xi_1,\xi_2,x_3)\right)\cdot\bm{e}_3\right\}.
	\end{aligned}
	\end{equation*}
	
	Thanks to~\eqref{eq:graph}, the mapping $\xi\mapsto\overline{F^{(\alpha)}(\xi)}$ has non-empty compact values and is continuous, and the height functions between curly brackets are continuous, so $L_n$ and $U_n$ are continuous by Berge's Maximum Theorem (cf., e.g., page~116 in~\cite{Ber63}) with $X=\overline{\pi_\kappa(\Omega^{(\alpha)})}$, $Y=\mathbb{R}$, $\Gamma \xi=\overline{F^{(\alpha)}(\xi)}$. Moreover, it results that:
	\begin{equation*}
		L_n(\xi)-U_n(\xi)\ge\dfrac{\delta}{2}>0,\qquad\xi\in\overline{\omega}.
	\end{equation*}
	
	Since the gap is positive, the function $(L_n+U_n)/2$ lies strictly between the two bodies on $\overline{\omega}$. On $\pi_\kappa(\Omega^{(1)})\setminus\omega$ only $\Omega^{(1)}$ is present, and on $\pi_\kappa(\Omega^{(2)})\setminus\omega$ only $\Omega^{(2)}$ is present, so each side region carries only one one-sided inequality. This gives enough room to extend $(L_n+U_n)/2$ continuously across $\partial\omega$ so as to obtain:
	\begin{equation*}
	\begin{aligned}
		\psi_n=\frac{L_n+U_n}{2}&,\quad\textup{ in }\overline{\omega},\\
		(\mathbb{I}_{\Omega^{(2)}}(x^{(2)})+\bm{u}^{(2)}_n(x^{(2)}))\cdot\bm{e}_3<\psi_n(x^{(2)}_1,x^{(2)}_2)&,\quad\textup{ on }\pi_\kappa(\Omega^{(2)})\setminus\omega,\\
		(\mathbb{I}_{\Omega^{(1)}}(x^{(1)})+\bm{u}^{(1)}_n(x^{(1)}))\cdot\bm{e}_3>\psi_n(x^{(1)}_1,x^{(1)}_2)&,\quad\textup{ on }\pi_\kappa(\Omega^{(1)})\setminus\omega.
	\end{aligned}
	\end{equation*}
	
	By~\eqref{kappa:compact}, it is ensured that $\pi_\kappa(\Omega^{(1)})\cup\pi_\kappa(\Omega^{(2)})\subset\subset\kappa$, so there exists a tubular neighbourhood of $\partial\kappa$ in which constraints are inactive. Cutting the continuous extension off smoothly there yields $\psi_n\in\mathcal{D}(\kappa)$ satisfying the required conditions. Hence, it results that $(\bm{u}^{(1)}_n,\bm{u}^{(2)}_n)\in\vec{\bm{U}}^{(2)}$, and since the pair converges strongly to $(\bm{u}^{(1)},\bm{u}^{(2)})$ in $\bm{V}(\Omega^{(1)})\times\bm{V}(\Omega^{(2)})$ as $n\to\infty$, the set $\vec{\bm{U}}^{(2)}$ is dense in $\vec{\bm{U}}^{(1)}$. This proves $(ii)$.
		
	$(iii)$ \emph{Let $\{\beta_k\}_{k=1}^\infty$ be an arbitrary sequence of positive real numbers such that $\beta_k\to 0$ as $k\to\infty$. There exists a sequence $\{n_k\}_{k=1}^\infty$ such that $n_k\to\infty$ as $k\to\infty$ and $(\beta_k\|\Delta\psi_{n_k}\|_{L^2(\kappa)}^2)\to 0$ as $k\to\infty$.} In what follows, to simplify the notation, we let $C_n:=\|\Delta\psi_n\|_{L^2(\kappa)}$. For each integer $k\ge 1$, define the set $M_k:=\{m\in\mathbb{N};C_m>\beta_k^{-1/4}\}$. If $k$ is such that $M_k=\emptyset$ then $C_m\le\beta_k^{-1/4}$ for all $m\in\mathbb{N}$ and thus, in particular, it results $C_k\le\beta_k^{-1/4}$ and we set $n_k:=k$.
	If, instead, the index $k$ is such that $M_k\neq\emptyset$, we let $n_k:=\max\{1,(\min M_k)-1\}$.
	
	We now show that the sequence $\{n_k\}_{k=1}^\infty$ that we defined by
	\begin{equation*}
		n_k:=
		\begin{cases}
			k&,\textup{ if }M_k=\emptyset,\\
			\max\{1,(\min M_k)-1\}&,\textup{ if }M_k\neq\emptyset,
		\end{cases}
	\end{equation*}
	tends to infinity as $k\to\infty$. Assume by contradiction that the sequence $\{n_k\}_{k=1}^\infty$ is bounded independently of $k$ or, equivalently, that there exists an integer $N\ge 1$ such that $n_k\le N$, for all $k\ge 1$. This implies that there are at most finitely many indices $k_1,\dots,k_r$ such that $M_{k_j}=\emptyset$ for all $1\le j\le r$. For each of the remaining indices $k$, for which the sets $M_k\neq\emptyset$, it results $\min M_k\le N+1$. This in turn implies that, for any integer $k\ge 1$, the index at which $\min M_k$ is attained cannot be larger than $N+1$. Therefore, it results that for each $k\in\mathbb{N}$ such that $M_k\neq\emptyset$
	\begin{equation*}
		\infty>\max\{C_1,\dots,C_{N+1}\}\ge C_{\min M_k}>\beta_k^{-1/4}\to\infty,\quad\textup{ as }k\to\infty,
	\end{equation*}
	providing the required contradiction.
	
	Let us now show that $\beta_k C_{n_k}^2\to 0$ as $k\to\infty$. If $M_k=\emptyset$, then $n_k=k\notin M_k$ and $C_{n_k}=C_k\le\beta_k^{-1/4}$.
	If $M_k\neq\emptyset$ and $\min M_k\ge2$, then $n_k=(\min M_k)-1\notin M_k$, so again $C_{n_k}\le\beta_k^{-1/4}$. The case $\min M_k=1$ can occur only when $C_1>\beta_k^{-1/4}$, that is, for finitely many instances of $k$. Therefore, for all $k$ sufficiently large, it results $C_{n_k}\le\beta_k^{-1/4}$, whence
	\begin{equation*}
		0\le\beta_k C_{n_k}^2\le\beta_k^{1/2}\to0,\quad\textup{ as }k\to\infty,
	\end{equation*}
	and the proof of $(iii)$ is complete.
	
	$(iv)$ \emph{Construction of the recovery sequence and completion of the proof.} Let $(\bm{u}^{(1)},\bm{u}^{(2)})\in\vec{\bm{U}}$ be given. Let $\{\beta_k\}_{k=1}^\infty$ be a sequence of positive real numbers such that $\beta_k\to 0$ as $k\to\infty$. 
	
	Since $\vec{\bm{U}}^{(1)}$ is dense in $\vec{\bm{U}}$ by $(i)$, there exists a sequence $\{(\bm{v}^{(1)}_m,\bm{v}^{(2)}_m)\}_{m=1}^\infty\subset\vec{\bm{U}}^{(1)}$ such that:
	\begin{equation*}
		(\bm{v}^{(1)}_m,\bm{v}^{(2)}_m)\to(\bm{u}^{(1)},\bm{u}^{(2)}),\quad\textup{ in }\bm{V}(\Omega^{(1)})\times\bm{V}(\Omega^{(2)})\textup{ as }m\to\infty.
	\end{equation*}
	
	For each integer $m\ge1$, item $(ii)$ applied to the element $(\bm{v}^{(1)}_m,\bm{v}^{(2)}_m)$ yields a sequence $\{(\bm{u}^{(1)}_{m,n},\bm{u}^{(2)}_{m,n})\}_{n=1}^\infty\subset\vec{\bm{U}}^{(2)}$ and a separating interface $\psi_{m,n}\in\mathcal{D}(\kappa)$ such that:
	\begin{equation*}
		(\bm{u}^{(1)}_{m,n},\bm{u}^{(2)}_{m,n})\to(\bm{v}^{(1)}_m,\bm{v}^{(2)}_m),\quad\textup{ in }\bm{V}(\Omega^{(1)})\times\bm{V}(\Omega^{(2)})\textup{ as }n\to\infty.
	\end{equation*}
	
	For each $m\ge1$ choose $n(m)\ge m$ such that
	\begin{equation*}
		\|(\bm{u}^{(1)}_{m,n(m)},\bm{u}^{(2)}_{m,n(m)})-(\bm{v}^{(1)}_m,\bm{v}^{(2)}_m)\|_{\bm{V}(\Omega^{(1)})\times\bm{V}(\Omega^{(2)})}\le\dfrac{1}{m},
	\end{equation*}
	and set $C_m:=\|\Delta\psi_{m,n(m)}\|_{L^2(\kappa)}$. Applying $(iii)$ to the sequence $\{C_m\}_{m=1}^\infty$ yields a sequence $m_k\to\infty$ with $\beta_k C_{m_k}^2\to0$ as $k\to\infty$.
	
	For every $k\ge1$, let $n_k:=n(m_k)$ and observe that the function $\Psi_{m_k,n_k}:=\psi_{m_k,n_k}\circ\pi$ belongs to $\mathcal{V}(K)$. The regularity is indeed clear, since $\psi_{m_k,n_k}\in\mathcal{D}(\kappa)\subset H^2_0(\kappa)$. For what concerns the validity of~\eqref{bdry:feas}, observe that $\bm{u}^{(\alpha)}_{m_k,n_k}=\bm{0}$ on $\Gamma_0^{(\alpha)}$. Moreover, since the inequalities of item $(ii)$ hold for a.a. points of $\Omega^{(\alpha)}$, and both sides are continuous up to $\overline{\Omega^{(\alpha)}}$ it results that these inequalities hold pointwise on $\overline{\Omega^{(\alpha)}}$. Restricting them to $\Gamma_0^{(\alpha)}$, where $\bm{u}^{(\alpha)}_{m_k,n_k}=\bm{0}$, gives exactly~\eqref{bdry:feas} with non-strict inequalities. Therefore the infimum in~\eqref{marginal} is bounded above by:
	\begin{equation*}
		\dfrac{\beta_k}{2}\|\Delta\Psi_{m_k,n_k}\|_{L^2(K)}^2=a\beta_k\|\Delta\psi_{m_k,n_k}\|_{L^2(\kappa)}^2=a\beta_k C_{m_k}^2.
	\end{equation*}
	
	Since the infimum in~\eqref{marginal} is non-negative, it results:
	\begin{equation}
		\label{eq:sandwich}
		\sum_{\alpha=1}^2 E_\alpha(\bm{u}^{(\alpha)}_{m_k,n_k})\le\tilde{\mathcal{F}}_{\beta_k}(\bm{u}^{(1)}_{m_k,n_k},\bm{u}^{(2)}_{m_k,n_k})\le a\beta_k C_{m_k}^2+\sum_{\alpha=1}^2 E_\alpha(\bm{u}^{(\alpha)}_{m_k,n_k}).
	\end{equation}
	
	Since, by construction $(\bm{u}^{(1)}_{m_k,n_k},\bm{u}^{(2)}_{m_k,n_k})\to(\bm{u}^{(1)},\bm{u}^{(2)})$ in $\bm{V}(\Omega^{(1)})\times\bm{V}(\Omega^{(2)})$ as $k\to\infty$, the continuity of $E_\alpha$ allows us to infer the following for the left-hand side of~\eqref{eq:sandwich}:
	\begin{equation*}
		\sum_{\alpha=1}^2 E_\alpha(\bm{u}^{(\alpha)}_{m_k,n_k})\to\sum_{\alpha=1}^2 E_\alpha(\bm{u}^{(\alpha)})=\tilde{\mathcal{F}}_0(\bm{u}^{(1)},\bm{u}^{(2)}),\quad\textup{ as }k\to\infty.
	\end{equation*}
	
	Since $(a\beta_k C_{m_k}^2)\to0$ as $k\to\infty$, we obtain that the right-hand side of~\eqref{eq:sandwich} satisfies
	\begin{equation*}
		a\beta_k C_{m_k}^2+\sum_{\alpha=1}^2 E_\alpha(\bm{u}^{(\alpha)}_{m_k,n_k})\to\sum_{\alpha=1}^2 E_\alpha(\bm{u}^{(\alpha)})=\tilde{\mathcal{F}}_0(\bm{u}^{(1)},\bm{u}^{(2)}),\quad\textup{ as }k\to\infty,
	\end{equation*}
	and the proof is complete.
\end{proof}

We are then ready to establish the $\Gamma$-convergence of the sequence of functionals $\{\tilde{\mathcal{F}}_n\}_{n=1}^\infty$ to $\tilde{\mathcal{F}}_0$.
\begin{theorem}
	\label{gamma-4}
	Let $\{\beta_n\}_{n=1}^\infty$ be an arbitrary sequence of positive real numbers such that $\beta_n\to 0$ as $n\to\infty$.
	Then, the sequence of functionals $\{\tilde{\mathcal{F}}_{\beta_n}\}_{n=1}^\infty$ $\Gamma$-converges to $\tilde{\mathcal{F}}_0$ in the weak topology of
	$\bm{V}(\Omega^{(1)})\times\bm{V}(\Omega^{(2)})$.
	
	Additionally, if $(\bm{u}^{(1)}_n,\bm{u}^{(2)}_n)$ denotes the unique minimiser of $\tilde{\mathcal{F}}_{\beta_n}$ constructed in~\eqref{id-min} and $(\bm{u}^{(1)}_0,\bm{u}^{(2)}_0)$ denotes the unique minimiser of $\tilde{\mathcal{F}}_0$ constructed in Theorem~\ref{th-ex-un-0}, it results that $(\bm{u}^{(1)}_n,\bm{u}^{(2)}_n)\rightharpoonup(\bm{u}^{(1)}_0,\bm{u}^{(2)}_0)$ in $\bm{V}(\Omega^{(1)})\times\bm{V}(\Omega^{(2)})$ as $n\to\infty$.
\end{theorem}
\begin{proof}
	Thanks to Lemma~\ref{gamma-1}, Lemma~\ref{gamma-2} and Lemma~\ref{gamma-3}, we are in position to apply Theorem~11.18-1 in~\cite{Cia25} (see also~\cite{DalMaso93}) to assert the $\Gamma$-convergence of the sequence of functionals $\{\tilde{\mathcal{F}}_{\beta_n}\}_{n=1}^\infty$ to $\tilde{\mathcal{F}}_0$. The convergence $(\bm{u}^{(1)}_n,\bm{u}^{(2)}_n)\rightharpoonup(\bm{u}^{(1)}_0,\bm{u}^{(2)}_0)$ in $\bm{V}(\Omega^{(1)})\times\bm{V}(\Omega^{(2)})$ as $n\to\infty$ also follows from Theorem~11.18-1 in~\cite{Cia25}.
\end{proof}

The $\Gamma$-convergence established in Theorem~\ref{gamma-4} allows us to describe the coincidence set of Problem~\ref{problem0} in terms of the optimal controls. For $\xi=(\xi_1,\xi_2)\in\omega$ define the envelopes of the limit states
\begin{equation*}
	\begin{aligned}
		\ell_0(\xi)&:=\textup{essinf}_{x_3\in F^{(1)}(\xi)}\left\{\left(\mathbb{I}_{\Omega^{(1)}}(\xi_1,\xi_2,x_3)+\bm{u}^{(1)}_0(\xi_1,\xi_2,x_3)\right)\cdot\bm{e}_3\right\},\\
		m_0(\xi)&:=\textup{esssup}_{x_3\in F^{(2)}(\xi)}\left\{\left(\mathbb{I}_{\Omega^{(2)}}(\xi_1,\xi_2,x_3)+\bm{u}^{(2)}_0(\xi_1,\xi_2,x_3)\right)\cdot\bm{e}_3\right\},
	\end{aligned}
\end{equation*}
and define the coincidence set by:
\begin{equation*}
	\Lambda_0:=\left\{\xi\in\omega; \ell_0(\xi)=m_0(\xi)\right\}.
\end{equation*}

For every integer $n\ge1$, the confinement conditions satisfied by the optimal states give the following estimates
\begin{equation}
	\label{coincidence-sandwich}
	[\mathbb{I}_{\Omega^{(2)}}(x^{(2)})+\bm{u}^{(2)}(\Psi_{\beta_n})(x^{(2)})]\cdot\bm{e}_3\le\Psi_{\beta_n}(\xi_1,\xi_2)\le[\mathbb{I}_{\Omega^{(1)}}(x^{(1)})+\bm{u}^{(1)}(\Psi_{\beta_n})(x^{(1)})]\cdot\bm{e}_3,
\end{equation}
for a.a. $x^{(1)}=(x_1,x_2,x_3^{(1)})\in\Omega^{(1)}$ with $(x_1,x_2)=\xi$, and a.a. $x^{(2)}=(x_1,x_2,x_3^{(2)})\in\Omega^{(2)}$ with $(x_1,x_2)=\xi$.

By Theorem~\ref{gamma-4}, the optimal states converge weakly in $\bm{V}(\Omega^{(1)})\times\bm{V}(\Omega^{(2)})$ to $(\bm{u}^{(1)}_0,\bm{u}^{(2)}_0)$. The compact embedding $\bm{H}^1(\Omega^{(\alpha)})\hookrightarrow\hookrightarrow\bm{L}^2(\Omega^{(\alpha)})$ allows us to apply the Rellich-Kondra\v{s}ov theorem (cf., e.g., Theorem~8.4-3 in~\cite{Cia25}) so as to infer that, up to passing to a subsequence that we do not relabel, it results
\begin{equation*}
	[\mathbb{I}_{\Omega^{(\alpha)}}+\bm{u}^{(\alpha)}(\Psi_{\beta_n})]\cdot\bm{e}_3\to[\mathbb{I}_{\Omega^{(\alpha)}}+\bm{u}^{(\alpha)}_0]\cdot\bm{e}_3,\quad\textup{ as }n\to\infty,
\end{equation*}
a.e. in $\Omega^{(\alpha)}$ for $\alpha\in\{1,2\}$. By Fubini's theorem (cf., e.g., Theorem~1.15-6$(b)$ in~\cite{Cia25}), for a.a. $\xi\in\omega$, the previous pointwise convergence holds on a.e. point of the two fibres $F^{(1)}(\xi)$ and $F^{(2)}(\xi)$.
Fix one such point $\xi\in\Lambda_0$ and let $y:=\ell_0(\xi)=m_0(\xi)$. For the $\xi$ under consideration, the left-hand side of~\eqref{coincidence-sandwich} converges to $m_0(\xi)=y$ as $n\to\infty$, and the right-hand side of~\eqref{coincidence-sandwich} converges to $\ell_0(\xi)=y$ as $n\to\infty$, and we thus obtain that:
\begin{equation*}
	\Psi_{\beta_n}(\xi_1,\xi_2)\to y,\quad \textup{ as }n\to\infty.
\end{equation*}

In conclusion, we have shown that, for $\mathrm{d}\mathcal{L}^2$-a.a. $\xi\in\Lambda_0$, the sequence $\{\Psi_{\beta_n}(\xi)\}_{n=1}^\infty$ converges to the contact height $y=\ell_0(\xi)=m_0(\xi)$ as $n\to\infty$.

We close this section with a comparison between the optimal control approach proposed here, and the classical penalty method.
\begin{remark}
\label{rem:comp}
To place the comparison on firm ground, we introduce the classical penalty method for Problem~\ref{problem0}. Recall that the admissible set $\vec{\bm{U}}$ is equivalently characterised by the single inequality $\Phi\ge0$ $\mathrm{d}\mu$-a.e. on $S$ (cf.~\eqref{eq:Phidef} and~\eqref{eq:Phi}). For each $\varepsilon>0$, the penalised problem reads
\begin{equation*}
	\min_{\substack{\bm{v}^{(1)}\in\bm{V}(\Omega^{(1)})\\\bm{v}^{(2)}\in\bm{V}(\Omega^{(2)})}}\mathcal{J}_\varepsilon(\bm{v}^{(1)},\bm{v}^{(2)}),
\end{equation*}
where
\begin{equation}
	\label{penalised-functional}
	\mathcal{J}_\varepsilon(\bm{v}^{(1)},\bm{v}^{(2)}):=\sum_{\alpha=1}^{2}E_\alpha(\bm{v}^{(\alpha)})+\frac{1}{2\varepsilon}\int_{S}\left|\left\{\Phi(\bm{v}^{(1)},\bm{v}^{(2)})\right\}^{-}\right|^2\,\mathrm{d}\mu.
\end{equation}

The penalised functional is strictly convex and coercive, the penalty term being convex and non-negative and the bulk energy being elliptic by~\eqref{eq:pos-def} and Theorem~\ref{KornCart}. Therefore, the functional $\mathcal{J}_\varepsilon$ admits a unique minimiser $(\bm{u}^{(1)}_\varepsilon,\bm{u}^{(2)}_\varepsilon)$, which converges to the solution of Problem~\ref{problem0} as $\varepsilon\to0$ in the weak topology of $\bm{V}(\Omega^{(1)})\times\bm{V}(\Omega^{(2)})$.

\emph{While the penalty method restores the non-penetration constraint in the limit as $\varepsilon\to0$, its numerical treatment is non-standard.}
This is due to the fact that the penalty term in~\eqref{penalised-functional} couples the two displacement fields through the four-dimensional fibre product $S$. Its discretisation requires the simultaneous assembly of both elasticity systems together with a fibre-wise integration over $S$ with respect to the product measure $\mathcal{L}^2\otimes\mathcal{L}^1\otimes\mathcal{L}^1$, and the conditioning of the resulting discrete system is expected to deteriorate as $\varepsilon\to0$.

By contrast, Problem~\ref{controlproblem} removes the coupling between the two bodies. For each fixed control $\Psi\in\mathcal{V}(K)$, the state problems~\ref{problem1} and~\ref{problem2} are two \emph{independent} single-body unilateral problems with a \emph{prescribed} obstacle, which is precisely the setting for which classical and robust numerical methods are available. The non-penetration requirement is enforced by the control itself, and the interface is now the primary unknown of this optimal control problem. The computational cost of this decoupling is the outer minimisation over the control, while its consistency as $\beta\to0$ is established via $\Gamma$-convergence (Theorem~\ref{gamma-4}).
In this sense, the optimal control formulation trades a six-field, tightly coupled and increasingly ill-conditioned non-linear system for a sequence of decoupled, well-understood single-body obstacle sub-problems.
\bqed
\end{remark}

\section{Existence and Regularity of Optimal Interfaces. Part~I}\label{Sec:3bis}

For the sake of not rendering the notation too heavy, in what follows we drop the dependence on $\beta$ in the cost functional considered in Problem~\ref{controlproblem} as from now onward $\beta$ is a fixed parameter.
To begin with, we relax the variational problems associated with the states, beginning with Problem~\ref{problem1}. Define the operator $\gls{beta1}:\bm{V}(\Omega^{(1)})\to\bm{V}(\Omega^{(1)})$ by:
\begin{equation*}
	\bm{\beta}_1(\bm{v};\Psi):=-\{[\mathbb{I}_{\Omega^{(1)}}+\bm{v}-\Psi\bm{e}_3]\cdot\bm{e}_3\}^{-}\bm{e}_3,\quad\textup{ for all }\bm{v}\in\bm{V}(\Omega^{(1)}).
\end{equation*}

In a similar fashion as Theorem~5.8-2 in~\cite{ABM14}, one can observe that the operator $\bm{\beta}_1(\cdot;\Psi)$ is Lipschitz continuous in $\bm{L}^2(\Omega^{(1)})$ with Lipschitz constant equal to one, and monotone.

In what follows, we denote by $\gls{rho}>0$ a penalty parameter which is meant to approach zero. We are now ready to state the relaxation of Problem~\ref{problem1}.
\begin{customprob}{$\mathcal{P}_\rho(\Omega^{(1)};\Psi)$}\label{problem1rho}
	Find $\bm{u}_\rho^{(1)}=(u_{\rho,i}^{(1)})\in\bm{V}(\Omega^{(1)})$ that satisfies the variational equations
	\begin{equation*}
		\int_{\Omega^{(1)}} A^{ijk\ell,1} e_{k\|\ell}(\bm{u}_\rho^{(1)}) e_{i\|j}(\bm{v}) \,\mathrm{d}x^{(1)}+\dfrac{1}{\rho}\int_{\Omega^{(1)}}\bm{\beta}_1(\bm{u}_\rho^{(1)};\Psi)\cdot\bm{v}\,\mathrm{d}x^{(1)} = \int_{\Omega^{(1)}} \bm{f}^{(1)} \cdot \bm{v} \,\mathrm{d}x^{(1)},
	\end{equation*}
	for all $\bm{v}=(v_i) \in \bm{V}(\Omega^{(1)})$.
	\bqed
\end{customprob}

The existence and uniqueness of solutions for Problem~\ref{problem1rho} follows from an application of Korn's inequality (Theorem~\ref{KornCart}) and of the Minty--Browder theorem (cf., e.g., Theorem~12.5-3 in~\cite{Cia25}), being the operator $\bm{\beta}_1(\cdot;\Psi)$ Lipschitz continuous and monotone.

\begin{theorem}
	\label{th:ex1rho}
	Let $\rho>0$ and $\bm{f}^{(1)}\in\bm{L}^2(\Omega^{(1)})$ be given.
	Then, Problem~\ref{problem1rho} admits a unique solution.
\end{theorem}
\begin{proof}
	Define the operator $\mathcal{A}_\rho^{(1)}:\bm{V}(\Omega^{(1)})\to\bm{V}^\ast(\Omega^{(1)})$ in a way such that:
	\begin{equation*}
		\langle \mathcal{A}_\rho^{(1)}\bm{u},\bm{v}\rangle_{\bm{V}^\ast(\Omega^{(1)}),\bm{V}(\Omega^{(1)})}=\int_{\Omega^{(1)}}A^{ijk\ell,1} e_{k\|\ell}(\bm{u})e_{i\|j}(\bm{v})\,\mathrm{d}x^{(1)} +\dfrac{1}{\rho}\int_{\Omega^{(1)}}\bm{\beta}_1(\bm{u};\Psi)\cdot\bm{v}\,\mathrm{d}x^{(1)},\quad\textup{ for all }\bm{u},\bm{v}\in\bm{V}(\Omega^{(1)}).
	\end{equation*}
	
	The operator $\mathcal{A}_\rho^{(1)}$ is continuous by H\"older's inequality (cf., e.g., \cite{Brez11}) and the boundedness of the components $A^{ijk\ell,1}$ of the fourth order three-dimensional elasticity tensor associated with the linearly elastic body $\Omega^{(1)}$. Besides, an application of the monotonicity of the operator $\bm{\beta}_1(\cdot;\Psi)$ and an application of Korn's inequality (Theorem~\ref{KornCart}) give that the operator $\mathcal{A}_\rho^{(1)}$ is also monotone and coercive. Therefore, the claimed existence and uniqueness follow from the Minty--Browder theorem (cf., e.g., Theorem~12.5-3 in~\cite{Cia25}).
\end{proof}

In order to relax Problem~\ref{problem2}, define the operator $\gls{beta2}:\bm{V}(\Omega^{(2)})\to\bm{V}(\Omega^{(2)})$ by:
\begin{equation*}
	\bm{\beta}_2(\bm{v};\Psi):=\{[\mathbb{I}_{\Omega^{(2)}}+\bm{v}-\Psi\bm{e}_3]\cdot\bm{e}_3\}^{+}\bm{e}_3,\quad\textup{ for all }\bm{v}\in\bm{V}(\Omega^{(2)}).
\end{equation*}

It is straightforward to observe that the operator $\bm{\beta}_2(\cdot;\Psi)$ is Lipschitz continuous in $\bm{L}^2(\Omega^{(2)})$ with Lipschitz constant equal to one, and monotone. We are now ready to state the relaxation of Problem~\ref{problem2}.

\begin{customprob}{$\mathcal{P}_\rho(\Omega^{(2)};\Psi)$}\label{problem2rho}
	Find $\bm{u}_\rho^{(2)}=(u_{\rho,i}^{(2)})\in\bm{V}(\Omega^{(2)})$ that satisfies the variational equations
	\begin{equation*}
		\int_{\Omega^{(2)}} A^{ijk\ell,2} e_{k\|\ell}(\bm{u}_\rho^{(2)}) e_{i\|j}(\bm{v}) \,\mathrm{d}x^{(2)}+\dfrac{1}{\rho}\int_{\Omega^{(2)}}\bm{\beta}_2(\bm{u}_\rho^{(2)};\Psi)\cdot\bm{v}\,\mathrm{d}x^{(2)} = \int_{\Omega^{(2)}} \bm{f}^{(2)} \cdot \bm{v} \,\mathrm{d}x^{(2)},
	\end{equation*}
	for all $\bm{v}=(v_i) \in \bm{V}(\Omega^{(2)})$.
	\bqed
\end{customprob}

The existence and uniqueness of solutions for Problem~\ref{problem2rho} follows from an application of Korn's inequality (Theorem~\ref{KornCart}) and of the Minty--Browder theorem (cf., e.g., Theorem~12.5-3 in~\cite{Cia25}), being the operator $\bm{\beta}_2(\cdot;\Psi)$ Lipschitz continuous and monotone. The proof is identical to the one of Theorem~\ref{th:ex1rho} and we thus omit it. We are thus in a position to introduce, for each $\rho>0$, a relaxation of Problem~\ref{controlproblem}, which we denote by $\mathscr{C}_\rho(K)$ in what follows.

\begin{customprob}{$\mathscr{C}_\rho(K)$}\label{cprho}
	Find $\Psi_\rho\in\mathcal{V}(K)$ that minimises the reduced cost functional
	\begin{equation*}
		\begin{aligned}
			&\gls{Jrho}(\Psi):=\dfrac{\beta}{2}\int_{K}|\Delta\Psi|^2\,\mathrm{d}x
			+\dfrac{1}{2}\int_{\Omega^{(1)}}A^{ijk\ell,1}e_{k\|\ell}(\bm{u}_\rho^{(1)}(\Psi))e_{i\|j}(\bm{u}_\rho^{(1)}(\Psi))\,\mathrm{d}x+\dfrac{1}{2\rho}\int_{\Omega^{(1)}}|\bm{\beta}_1(\bm{u}_\rho^{(1)}(\Psi);\Psi)|^2\,\mathrm{d}x\\
			&\qquad-\int_{\Omega^{(1)}}\bm{f}^{(1)}\cdot\bm{u}_\rho^{(1)}(\Psi)\,\mathrm{d}x\\
			&\quad+\dfrac{1}{2}\int_{\Omega^{(2)}}A^{ijk\ell,2}e_{k\|\ell}(\bm{u}_\rho^{(2)}(\Psi))e_{i\|j}(\bm{u}_\rho^{(2)}(\Psi))\,\mathrm{d}x+\dfrac{1}{2\rho}\int_{\Omega^{(2)}}|\bm{\beta}_2(\bm{u}_\rho^{(2)}(\Psi);\Psi)|^2\,\mathrm{d}x-\int_{\Omega^{(2)}}\bm{f}^{(2)}\cdot\bm{u}_\rho^{(2)}(\Psi)\,\mathrm{d}x,
		\end{aligned}
	\end{equation*}
	subject to $\bm{u}_\rho^{(1)}(\Psi)$ being the solution of Problem~\ref{problem1rho} and to $\bm{u}_\rho^{(2)}(\Psi)$ being the solution of Problem~\ref{problem2rho}.
	\bqed	
\end{customprob}

Observe that the penalty term in Problem~\ref{cprho} makes the penalised elastic energy appearing in the expression of $\mathcal{J}_\rho$ a convex function of $\Psi$ (cf., e.g., Section~3.2.5 in~\cite{BV04}). Combined with the strict convexity of the term $\frac{\beta}{2}\|\Delta\Psi\|^2_{L^2(K)}$, this yields the strict convexity of the reduced cost functional $\mathcal{J}_\rho$ (Theorem~\ref{ex-un-op-con-rho}).

By contrast, if only the elastic energy $E_\alpha(\bm{u}_\rho^{(\alpha)}(\Psi))$ were kept in the expression of the cost functional, convexity with respect to $\Psi$ could no longer be guaranteed. This is due to the fact that the control-to-state mapping $\Psi\mapsto\bm{u}_\rho^{(\alpha)}(\Psi)$ is in general non-linear and, therefore, the composition of the convex functional $E_\alpha$ with this mapping may not be convex.

The existence and uniqueness of solutions for Problem~\ref{cprho} hinges on two preparatory lemmas. The first lemma consists in establishing the continuity of the control-to-state operator for Problem~\ref{controlproblem}. We observe that this fact cannot be established directly since the function $\Psi$ associated with the obstacle enters the formulation of Problem~\ref{problem1} and Problem~\ref{problem2} only in the definition of the sets $\bm{K}(\Omega^{(1)};\Psi)$ and $\bm{K}(\Omega^{(2)};\Psi)$ of admissible competitors. In order to establish the continuity of the control-to-state operator, we first show that for a given $\Psi\in\mathcal{V}(K)$ the family of solutions $\bm{u}_\rho^{(1)}(\Psi)$ of Problem~\ref{problem1rho} converges to the solution $\bm{u}^{(1)}(\Psi)$ of Problem~\ref{problem1} in $\bm{V}(\Omega^{(1)})$ as $\rho\to0$.
\begin{lemma}
	\label{lem:conv}
	Let $\bm{f}^{(\alpha)}\in\bm{L}^2(\Omega^{(\alpha)})$ and let $\Psi\in\mathcal{V}(K)$ be given.
	For each $\rho>0$, let $\bm{u}_\rho^{(\alpha)}(\Psi)$ be the unique solution of Problem~\ref{problem1rho} when $\alpha=1$ and of Problem~\ref{problem2rho} when $\alpha=2$. Then, it results that $\bm{u}_\rho^{(\alpha)}(\Psi)\to\bm{u}^{(\alpha)}(\Psi)$ in $\bm{V}(\Omega^{(\alpha)})$ as $\rho\to0$, where $\bm{u}^{(\alpha)}(\Psi)$ is the unique solution of Problem~\ref{problem1}when $\alpha=1$ and of Problem~\ref{problem2} when $\alpha=2$.
\end{lemma}
\begin{proof}
	We will focus on establishing the sought convergence for $\alpha=1$ as the proof for the case $\alpha=2$ is identical.
	The proof is divided into three parts, numbered $(i)$--$(iii)$.
	
	$(i)$ \emph{The family $\{\bm{u}_\rho^{(1)}(\Psi)\}_{\rho>0}$ is bounded in $\bm{V}(\Omega^{(1)})$ independently of $\rho$, and the family $\{\rho^{-1}\bm{\beta}_1(\bm{u}_\rho^{(1)}(\Psi);\Psi)\}_{\rho>0}$ is bounded in $\bm{V}^\ast(\Omega^{(1)})$ independently of $\rho$.} Define the vector field
	\begin{equation*}
		\bm{w}_0^{(1)}(\Psi):=\left\{(\mathbb{I}_{\Omega^{(1)}}-\Psi\bm{e}_3)\cdot\bm{e}_3\right\}^{-}\bm{e}_3,
	\end{equation*}
	and observe that, thanks to Stampacchia's theorem (cf., e.g., Theorem~5.8-2 in~\cite{ABM14}) and~\eqref{bdry:feas}, it results $\bm{w}_0^{(1)}(\Psi)\in\bm{V}(\Omega^{(1)})$. Additionally, observe that
	\begin{equation*}
		\bm{\beta}_1(\bm{w}_0^{(1)}(\Psi);\Psi)=-\left\{(\mathbb{I}_{\Omega^{(1)}}+\bm{w}_0^{(1)}(\Psi)-\Psi\bm{e}_3)\cdot\bm{e}_3\right\}^{-}\bm{e}_3
		=-\left\{\left\{(\mathbb{I}_{\Omega^{(1)}}-\Psi\bm{e}_3)\cdot\bm{e}_3\right\}^{+}\right\}^{-}\bm{e}_3=\bm{0},
	\end{equation*}
	so that $\bm{w}_0^{(1)}(\Psi)\in\bm{K}(\Omega^{(1)};\Psi)$. Specialise $\bm{v}=(\bm{u}_\rho^{(1)}(\Psi)-\bm{w}_0^{(1)}(\Psi))$ in the variational equations of Problem~\ref{problem1rho}. We first observe that:
	\begin{equation}
		\label{est:penalty}
		\begin{aligned}
			&\int_{\Omega^{(1)}}\bm{\beta}_1(\bm{u}_\rho^{(1)}(\Psi);\Psi)\cdot(\bm{u}_\rho^{(1)}(\Psi)-\bm{w}_0^{(1)}(\Psi))\,\mathrm{d}x^{(1)}\\
			&=\int_{\Omega^{(1)}}\left\{(\mathbb{I}_{\Omega^{(1)}}+\bm{u}_\rho^{(1)}(\Psi)-\Psi\bm{e}_3)\cdot\bm{e}_3\right\}^{-}\left\{(\mathbb{I}_{\Omega^{(1)}}-\Psi\bm{e}_3)\cdot\bm{e}_3\right\}^{+}\,\mathrm{d}x^{(1)}\\
			&\quad-\int_{\Omega^{(1)}}\left\{(\mathbb{I}_{\Omega^{(1)}}+\bm{u}_\rho^{(1)}(\Psi)-\Psi\bm{e}_3)\cdot\bm{e}_3\right\}^{-}\left((\mathbb{I}_{\Omega^{(1)}}+\bm{u}_\rho^{(1)}(\Psi)-\Psi\bm{e}_3)\cdot\bm{e}_3\right)\,\mathrm{d}x^{(1)}\\
			&\ge\int_{\Omega^{(1)}}|\bm{\beta}_1(\bm{u}_\rho^{(1)}(\Psi);\Psi)|^2\,\mathrm{d}x^{(1)}.
		\end{aligned}
	\end{equation}
	
	Second, we obtain that an application of the uniform positive-definiteness of the fourth order three-dimensional tensor $\{A^{ijk\ell,1}\}$ (viz. \eqref{eq:pos-def}), Korn's inequality (Theorem~\ref{KornCart}), Young's inequality~\cite{Young1912} and~\eqref{est:penalty} give
	\begin{equation*}
		\begin{aligned}
			&\dfrac{1}{2C_e^{(1)}|c_0^{(1)}|^2}\|\bm{u}_\rho^{(1)}(\Psi)\|_{\bm{V}(\Omega^{(1)})}^2+\dfrac{1}{\rho}\int_{\Omega^{(1)}}|\bm{\beta}_1(\bm{u}_\rho^{(1)}(\Psi);\Psi)|^2\,\mathrm{d}x^{(1)}\\
			&\le\|\bm{f}^{(1)}\|_{\bm{L}^2(\Omega^{(1)})}\|\bm{w}_0^{(1)}(\Psi)\|_{\bm{L}^2(\Omega^{(1)})}\\
			&\quad+C_e^{(1)}|c_0^{(1)}|^2\left(\|\bm{f}^{(1)}\|_{\bm{L}^2(\Omega^{(1)})}^2+\max_{i,j,k,\ell}\|A^{ijk\ell,1}\|_{\mathcal{C}^0(\overline{\Omega^{(1)}})}^2\|\bm{w}_0^{(1)}(\Psi)\|_{\bm{V}(\Omega^{(1)})}^2\right)=:D^{(1)}(\Psi,\bm{f}^{(1)}),
		\end{aligned}
	\end{equation*}
	from which the following estimates follow
	\begin{equation}
		\label{eq:bdd:sol}
		\begin{aligned}
			\|\bm{u}_\rho^{(1)}(\Psi)\|_{\bm{V}(\Omega^{(1)})}&\le\sqrt{2C_e^{(1)}|c_0^{(1)}|^2 D^{(1)}(\Psi,\bm{f}^{(1)})},\\
			\|\bm{\beta}_1(\bm{u}_\rho^{(1)}(\Psi);\Psi)\|_{\bm{L}^2(\Omega^{(1)})}&\le\sqrt{\rho D^{(1)}(\Psi,\bm{f}^{(1)})},
		\end{aligned}
	\end{equation}
	thus showing that the family of solutions $\{\bm{u}_\rho^{(1)}\}_{\rho>0}$ is bounded in $\bm{V}(\Omega^{(1)})$ independently of $\rho$.
	Therefore, the variational equations of Problem~\ref{problem1rho} in turn give:
	\begin{equation}
		\label{eq:bdd:penalty}
		\dfrac{1}{\rho}\left\|\bm{\beta}_1(\bm{u}_\rho^{(1)}(\Psi);\Psi)\right\|_{\bm{V}^\ast(\Omega^{(1)})}\le\|\bm{f}^{(1)}\|_{\bm{L}^2(\Omega^{(1)})}+\sup_{\substack{\bm{v}\in\bm{V}(\Omega^{(1)})\\\|\bm{v}\|_{\bm{V}(\Omega^{(1)})}=1}}\left|\int_{\Omega^{(1)}}A^{ijk\ell,1}e_{k\|\ell}(\bm{u}_\rho^{(1)}(\Psi))e_{i\|j}(\bm{v})\,\mathrm{d}x^{(1)}\right|.
	\end{equation}
	
	Therefore, we obtain that an application of~\eqref{eq:bdd:sol} gives that the left-hand side of~\eqref{eq:bdd:penalty} is bounded independently of $\rho$, as it was to be proved.
	
	$(ii)$ \emph{Every sequence $\{\rho_n\}_{n=1}^\infty$ such that $\rho_n\to0$ admits a subsequence, whose iterates are still denoted by $\rho_n$, such that $\bm{u}_{\rho_n}^{(1)}(\Psi)\rightharpoonup\bm{u}^{(1)}(\Psi)$ in $\bm{V}(\Omega^{(1)})$ as $n\to\infty$, where $\bm{u}^{(1)}(\Psi)$ is the unique solution of Problem~\ref{problem1}.} An application of~\eqref{eq:bdd:sol} gives that, up to passing to a subsequence $\{\rho_n\}_{n=1}^\infty$ such that $\rho_n\to 0$ as $n\to\infty$, it results that $\bm{u}_{\rho_n}^{(1)}(\Psi)\rightharpoonup\bm{u}^{(1)}$ in $\bm{V}(\Omega^{(1)})$ as $n\to\infty$, and that $\bm{u}^{(1)}\in\bm{K}(\Omega^{(1)};\Psi)$.
	
	Testing the variational equations of Problem~\ref{problem1rho} at $(\bm{v}-\bm{u}_{\rho_n}^{(1)}(\Psi))$, where $\bm{v}$ is arbitrarily chosen in $\bm{K}(\Omega^{(1)};\Psi)$, gives:
	\begin{equation}
		\label{eq:pen}
		\begin{aligned}
			&\int_{\Omega^{(1)}}A^{ijk\ell,1}e_{k\|\ell}(\bm{u}_{\rho_n}^{(1)}(\Psi))e_{i\|j}(\bm{v}-\bm{u}_{\rho_n}^{(1)}(\Psi))\,\mathrm{d}x^{(1)}+\dfrac{1}{\rho_n}\int_{\Omega^{(1)}}\bm{\beta}_1(\bm{u}_{\rho_n}^{(1)}(\Psi);\Psi)\cdot(\bm{v}-\bm{u}_{\rho_n}^{(1)}(\Psi))\,\mathrm{d}x^{(1)}\\
			&=\int_{\Omega^{(1)}}\bm{f}^{(1)}\cdot(\bm{v}-\bm{u}_{\rho_n}^{(1)}(\Psi))\,\mathrm{d}x^{(1)}.
		\end{aligned}
	\end{equation}
	
	Exploiting the monotonicity of $\bm{\beta}_1(\cdot;\Psi)$, the fact that $\bm{v}\in\bm{K}(\Omega^{(1)};\Psi)$, and the continuity of the operator $\mathcal{A}$, letting $n\to\infty$ in~\eqref{eq:pen} gives:
	\begin{equation}
		\label{eq:orig}
		\int_{\Omega^{(1)}}A^{ijk\ell,1}e_{k\|\ell}(\bm{u}^{(1)})e_{i\|j}(\bm{v}-\bm{u}^{(1)})\,\mathrm{d}x^{(1)} \ge \int_{\Omega^{(1)}}\bm{f}^{(1)}\cdot(\bm{v}-\bm{u}^{(1)})\,\mathrm{d}x^{(1)},\quad\textup{ for all }\bm{v}\in\bm{K}(\Omega^{(1)};\Psi).
	\end{equation}
	
	Since~\eqref{eq:orig} admits a unique solution (Theorem~\ref{th:ex1}), it must be $\bm{u}^{(1)}=\bm{u}^{(1)}(\Psi)$.
	
	$(iii)$ \emph{Conclusion of the proof.}
	Consider the vector field $\bm{v}_\rho:=\bm{u}_\rho^{(1)}(\Psi)-\bm{\beta}_1(\bm{u}_\rho^{(1)}(\Psi);\Psi)$, and observe that, in light of Stampacchia's theorem (cf., e.g., Theorem~5.8.2 of~\cite{ABM14}) and~\eqref{bdry:feas}, it results $\bm{v}_\rho\in\bm{V}(\Omega^{(1)})$, being both $\bm{u}_\rho^{(1)}(\Psi)\in\bm{V}(\Omega^{(1)})$ and $\bm{\beta}_1(\bm{u}_\rho^{(1)}(\Psi);\Psi)\in\bm{V}(\Omega^{(1)})$. Besides, observe that
	\begin{equation*}
		\begin{aligned}
			&(\mathbb{I}_{\Omega^{(1)}}+\bm{v}_\rho)\cdot\bm{e}_3-\Psi
			=(\mathbb{I}_{\Omega^{(1)}}+\bm{u}_\rho^{(1)}(\Psi)-\bm{\beta}_1(\bm{u}_\rho^{(1)}(\Psi);\Psi))\cdot\bm{e}_3-\Psi\\
			&=(\mathbb{I}_{\Omega^{(1)}}+\bm{u}_\rho^{(1)}(\Psi)-\Psi\bm{e}_3)\cdot\bm{e}_3+\{[\mathbb{I}_{\Omega^{(1)}}+\bm{u}_\rho^{(1)}(\Psi)-\Psi\bm{e}_3]\cdot\bm{e}_3\}^{-}\\
			&=\{[\mathbb{I}_{\Omega^{(1)}}+\bm{u}_\rho^{(1)}(\Psi)-\Psi\bm{e}_3]\cdot\bm{e}_3\}^{+}\ge0,\quad\textup{ a.e. in }\Omega^{(1)},
		\end{aligned}
	\end{equation*}
	which in turn implies that $\bm{v}_\rho\in\bm{K}(\Omega^{(1)};\Psi)$. Let $\bm{u}^{(1)}$ denote the unique solution of Problem~\ref{problem1}. Subtracting the variational inequalities governing Problem~\ref{problem1} tested at $\bm{v}_\rho$ from the variational equations governing Problem~\ref{problem1rho} tested at $(\bm{v}_\rho-\bm{u}^{(1)})$ gives:
	\begin{equation}
		\label{eq:1}
		\begin{aligned}
			&\int_{\Omega^{(1)}}A^{ijk\ell,1}e_{k\|\ell}(\bm{u}_\rho^{(1)}(\Psi)-\bm{u}^{(1)})e_{i\|j}(\bm{u}_\rho^{(1)}(\Psi)-\bm{u}^{(1)})\,\mathrm{d}x^{(1)}\\
			&\quad+\int_{\Omega^{(1)}}A^{ijk\ell,1}e_{k\|\ell}(\bm{u}_\rho^{(1)}(\Psi)-\bm{u}^{(1)})e_{i\|j}(\{(\mathbb{I}_{\Omega^{(1)}}+\bm{u}_\rho^{(1)}(\Psi)-\Psi\bm{e_3})\cdot\bm{e}_3\}^{-}\bm{e}_3)\,\mathrm{d}x^{(1)}\\
			&\quad+\dfrac{1}{\rho}\int_{\Omega^{(1)}}\bm{\beta}_1(\bm{u}_\rho^{(1)}(\Psi);\Psi)\cdot(\bm{u}_\rho^{(1)}(\Psi)-\bm{u}^{(1)})\,\mathrm{d}x\le\dfrac{1}{\rho}\int_{\Omega^{(1)}}|\bm{\beta}_1(\bm{u}_\rho^{(1)}(\Psi);\Psi)|^2\,\mathrm{d}x^{(1)}.
		\end{aligned}
	\end{equation}
	
	Combining an application of Korn's inequality, the fact that $\bm{\beta}_1(\bm{u}_\rho^{(1)}(\Psi);\Psi)\cdot(\bm{u}_\rho^{(1)}(\Psi)-\bm{u}^{(1)})\ge 0$ a.e. in $\Omega^{(1)}$, and the variational equations of Problem~\ref{problem1rho} transforms~\eqref{eq:1} into,
	\begin{equation*}
		\begin{aligned}
			&\dfrac{\|\bm{u}_\rho^{(1)}(\Psi)-\bm{u}^{(1)}\|_{\bm{V}(\Omega^{(1)})}^2}{C_e^{(1)} |c_0^{(1)}|^2}\le\int_{\Omega^{(1)}}\bm{f}^{(1)}\cdot\bm{\beta}_1(\bm{u}_\rho^{(1)}(\Psi);\Psi)\,\mathrm{d}x^{(1)}-\int_{\Omega^{(1)}}A^{ijk\ell,1}e_{k\|\ell}(\bm{u}^{(1)})e_{i\|j}(\bm{\beta}_1(\bm{u}_\rho^{(1)}(\Psi);\Psi))\,\mathrm{d}x^{(1)}\to 0,
		\end{aligned}
	\end{equation*}
	as $\rho\to0$. Additionally, since $\bm{u}^{(1)}(\Psi)$ is the unique solution of Problem~\ref{problem1}, the entire family $\{\bm{u}_\rho^{(1)}(\Psi)\}_{\rho>0}$ converges to $\bm{u}^{(1)}(\Psi)$ strongly in $\bm{V}(\Omega^{(1)})$, and the proof is complete.
\end{proof}

We are now ready to establish the existence and uniqueness of solutions for Problem~\ref{cprho}.
\begin{theorem}
	\label{ex-un-op-con-rho}
	Let $\bm{f}^{(\alpha)}\in\bm{L}^2(\Omega^{(\alpha)})$ be given. Then Problem~\ref{cprho} admits one and only one solution $\Psi_\rho\in\mathcal{V}(K)$.
\end{theorem}
\begin{proof}
	For $(\Psi,\bm{v}^{(1)},\bm{v}^{(2)})\in\mathcal{V}(K)\times\bm{V}(\Omega^{(1)}) \times\bm{V}(\Omega^{(2)})$ define the full functional
	\begin{equation}
		\label{eq:full-rho}
		\begin{aligned}
			&\widehat{\mathcal{J}}_\rho(\Psi,\bm{v}^{(1)},\bm{v}^{(2)}):=\dfrac{\beta}{2}\int_K|\Delta\Psi|^2\,\mathrm{d}x+\sum_{\alpha=1}^{2}E_\alpha(\bm{v}^{(\alpha)})\\
			&\quad+\dfrac{1}{2\rho}\int_{\Omega^{(1)}}|\bm{\beta}_1(\bm{v}^{(1)};\Psi)|^2\,\mathrm{d}x+\frac{1}{2\rho}\int_{\Omega^{(2)}}|\bm{\beta}_2(\bm{v}^{(2)};\Psi)|^2\,\mathrm{d}x,
		\end{aligned}
	\end{equation}
	where $E_\alpha(\bm{v}^{(\alpha)}):=\frac{1}{2}\int_{\Omega^{(\alpha)}}A^{ijk\ell,\alpha}e_{k\|\ell}(\bm{v}^{(\alpha)})e_{i\|j}(\bm{v}^{(\alpha)})\,\mathrm{d}x-\int_{\Omega^{(\alpha)}}\bm{f}^{(\alpha)}\cdot\bm{v}^{(\alpha)}\,\mathrm{d}x$.
	
	The proof of existence and uniqueness of solutions for Problem~\ref{cprho} is broken down into four steps, numbered $(i)$--$(iv)$.
	
	$(i)$ \emph{The functional $\widehat{\mathcal{J}}_\rho$ is strictly convex in the triple $(\Psi,\bm{v}^{(1)},\bm{v}^{(2)})$.} Thanks to~\eqref{eq:pos-def} and Korn's inequality (Theorem~\ref{KornCart}), each $E_\alpha$ is strictly convex on $\bm{V}(\Omega^{(\alpha)})$. Writing each $\Psi\in\mathcal{V}(K)$ as $\Psi=\psi\circ\pi$ with $\psi\in H^2_0(\kappa)$, an application of Theorem~8.7-1 of ~\cite{Cia25} in turn implies that the map $\Psi\mapsto\|\Delta\Psi\|_{L^2(K)}^2$ is strictly convex on $\mathcal{V}(K)$.
	
	Finally, each penalty term is convex, since the mappings $\bm{v}^{(\alpha)}\mapsto|\bm{\beta}_\alpha(\bm{v}^{(\alpha)};\Psi)|^2$ are obtained from the composition of the affine maps $\bm{v}^{(\alpha)}\mapsto(\mathbb{I}_{\Omega^{(\alpha)}}+\bm{v}^{(\alpha)}-\Psi\bm{e}_3)\cdot\bm{e}_3$ with the convex functions $t\mapsto|\{t\}^{\pm}|^2$. Hence $\widehat{\mathcal{J}}_\rho$ is strictly convex in the triple $(\Psi,\bm{v}^{(1)},\bm{v}^{(2)})$.
	
	$(ii)$ \emph{The functional $\widehat{\mathcal{J}}_\rho$ is coercive in the triple $(\Psi,\bm{v}^{(1)},\bm{v}^{(2)})$.} This descends immediately from Theorem~8.7-1 in~\cite{Cia25}, Korn's inequality (Theorem~\ref{KornCart}) and the non-negativity of the penalty terms.
	
	$(iii)$ \emph{The functional $\widehat{\mathcal{J}}_\rho$ admits a unique minimiser in $\mathcal{V}(K)\times\bm{V}(\Omega^{(1)})\times\bm{V}(\Omega^{(2)})$.} Each term of~\eqref{eq:full-rho} is convex and continuous and, hence, sequentially weakly lower semi-continuous. An application of Theorems~11.7-1 and~11.10-2 in~\cite{Cia25} allows to infer the existence and uniqueness of a minimiser for the functional $\widehat{\mathcal{J}}_\rho$ in $\mathcal{V}(K)\times\bm{V}(\Omega^{(1)})\times\bm{V}(\Omega^{(2)})$.
	
	$(iv)$ \emph{Existence and uniqueness of minimisers for the cost functional $J_\rho$}. For each fixed $\Psi\in\mathcal{V}(K)$, the minimisation of $\widehat{\mathcal{J}}_\rho(\Psi,\cdot,\cdot)$ decouples into the two \emph{independent} problems
	\begin{equation}
		\label{eq:part-min}
		\begin{aligned}
			&\inf_{\bm{v}^{(1)}\in\bm{V}(\Omega^{(1)})}\left\{E_1(\bm{v}^{(1)})+\frac{1}{2\rho}\|\bm{\beta}_1(\bm{v}^{(1)};\Psi)\|_{L^2(\Omega^{(1)})}^2\right\},\\
			&\inf_{\bm{v}^{(2)}\in\bm{V}(\Omega^{(2)})}\left\{E_2(\bm{v}^{(2)})+\frac{1}{2\rho}\|\bm{\beta}_2(\bm{v}^{(2)};\Psi)\|_{L^2(\Omega^{(2)})}^2\right\}.
		\end{aligned}
	\end{equation}
	
	Each of the functionals in~\eqref{eq:part-min} is strictly convex and coercive. Therefore, each of the functionals in~\eqref{eq:part-min} admits a unique minimiser $\bm{u}_\rho^{(\alpha)}\in\bm{V}(\Omega^{(\alpha)})$. Additionally, observe that the Euler--Lagrange equations associated with~\eqref{eq:part-min} are precisely the state equations in Problem~\ref{problem1rho} and Problem~\ref{problem2rho}, respectively. Therefore, the unique minimisers for the functionals in~\eqref{eq:part-min} are the penalised states $\bm{u}_\rho^{(1)}(\Psi)$ and $\bm{u}_\rho^{(2)}(\Psi)$, respectively. Besides, the unique minimiser of $\widehat{\mathcal{J}}_\rho$ satisfies:
	\begin{equation*}
		\bm{u}_\rho^{(1)}=\bm{u}_\rho^{(1)}(\Psi_\rho)\quad\textup{ and }\quad\bm{u}_\rho^{(2)}=\bm{u}_\rho^{(2)}(\Psi_\rho).
	\end{equation*}
	
	Consequently, for every $\Psi\in\mathcal{V}(K)$, it results
	\begin{equation*}
		\mathcal{J}_\rho(\Psi)=\widehat{\mathcal{J}}_\rho\left(\Psi,\bm{u}_\rho^{(1)}(\Psi),\bm{u}_\rho^{(2)}(\Psi)\right)=\inf_{(\bm{v}^{(1)},\bm{v}^{(2)})}\widehat{\mathcal{J}}_\rho(\Psi,\bm{v}^{(1)},\bm{v}^{(2)}),
	\end{equation*}
	and $\Psi_\rho$ is thus a minimiser of $\mathcal{J}_\rho$. Uniqueness of minimisers for $\mathcal{J}_\rho$ follows from uniqueness of minimisers for $\widehat{\mathcal{J}}_\rho$. Indeed, if $\Psi'_\rho\in\mathcal{V}(K)$ also minimised $\mathcal{J}_\rho$, then the triplet $\left(\Psi'_\rho,\bm{u}_\rho^{(1)}(\Psi'_\rho),\bm{u}_\rho^{(2)}(\Psi'_\rho)\right)$ would be another minimiser of $\widehat{\mathcal{J}}_\rho$, contradicting item $(iii)$. This completes the proof.
\end{proof}

Thanks to Lemma~\ref{lem:conv}, we are also able to show that minimisers of Problem~\ref{cprho} converge to minimisers of Problem~\ref{controlproblem} weakly in $H^2(K)$ as $\rho\to0$ without resorting to $\Gamma$-convergence~\cite{ABM14,DalMaso93}.
\begin{theorem}
	\label{th:optimal-1}
	Let $\bm{f}^{(\alpha)}\in\bm{L}^2(\Omega^{(\alpha)})$ be given. Then, the entire family $\{\Psi_\rho\}_{\rho>0}$ satisfies $\Psi_\rho\rightharpoonup\Psi$ in $H^2(K)$ as $\rho\to0$, where $\Psi$ is the unique solution of Problem~\ref{controlproblem}.
\end{theorem}
\begin{proof}
	The proof is divided into three steps, numbered $(i)$--$(iii)$.
	
	$(i)$ \emph{The function
		\begin{equation*}
			\Psi\mapsto\|\Delta\Psi\|_{L^2(K)},
		\end{equation*}
		is a norm equivalent to $\|\cdot\|_{H^2(K)}$ along the elements of the space $\{\Phi\in H^2(K) \textup{ such that } \Phi=\partial_{\bm{\nu}}\Phi=0\,\mathrm{d}\mathcal{H}^2\textup{-a.e. on }\partial\kappa\times(-a,a), \partial_3\Phi=0 \textup{ a.e. in }K\}$.}
	Since any $\Psi\in\{\Phi\in H^2(K) \textup{ such that } \Phi=\partial_{\bm{\nu}}\Phi=0\,\mathrm{d}\mathcal{H}^2\textup{-a.e. on }\partial\kappa\times(-a,a), \partial_3\Phi=0 \textup{ a.e. in }K\}$ can be identified with a unique function $\tilde{\psi}\in H^2_0(\kappa)$, it results that there exists a constant $C>0$ such that
	\begin{equation*}
		\begin{aligned}
		&|\Psi|_{2,K}^2:=\int_{K}\left(\sum_{i=1}^3|\partial_{ii}\Psi|^2+\sum_{\substack{i,j=1\\i\neq j}}^3|\partial_{ij}\Psi|^2\right)\,\mathrm{d}x\\
		&=2a\int_{\kappa}\left(|\partial_{11}\tilde{\psi}|^2+|\partial_{22}\tilde{\psi}|^2+2|\partial_{12}\tilde{\psi}|^2\right)\,\mathrm{d}x_1 \,\mathrm{d}x_2\\
		&=2a|\tilde{\psi}|_{2,\kappa}^2\ge 2aC^2\|\Delta\tilde{\psi}\|_{L^2(\kappa)}^2=C^2\|\Delta\Psi\|_{L^2(K)}^2,
		\end{aligned}
	\end{equation*}
	where the estimate holds thanks to Theorem~8.7-1 of ~\cite{Cia25}. In a similar fashion, as a result of an application of Theorem~8.7-1 of~\cite{Cia25}, we can show that there exists a constant $c_0>0$ such that
	\begin{equation*}
		|\Psi|_{2,K}^2=2a|\tilde{\psi}|_{2,\kappa}^2\le 2ac_0^2\|\Delta\tilde{\psi}\|_{L^2(\kappa)}^2=c_0^2\|\Delta\Psi\|_{L^2(K)}^2,
	\end{equation*}
	thus showing that $|\cdot|_{2,K}$ and the mapping $\Psi\mapsto\|\Delta\Psi\|_{L^2(K)}$ are equivalent along the elements of $\{\Phi\in H^2(K) \textup{ such that } \Phi=\partial_{\bm{\nu}}\Phi=0\, \mathrm{d}\mathcal{H}^2\textup{-a.e. on }\partial\kappa\times(-a,a), \partial_3\Phi=0 \textup{ a.e. in }K\}$.
	
	$(ii)$ \emph{Let $\{\Psi_\rho\}_{\rho>0}$ be a sequence of minimisers for Problem~\ref{cprho}. Let $\bm{u}^{(\alpha)}_\rho(\Psi_\rho)$ be the unique solution of either Problem~\ref{problem1rho} or Problem~\ref{problem2rho}. Then, the entire family $\{\bm{u}^{(\alpha)}_\rho(\Psi_\rho)\}_{\rho>0}$ weakly converges in $\bm{V}(\Omega^{(\alpha)})$ to the unique solution of either Problem~\ref{problem1} or Problem~\ref{problem2} as $\rho\to 0$}.
	Denote the three-dimensional linearly elastic energy by
	\begin{equation*}
		E_\alpha(\bm{v}^{(\alpha)}):=\frac{1}{2}\int_{\Omega^{(\alpha)}}A^{ijk\ell,\alpha}e_{k\|\ell}(\bm{v}^{(\alpha)})e_{i\|j}(\bm{v}^{(\alpha)})\,\mathrm{d}x^{(\alpha)}-\int_{\Omega^{(\alpha)}}\bm{f}^{(\alpha)}\cdot\bm{v}^{(\alpha)}\,\mathrm{d}x^{(\alpha)},
	\end{equation*}
	for all $\bm{v}^{(\alpha)}\in\bm{V}(\Omega^{(\alpha)})$.
	Since $\Psi_\rho$ is the unique minimiser of $\mathcal{J}_\rho$ over $\mathcal{V}(K)$, it results that:
	\begin{equation*}
		\mathcal{J}_\rho(\Psi_\rho)\le\mathcal{J}_\rho(0),\quad\textup{ for all }\rho>0.
	\end{equation*}
	
	An application of estimate~\eqref{eq:bdd:sol} in Lemma~\ref{lem:conv} with $\Psi=0$ gives that there exists a constant $D^{(\alpha)}(0,\bm{f}^{(\alpha)})$ independent of $\rho$ such that:
	\begin{equation*}
		\begin{aligned}
			\|\bm{\beta}_\alpha(\bm{u}_\rho^{(\alpha)}(0);0)\|_{\bm{L}^2(\Omega^{(\alpha)})}^2&\le\rho D^{(\alpha)}(0,\bm{f}^{(\alpha)}),\\
			\|\bm{u}_\rho^{(\alpha)}(0)\|_{\bm{V}(\Omega^{(\alpha)})}&\le\sqrt{2C_e^{(\alpha)}|c_0^{(\alpha)}|^2 D^{(\alpha)}(0,\bm{f}^{(\alpha)})}.
		\end{aligned}
	\end{equation*}
	
	Consequently the penalty terms satisfy
	\begin{equation*}
		\dfrac{1}{2\rho}\int_{\Omega^{(\alpha)}}|\bm{\beta}_\alpha(\bm{u}_\rho^{(\alpha)}(0);0)|^2\,\mathrm{d}x^{(\alpha)}\le\dfrac{1}{2} D^{(\alpha)}(0,\bm{f}^{(\alpha)}),
	\end{equation*}
	while the boundedness of $\bm{u}_\rho^{(\alpha)}(0)$ in $\bm{V}(\Omega^{(\alpha)})$ independently of $\rho$ and the continuity of $E_\alpha$ imply that the terms $E_\alpha(\bm{u}_\rho^{(\alpha)}(0))$ are bounded independently of $\rho$. Therefore, there exists a constant $C>0$ independent of $\rho$ such that:
	\begin{equation}
		\label{eq:bdd-J}
		\mathcal{J}_\rho(\Psi_\rho)\le\mathcal{J}_\rho(0)\le C.
	\end{equation}
	
	Combining~\eqref{eq:bdd-J} with an application of step~$(i)$ we obtain that the family $\{\Psi_\rho\}_{\rho>0}$ is bounded in $H^2(K)$ independently of $\rho$. Therefore, up to passing to a subsequence $\{\rho_n\}_{n=1}^\infty$ of positive numbers such that $\rho_n \to 0$ as $n\to\infty$, we obtain
	\begin{equation*}
		\Psi_{\rho_n}\rightharpoonup\Psi,\quad\textup{ in }H^2(K)\textup{ as }n\to\infty,
	\end{equation*}
	and $\Psi\in\mathcal{V}(K)$. Let $\bm{v}\in\bm{K}(\Omega^{(1)};\Psi)$ be given. For each $\rho>0$, define the vector field $\bm{v}_\rho\in\bm{V}(\Omega^{(1)})$ as follows:
	\begin{equation*}
		\bm{v}_\rho:=\bm{v}+\left\{(\mathbb{I}_{\Omega^{(1)}}+\bm{v})\cdot\bm{e}_3-\Psi_\rho\right\}^{-}\bm{e}_3.
	\end{equation*}
	
	Observe that
	\begin{equation}
		\label{constraint:psirho}
		\begin{aligned}
			&(\mathbb{I}_{\Omega^{(1)}}+\bm{v}_\rho)\cdot\bm{e}_3-\Psi_\rho
			=(\mathbb{I}_{\Omega^{(1)}}+\bm{v})\cdot\bm{e}_3-\Psi_\rho+\left\{(\mathbb{I}_{\Omega^{(1)}}+\bm{v})\cdot\bm{e}_3-\Psi_\rho\right\}^{-}\\
			&=\left\{(\mathbb{I}_{\Omega^{(1)}}+\bm{v})\cdot\bm{e}_3-\Psi_\rho\right\}^{+}\ge 0,
		\end{aligned}
	\end{equation}
	and that, as a consequence of Theorem~5.8.2 in~\cite{ABM14} and the convergence $\Psi_{\rho_n}\to\Psi$ in $H^1(K)$ as $n\to\infty$, which holds due to Rellich-Kondra\v{s}ov theorem (cf., e.g., Theorem~8.4-3 in~\cite{Cia25}), we obtain $\bm{v}_{\rho_n}\to\bm{v}$ in $\bm{V}(\Omega^{(1)})$ as $n\to\infty$. In correspondence of a \emph{given} minimiser $\Psi_\rho$ for the reduced cost functional $\mathcal{J}_\rho$, we denote by $\bm{u}^{(1)}_\rho(\Psi_\rho)$ the unique solution of Problem~\ref{problem1rho}.
	
	Define the function $\bm{r}_\rho^{(1)}:=\{(\mathbb{I}_{\Omega^{(1)}}-\Psi_\rho\bm{e}_3)\cdot\bm{e}_3\}^{-}\bm{e}_3$ and observe that~\eqref{bdry:feas} implies that $\bm{r}_\rho^{(1)}\in\bm{V}(\Omega^{(1)})$. Additionally, a straightforward calculation shows that $\bm{r}_\rho^{(1)}\in\bm{K}(\Omega^{(1)};\Psi_\rho)$. Therefore, we obtain that:
	\begin{equation}
		\label{energy-bound}
		E_1(\bm{u}_\rho^{(1)}(\Psi_\rho))+\dfrac{1}{2\rho}\|\bm{\beta}_1(\bm{u}_\rho^{(1)}(\Psi_\rho);\Psi_\rho)\|_{\bm{L}^2(\Omega^{(1)})}^2\le E_1(\bm{r}_\rho^{(1)}),\quad\textup{ for all }\rho>0.
	\end{equation}
	
	An application of Korn's inequality (Theorem~\ref{KornCart}), Young's inequality~\cite{Young1912} and the boundedness of the family $\{\Psi_\rho\}_{\rho>0}$ to~\eqref{energy-bound} gives that there exists a constant $C>0$ independent of $\rho$ such that $\|\bm{u}_\rho^{(1)}(\Psi_\rho)\|_{\bm{V}(\Omega^{(1)})}\le C$, for all $\rho>0$. Therefore, it results $\bm{u}_{\rho_n}^{(1)}(\Psi_{\rho_n})\rightharpoonup\bm{u}_0^{(1)}$ in $\bm{V}(\Omega^{(1)})$ as $n\to\infty$. Up to passing to a further subsequence that we do not relabel, an application of the Rellich-Kondra\v{s}ov theorem (Theorem~8.4-3 in~\cite{Cia25}) gives $\Psi_{\rho_n}\to\Psi$ in $H^1(K)$ as $n\to\infty$, and $\bm{u}_{\rho_n}^{(1)}(\Psi_{\rho_n})\to\bm{u}_0^{(1)}$ in $\bm{L}^2(\Omega^{(1)})$ as $n\to\infty$. Since by~\eqref{energy-bound} it results $\bm{\beta}_1(\bm{u}_{\rho_n}^{(1)}(\Psi_{\rho_n});\Psi_{\rho_n})\to\bm{0}$ as $n\to\infty$, we obtain that $\bm{u}_0^{(1)}\in\bm{K}(\Omega^{(1)};\Psi)$.
	Testing the variational equations of Problem~\ref{problem1rho} at $(\bm{v}_\rho-\bm{u}^{(1)}_\rho(\Psi_\rho))$ gives:
	\begin{equation}
		\label{eq:test-1}
		\begin{aligned}
			&\int_{\Omega^{(1)}} A^{ijk\ell,1} e_{k\|\ell}(\bm{u}_\rho^{(1)}(\Psi_\rho)) e_{i\|j}(\bm{v}_\rho-\bm{u}^{(1)}_\rho(\Psi_\rho)) \,\mathrm{d}x^{(1)}+\dfrac{1}{\rho}\int_{\Omega^{(1)}}\bm{\beta}_1(\bm{u}_\rho^{(1)}(\Psi_\rho);\Psi_\rho)\cdot(\bm{v}_\rho-\bm{u}^{(1)}_\rho(\Psi_\rho))\,\mathrm{d}x^{(1)}\\
			&=\int_{\Omega^{(1)}} \bm{f}^{(1)} \cdot (\bm{v}_\rho-\bm{u}^{(1)}_\rho(\Psi_\rho)) \,\mathrm{d}x^{(1)}.
		\end{aligned}
	\end{equation}
	
	An application of~\eqref{constraint:psirho} in turn gives that $\bm{\beta}_1(\bm{v}_\rho;\Psi_\rho)=\bm{0}$ a.e. in $\Omega^{(1)}$. Thanks to the monotonicity of $\bm{\beta}_1(\cdot;\Psi_\rho)$, the validity of the convergence $\bm{v}_{\rho_n}\to\bm{v}$ in $\bm{V}(\Omega^{(1)})$ as $n\to\infty$, the boundedness of the sequence $\{\bm{u}^{(1)}_{\rho_n}(\Psi_{\rho_n})\}_{n=1}^\infty$ in $\bm{V}(\Omega^{(1)})$ and the sequential weak lower semicontinuity of the mapping
	\begin{equation*}
		\bm{v}\in\bm{V}(\Omega^{(1)})\mapsto\int_{\Omega^{(1)}}A^{ijk\ell,1} e_{k\|\ell}(\bm{v}) e_{i\|j}(\bm{v}) \,\mathrm{d}x^{(1)},
	\end{equation*}
	we obtain that a limit passage as $n\to\infty$ in~\eqref{eq:test-1} along the sequence $\{\rho_n\}_{n=1}^\infty$ gives:
	\begin{equation*}
		\int_{\Omega^{(1)}} A^{ijk\ell,1} e_{k\|\ell}(\bm{u}_0^{(1)}) e_{i\|j}(\bm{v}-\bm{u}_0^{(1)}) \,\mathrm{d}x^{(1)}\ge\int_{\Omega^{(1)}} \bm{f}^{(1)} \cdot (\bm{v}-\bm{u}_0^{(1)}) \,\mathrm{d}x^{(1)}.
	\end{equation*}
	
	Since $\bm{v}$ was chosen arbitrarily in $\bm{K}(\Omega^{(1)};\Psi)$, and since we know that Problem~\ref{problem1} admits a unique solution, then it must be $\bm{u}_0^{(1)}=\bm{u}^{(1)}(\Psi)$. The uniqueness of the cluster point $\bm{u}^{(1)}(\Psi)$ (Theorem~\ref{th:ex1}) for the sequence $\{\bm{u}_{\rho_n}^{(1)}(\Psi_{\rho_n})\}_{n=1}^\infty$  implies that the entire family $\{\bm{u}_\rho^{(1)}(\Psi_\rho)\}_{\rho>0}$ weakly converges to $\bm{u}^{(1)}(\Psi)$ in $\bm{V}(\Omega^{(1)})$ as $\rho\to 0$.
	In a similar fashion, we can show that $\bm{u}^{(2)}_\rho(\Psi_\rho)\rightharpoonup\bm{u}^{(2)}(\Psi)$ in $\bm{V}(\Omega^{(2)})$ as $\rho\to0$, with $\bm{u}^{(2)}(\Psi)$ unique solution of Problem~\ref{problem2} by testing the variational equations of Problem~\ref{problem2rho} at $(\bm{w}_\rho-\bm{u}^{(2)}_\rho(\Psi_\rho))$, where
	\begin{equation*}
		\bm{w}_\rho:=\bm{w}-\left\{(\mathbb{I}_{\Omega^{(2)}}+\bm{w})\cdot\bm{e}_3-\Psi_\rho\right\}^{+}\bm{e}_3,
	\end{equation*}
	and the element $\bm{w}$ is chosen arbitrarily in $\bm{K}(\Omega^{(2)};\Psi)$.
	
	$(iii)$ \emph{Optimality of the weak limit and conclusion of the proof.} We work along the subsequence $\{\rho_n\}_{n=1}^\infty$ extracted in part $(ii)$, for which:
	\begin{equation*}
		\begin{aligned}
			\Psi_{\rho_n}\rightharpoonup\Psi&,\textup{ in }H^2(K)\textup{ as }n\to\infty,\\
			\bm{u}_{\rho_n}^{(\alpha)}(\Psi_{\rho_n})\rightharpoonup\bm{u}^{(\alpha)}(\Psi)&,\textup{ in }\bm{V}(\Omega^{(\alpha)})\textup{ as }n\to\infty\textup{ for each }\alpha\in\{1,2\}.
		\end{aligned}
	\end{equation*}
	
	Fix $\Phi\in\mathcal{V}(K)$. An application of Lemma~\ref{lem:conv} and the continuity of $E_\alpha$ give, respectively:
	\begin{equation}
		\label{conv1}
		\begin{aligned}
			\bm{u}_\rho^{(\alpha)}(\Phi)\to\bm{u}^{(\alpha)}(\Phi)&,\quad\textup{ in }\bm{V}(\Omega^{(\alpha)})\textup{ as }\rho\to0,\\
			E_\alpha(\bm{u}_\rho^{(\alpha)}(\Phi))\to E_\alpha(\bm{u}^{(\alpha)}(\Phi))&,\quad\textup{ as }\rho\to 0.
		\end{aligned}
	\end{equation}

	The fact that $\bm{u}^{(\alpha)}(\Phi)\in\bm{K}(\Omega^{(\alpha)};\Phi)$ implies that  $\bm{\beta}_\alpha(\bm{u}^{(\alpha)}(\Phi);\Phi)=\bm{0}$. The fact that $\bm{u}_\rho^{(\alpha)}(\Phi)$ is the unique solution of Problem~\ref{problem1rho} when $\alpha=1$ and of Problem~\ref{problem2rho} when $\alpha=2$, gives:
	\begin{equation}
		\label{eq:seb}
		E_\alpha(\bm{u}_\rho^{(\alpha)}(\Phi))+\dfrac{1}{2\rho}\left\|\bm{\beta}_\alpha(\bm{u}_\rho^{(\alpha)}(\Phi);\Phi)\right\|^2_{\bm{L}^2(\Omega^{(\alpha)})}
		\le E_\alpha(\bm{u}^{(\alpha)}(\Phi)).
	\end{equation}
	
	An application of the second convergence in~\eqref{conv1} gives
	\begin{equation*}
		0\le\dfrac{1}{2\rho}\left\|\bm{\beta}_\alpha(\bm{u}_\rho^{(\alpha)}(\Phi);\Phi)\right\|^2_{\bm{L}^2(\Omega^{(\alpha)})}\le E_\alpha(\bm{u}^{(\alpha)}(\Phi))-E_\alpha(\bm{u}_\rho^{(\alpha)}(\Phi))\to0,\quad\textup{ as }\rho\to 0,
	\end{equation*}
	and, consequently, it results:
	\begin{equation}
		\label{eq:conv-J-Phi}
		\lim_{\rho\to0}\mathcal{J}_\rho(\Phi)=\dfrac{\beta}{2}\|\Delta\Phi\|^2_{L^2(K)}+\sum_{\alpha=1}^2E_\alpha(\bm{u}^{(\alpha)}(\Phi))=\mathcal{J}(\Phi).
	\end{equation}
	
	The map $\Psi\mapsto\|\Delta\Psi\|^2_{L^2(K)}$ is convex and continuous and, thus, weakly lower semi-continuous. Besides, each $E_\alpha$ is weakly lower semi-continuous on $\bm{V}(\Omega^{(\alpha)})$ as the sum of a coercive continuous quadratic form and a continuous linear form. Exploiting the weak convergences established in item $(ii)$ allows us to infer that:
	\begin{equation}
		\label{liminf-ineq}
		\begin{aligned}
			&\mathcal{J}(\Psi)=\dfrac{\beta}{2}\|\Delta\Psi\|^2_{L^2(K)}+\sum_{\alpha=1}^2E_\alpha(\bm{u}^{(\alpha)}(\Psi))\\
			&\le\liminf_{n\to\infty}\left(\frac{\beta}{2}\|\Delta\Psi_{\rho_n}\|^2_{L^2(K)}+\sum_{\alpha=1}^2 E_\alpha(\bm{u}_{\rho_n}^{(\alpha)}(\Psi_{\rho_n}))\right)\\
			&\le\liminf_{n\to\infty}\mathcal{J}_{\rho_n}(\Psi_{\rho_n}).
		\end{aligned}
	\end{equation}
	
	Since $\Psi_{\rho_n}$ minimises $\mathcal{J}_{\rho_n}$ over $\mathcal{V}(K)$, we obtain that for every $\Phi\in\mathcal{V}(K)$ it results $\mathcal{J}_{\rho_n}(\Psi_{\rho_n})\le\mathcal{J}_{\rho_n}(\Phi)$. Combining~\eqref{liminf-ineq} with~\eqref{eq:conv-J-Phi} gives
	\begin{equation*}
		\mathcal{J}(\Psi)\le\liminf_{n\to\infty}\mathcal{J}_{\rho_n}(\Psi_{\rho_n})\le\liminf_{n\to\infty}\mathcal{J}_{\rho_n}(\Phi)=\mathcal{J}(\Phi),\quad\textup{ for all }\Phi\in\mathcal{V}(K),
	\end{equation*}
	thus showing that $\Psi$ minimises $\mathcal{J}$ over $\mathcal{V}(K)$. The fact that Problem~\ref{controlproblem} admits a unique solution (Theorem~\ref{ex-un-op-con}) implies that the weak limit $\Psi$ is actually the unique solution of Problem~\ref{controlproblem}.
	
	Observing that the argument above applies to \emph{every subsequence} of the family $\{\Psi_\rho\}_{\rho>0}$, it results that the entire family $\{\Psi_\rho\}_{\rho>0}$ satisfies
	\begin{equation*}
		\Psi_\rho\rightharpoonup\Psi,\quad\textup{ in }H^2(K)\textup{ as }\rho\to 0,
	\end{equation*}
	and, by part $(ii)$, the families $\{\bm{u}_\rho^{(\alpha)}(\Psi_\rho)\}_{\rho>0}$ satisfy
	\begin{equation*}
		\bm{u}_\rho^{(\alpha)}(\Psi_\rho)\rightharpoonup\bm{u}^{(\alpha)}(\Psi),\quad\textup{ in }\bm{V}(\Omega^{(\alpha)})\textup{ as }\rho\to 0\textup{ for each }\alpha\in\{1,2\},
	\end{equation*}
	and the proof is complete.
\end{proof}

\begin{remark}
	In the proof of Theorem~\ref{th:optimal-1}, the $H^2(K)$ regularity for the admissible controls is key to assert the convergence $\bm{v}_{\rho_n}\to\bm{v}$ in $\bm{V}(\Omega^{(\alpha)})$ as $n\to\infty$.
	\bqed
\end{remark}

\section{Existence and Regularity of Optimal Interfaces. Part~II}\label{Sec:4}
The purpose of this section is to provide a first order necessary condition for approximating solutions of Problem~\ref{controlproblem}. To this aim, we will further penalise Problem~\ref{cprho} by smoothening the penalty term. Indeed, the non-linear penalty term considered in the formulation of Problem~\ref{cprho} is, in general, not differentiable. For smoothening the non-linear penalty term in Problem~\ref{cprho}, we employ the Moreau--Yosida approximation (cf., e.g., Lemma~5.8.1 in~\cite{ABM14}). Let $\gls{c}>0$ be another penalty parameter meant to approach zero and consider the function $b_{1,c}:\mathbb{R}\to\mathbb{R}$ defined by
\begin{equation*}
	b_{1,c}(r):=
	\begin{cases}
		r+\dfrac{c}{2}&,\textup{ if }r\le -c,\\
		\\
		-\dfrac{r^2}{2c}&,\textup{ if } -c<r\le 0,\\
		\\
		0&,\textup{ if } r>0,
	\end{cases}
\end{equation*}
and observe that $b_{1,c}\in\mathcal{C}^1(\mathbb{R})$, and that:
\begin{equation*}
	b_{1,c}'(r):=
	\begin{cases}
		1&,\textup{ if }r\le -c,\\
		\\
		-\dfrac{r}{c}&,\textup{ if } -c<r\le 0,\\
		\\
		0&,\textup{ if } r>0.
	\end{cases}
\end{equation*}

Since $b_{1,c}'\ge 0$ and $b_{1,c}'$ is bounded in $\mathbb{R}$ independently of $c$, we obtain that $b_{1,c}$ is monotone and Lipschitz continuous, with Lipschitz constant equal to one. Given $c>0$ and given $\Psi\in\mathcal{V}(K)$, define $\gls{beta1c}$ by:
\begin{equation*}
	\bm{\beta}_{1,c}(\bm{v};\Psi):=b_{1,c}((\mathbb{I}_{\Omega^{(1)}}+\bm{v})\cdot\bm{e}_3-\Psi)\bm{e}_3.
\end{equation*}

Since $b_{1,c}$ is Lipschitz continuous in $\bm{L}^2(\Omega^{(1)})$ with Lipschitz constant equal to one and monotone, then $\bm{\beta}_{1,c}(\cdot;\Psi):\bm{V}(\Omega^{(1)})\to\bm{V}(\Omega^{(1)})$ is Lipschitz continuous with Lipschitz constant equal to one and monotone.
We are now ready to state the regularised version of Problem~\ref{problem1rho}.
\begin{customprob}{$\mathcal{P}_{\rho,c}(\Omega^{(1)};\Psi)$}\label{problem1rhoc}
	Find $\bm{u}_{\rho,c}^{(1)}(\Psi)=(u_{\rho,c,i}^{(1)}(\Psi))\in\bm{V}(\Omega^{(1)})$ that satisfies the variational equations
	\begin{equation*}
		\int_{\Omega^{(1)}} A^{ijk\ell,1} e_{k\|\ell}(\bm{u}_{\rho,c}^{(1)}(\Psi)) e_{i\|j}(\bm{v}) \,\mathrm{d}x^{(1)}+\dfrac{1}{\rho}\int_{\Omega^{(1)}}\bm{\beta}_{1,c}(\bm{u}_{\rho,c}^{(1)}(\Psi);\Psi)\cdot\bm{v}\,\mathrm{d}x^{(1)}=\int_{\Omega^{(1)}}\bm{f}^{(1)}\cdot\bm{v}\,\mathrm{d}x^{(1)},
	\end{equation*}
	for all $\bm{v}=(v_i)\in\bm{V}(\Omega^{(1)})$.
	\bqed
\end{customprob}

In the next theorem, for given $\rho>0$, $c>0$ and $\Psi\in \mathcal{V}(K)$, we establish the existence and uniqueness of solutions for Problem~\ref{problem1rhoc}.
\begin{theorem}
	\label{th:exun:reg}
	Let $\rho>0$, $c>0$ and $\Psi\in \mathcal{V}(K)$ be given. Then, Problem~\ref{problem1rhoc} admits a unique solution. Besides, we have
	\begin{equation*}
		\bm{u}^{(1)}_{\rho,c}(\Psi)\to\bm{u}^{(1)}_\rho(\Psi),\quad\textup{ in }\bm{V}(\Omega^{(1)}) \textup{ as }c\to 0,
	\end{equation*}
	where $\bm{u}^{(1)}_\rho(\Psi)$ is the unique solution of Problem~\ref{problem1rho}.
\end{theorem}
\begin{proof}
	That Problem~\ref{problem1rhoc} admits a unique solution follows from a direct application of the Minty--Browder theorem (cf., e.g., Theorem~12.5.3 in~\cite{Cia25}), being the operator $\bm{\beta}_{1,c}(\cdot;\Psi)$ monotone.
	
	The proof of the convergence $\bm{u}^{(1)}_{\rho,c}(\Psi)\to\bm{u}^{(1)}_\rho(\Psi)$ in $\bm{V}(\Omega^{(1)})$ as $c\to 0$, where $\bm{u}^{(1)}_\rho(\Psi)$ is the unique solution of Problem~\ref{problem1rho} is divided into four steps, numbered $(i)$--$(iv)$.
	
	$(i)$ \emph{It results that $\bm{\beta}_{1,c}(\bm{v};\Psi)\to\bm{\beta}_1(\bm{v};\Psi)$ in $\bm{L}^2(\Omega^{(1)})$ as $c\to 0$, for every $\bm{v}\in\bm{L}^2(\Omega^{(1)})$.}
	Fix $\bm{v}\in\bm{L}^2(\Omega^{(1)})$ and observe that the estimate $|b_{1,c}(s)-\min\{s,0\}|\le c$ for all $s\in\mathbb{R}$ gives
	\begin{equation*}
		\|\bm{\beta}_{1,c}(\bm{v};\Psi)-\bm{\beta}_1(\bm{v};\Psi)\|_{\bm{L}^2(\Omega^{(1)})}\le c|\Omega^{(1)}|^{1/2}\to 0,\quad\textup{ as }c\to 0,
	\end{equation*}
	which proves the claim.
	
	$(ii)$ \emph{It results that $\bm{\beta}_{1,c}(\bm{v}_c;\Psi)\to\bm{\beta}_1(\bm{v};\Psi)$ in $\bm{L}^2(\Omega^{(1)})$ as $c\to 0$, for every family $\{\bm{v}_c\}_{c>0}\subset\bm{L}^2(\Omega^{(1)})$ such that $\bm{v}_c\to\bm{v}$ in $\bm{L}^2(\Omega^{(1)})$ as $c\to 0$.}
	Given a family $\{\bm{v}_c\}_{c>0}\subset\bm{L}^2(\Omega^{(1)})$ such that $\bm{v}_c\to\bm{v}$ in $\bm{L}^2(\Omega^{(1)})$ as $c\to 0$, evaluate
	\begin{equation*}
		\begin{aligned}
		&\|\bm{\beta}_{1,c}(\bm{v}_c;\Psi)-\bm{\beta}_1(\bm{v};\Psi)\|_{\bm{L}^2(\Omega^{(1)})}
		\le\|\bm{\beta}_{1,c}(\bm{v}_c;\Psi)-\bm{\beta}_{1,c}(\bm{v};\Psi)\|_{\bm{L}^2(\Omega^{(1)})}+\|\bm{\beta}_{1,c}(\bm{v};\Psi)-\bm{\beta}_1(\bm{v};\Psi)\|_{\bm{L}^2(\Omega^{(1)})}\\
		&\le\|\bm{v}_c-\bm{v}\|_{\bm{L}^2(\Omega^{(1)})}+\|\bm{\beta}_{1,c}(\bm{v};\Psi)-\bm{\beta}_1(\bm{v};\Psi)\|_{\bm{L}^2(\Omega^{(1)})},
		\end{aligned}
	\end{equation*}
	and the result follows by combining the assumed convergence $\bm{v}_c\to\bm{v}$ in $\bm{L}^2(\Omega^{(1)})$ as $c\to 0$ with item $(i)$.
	
	$(iii)$ \emph{It results, up to passing to a subsequence, that $\bm{u}^{(1)}_{\rho,c}(\Psi)\rightharpoonup\bm{u}^{(1)}_\rho(\Psi)$ in $\bm{V}(\Omega^{(1)})$ as $c\to 0$, where $\bm{u}^{(1)}_\rho(\Psi)$ is the unique solution of Problem~\ref{problem1rho}.} Testing the variational equations of Problem~\ref{problem1rhoc} at $\bm{u}^{(1)}_{\rho,c}(\Psi)$, we obtain that an application of Korn's inequality (Theorem~\ref{KornCart}) and an exploitation of the uniform positive-definiteness of the fourth order three-dimensional elasticity tensor $\{A^{ijk\ell,1}\}$, of the monotonicity of $\bm{\beta}_{1,c}(\cdot;\Psi)$ and the Cauchy--Schwarz inequality give
	\begin{equation}
		\label{eq:bdd:sol-c}
		\begin{aligned}
			&\dfrac{\|\bm{u}^{(1)}_{\rho,c}(\Psi)\|_{\bm{V}(\Omega^{(1)})}^2}{C_e^{(1)} |c_0^{(1)}|^2}\le \int_{\Omega^{(1)}} A^{ijk\ell,1} e_{k\|\ell}(\bm{u}^{(1)}_{\rho,c}(\Psi)) e_{i\|j}(\bm{u}^{(1)}_{\rho,c}(\Psi)) \,\mathrm{d}x^{(1)}\\
			&\quad+\dfrac{1}{\rho}\int_{\Omega^{(1)}}(\bm{\beta}_{1,c}(\bm{u}^{(1)}_{\rho,c}(\Psi);\Psi)-\bm{\beta}_{1,c}(\bm{0};\Psi))\cdot\bm{u}^{(1)}_{\rho,c}(\Psi)\,\mathrm{d}x^{(1)}\\
			&= \int_{\Omega^{(1)}} \bm{f}^{(1)} \cdot \bm{u}^{(1)}_{\rho,c}(\Psi) \,\mathrm{d}x^{(1)}-\dfrac{1}{\rho}\int_{\Omega^{(1)}}\bm{\beta}_{1,c}(\bm{0};\Psi)\cdot\bm{u}^{(1)}_{\rho,c}(\Psi)\,\mathrm{d}x^{(1)}\\
			&\le\left(\|\bm{f}^{(1)}\|_{\bm{L}^2(\Omega^{(1)})}+\dfrac{\|\bm{\beta}_{1,c}(\bm{0};\Psi)\|_{\bm{L}^2(\Omega^{(1)})}}{\rho}\right)\|\bm{u}^{(1)}_{\rho,c}(\Psi)\|_{\bm{V}(\Omega^{(1)})}\\
			&\le\left(\|\bm{f}^{(1)}\|_{\bm{L}^2(\Omega^{(1)})}+\dfrac{c|\Omega^{(1)}|^{1/2}+\|\bm{\beta}_1(\bm{0};\Psi)\|_{\bm{L}^2(\Omega^{(1)})}}{\rho}\right)\|\bm{u}^{(1)}_{\rho,c}(\Psi)\|_{\bm{V}(\Omega^{(1)})},
		\end{aligned}
	\end{equation}
	so that it results $\|\bm{u}^{(1)}_{\rho,c}(\Psi)\|_{\bm{V}(\Omega^{(1)})}\le C_e^{(1)} |c_0^{(1)}|^2(\|\bm{f}^{(1)}\|_{\bm{L}^2(\Omega^{(1)})}+\rho^{-1}c|\Omega^{(1)}|^{1/2}+\rho^{-1}\|\bm{\beta}_1(\bm{0};\Psi)\|_{\bm{L}^2(\Omega^{(1)})})$ which is bounded independently of $c$.
	Therefore, given any subsequence $\{\bm{u}^{(1)}_{\rho,c_n}(\Psi)\}_{n=1}^\infty$, where $\{c_n\}_{n=1}^\infty$ is a suitable sequence of positive numbers such that $c_n\to 0$ as $n\to\infty$, we can extract a further subsequence that weakly converges to a limit denoted by $\bm{u}^{(1)}_{\rho,0}\in\bm{V}(\Omega^{(1)})$ as $n\to\infty$. Note that an application of the Rellich-Kondra\v{s}ov theorem (cf., e.g., Theorem~8.4-3 in~\cite{Cia25}) gives, up to passing to an even further subsequence that we do not relabel, that $\bm{u}^{(1)}_{\rho,c_n}(\Psi)\to\bm{u}^{(1)}_{\rho,0}$ in $\bm{L}^2(\Omega^{(1)})$ as $n\to\infty$. To verify that $\bm{u}^{(1)}_{\rho,0}$ is the unique solution of Problem~\ref{problem1rho} note that, for any $\bm{v}\in\bm{V}(\Omega^{(1)})$, an application of the convergence $\bm{u}^{(1)}_{\rho,c_n}(\Psi)\rightharpoonup\bm{u}^{(1)}_{\rho,0}$ in $\bm{V}(\Omega^{(1)})$ as $n\to\infty$ and item $(ii)$ give
	\begin{equation*}
		\begin{aligned}
			\int_{\Omega^{(1)}} A^{ijk\ell,1} e_{k\|\ell}(\bm{u}_{\rho,c_n}^{(1)}(\Psi)) e_{i\|j}(\bm{v}) \,\mathrm{d}x^{(1)}\to&&\int_{\Omega^{(1)}} A^{ijk\ell,1} e_{k\|\ell}(\bm{u}_{\rho,0}^{(1)}) e_{i\|j}(\bm{v}) \,\mathrm{d}x^{(1)}&,\quad\textup{ as }n\to\infty\\
			\dfrac{1}{\rho}\int_{\Omega^{(1)}}\bm{\beta}_{1,c_n}(\bm{u}_{\rho,c_n}^{(1)}(\Psi);\Psi)\cdot\bm{v}\,\mathrm{d}x^{(1)}\to&&\dfrac{1}{\rho}\int_{\Omega^{(1)}}\bm{\beta}_{1}(\bm{u}_{\rho,0}^{(1)};\Psi)\cdot\bm{v}\,\mathrm{d}x^{(1)}&,\quad\textup{ as }n\to\infty,
		\end{aligned}
	\end{equation*}
	and the uniqueness of solutions for Problem~\ref{problem1rho} allows us to infer that $\bm{u}_{\rho,0}^{(1)}=\bm{u}_{\rho}^{(1)}(\Psi)$.
	
	$(iv)$ \emph{The weak convergence established in step $(iii)$ is actually strong. Besides, the entire family $\{\bm{u}^{(1)}_{\rho,c}(\Psi)\}_{c>0}$ converges to  $\bm{u}^{(1)}_\rho(\Psi)$ in $\bm{V}(\Omega^{(1)})$ as $c\to 0$.} Subtracting the variational equations of Problem~\ref{problem1rho} tested at $(\bm{u}^{(1)}_{\rho,c}(\Psi)-\bm{u}^{(1)}_\rho(\Psi))$ from the variational equations of Problem~\ref{problem1rhoc} tested at $(\bm{u}^{(1)}_{\rho,c}(\Psi)-\bm{u}^{(1)}_\rho(\Psi))$, we obtain:
	\begin{equation*}
		\begin{aligned}
		&\int_{\Omega^{(1)}}A^{ijk\ell,1}e_{k\|\ell}(\bm{u}^{(1)}_{\rho,c}(\Psi)-\bm{u}^{(1)}_\rho(\Psi))e_{i\|j}(\bm{u}^{(1)}_{\rho,c}(\Psi)-\bm{u}^{(1)}_\rho(\Psi))\,\mathrm{d}x^{(1)}\\
		&\quad+\dfrac{1}{\rho}\int_{\Omega^{(1)}}\left(\bm{\beta}_{1,c}(\bm{u}^{(1)}_{\rho,c}(\Psi);\Psi)-\bm{\beta}_1(\bm{u}^{(1)}_{\rho}(\Psi);\Psi)\right)\cdot (\bm{u}^{(1)}_{\rho,c}(\Psi)-\bm{u}^{(1)}_\rho(\Psi))\,\mathrm{d}x^{(1)}=0.
		\end{aligned}
	\end{equation*}
	
	Observe that from the variational equations of Problem~\ref{problem1rhoc} we infer that
	\begin{equation}
		\label{eq:bdd:penalty-c}
		\begin{aligned}
			&\dfrac{1}{\rho}\left\|\bm{\beta}_{1,c}(\bm{u}_{\rho,c}^{(1)}(\Psi);\Psi)\right\|_{\bm{V}^\ast(\Omega^{(1)})}\le\|\bm{f}^{(1)}\|_{\bm{L}^2(\Omega^{(1)})}+\sup_{\substack{\bm{v}\in\bm{V}(\Omega^{(1)})\\\|\bm{v}\|_{\bm{V}(\Omega^{(1)})}=1}}\left|\int_{\Omega^{(1)}}A^{ijk\ell,1}e_{k\|\ell}(\bm{u}_{\rho,c}^{(1)}(\Psi))e_{i\|j}(\bm{v})\,\mathrm{d}x^{(1)}\right|\\
			&\le\|\bm{f}^{(1)}\|_{\bm{L}^2(\Omega^{(1)})}+\max_{i,j,k,\ell}\|A^{ijk\ell,1}\|_{\mathcal{C}^0(\overline{\Omega^{(1)}})}\|\bm{u}_{\rho,c}^{(1)}(\Psi)\|_{\bm{V}(\Omega^{(1)})},
		\end{aligned}
	\end{equation}
	and we observe that an application of~\eqref{eq:bdd:sol-c} gives that the right-hand side of~\eqref{eq:bdd:penalty-c} is bounded independently of $c$.
	Thanks to the the uniform positive-definiteness of the fourth order three-dimensional elasticity tensor $\{A^{ijk\ell,1}\}$, Korn's inequality (Theorem~\ref{KornCart}), steps $(ii)$ and $(iii)$, and the convergence $\bm{u}^{(1)}_{\rho,c_n}(\Psi)\to\bm{u}^{(1)}_\rho(\Psi)$ in $\bm{L}^2(\Omega^{(1)})$ as $n\to\infty$ along a sequence $\{c_n\}_{n=1}^\infty$ of positive numbers such that $c_n\to 0$ as $n\to\infty$, we obtain that
	\begin{equation*}
		\|\bm{u}^{(1)}_{\rho,c_n}(\Psi)-\bm{u}^{(1)}_\rho(\Psi)\|_{\bm{V}(\Omega^{(1)})}^2\le-\dfrac{C_e^{(1)}|c_0^{(1)}|^2}{\rho}\int_{\Omega^{(1)}}(\bm{\beta}_{1,c_n}(\bm{u}^{(1)}_{\rho,c_n}(\Psi);\Psi)-\bm{\beta}_1(\bm{u}^{(1)}_{\rho}(\Psi);\Psi))\cdot (\bm{u}^{(1)}_{\rho,c_n}(\Psi)-\bm{u}^{(1)}_\rho(\Psi))\,\mathrm{d}x^{(1)} \to 0,
	\end{equation*}
	as $n\to\infty$, thus showing that $\bm{u}^{(1)}_{\rho,c_n}(\Psi)\to\bm{u}^{(1)}_\rho(\Psi)$ in $\bm{V}(\Omega^{(1)})$ as $n\to\infty$.
	Since, for a given $\rho>0$ and for a given $\Psi\in\mathcal{V}(K)$, the solution of Problem~\ref{problem1rho} is uniquely determined, we obtain that the \emph{entire family} $\{\bm{u}^{(1)}_{\rho,c}(\Psi)\}_{c>0}$ converges.
	This completes the proof.
\end{proof}

In a similar fashion, consider the function $b_{2,c}:\mathbb{R}\to\mathbb{R}$ defined by
\begin{equation*}
	b_{2,c}(r):=
	\begin{cases}
		r-\dfrac{c}{2}&,\textup{ if }r\ge c,\\
		\\
		\dfrac{r^2}{2c}&,\textup{ if } 0<r\le c,\\
		\\
		0&,\textup{ if } r<0,
	\end{cases}
\end{equation*}
and observe that $b_{2,c}\in\mathcal{C}^1(\mathbb{R})$, and that:
\begin{equation*}
	b_{2,c}'(r):=
	\begin{cases}
		1&,\textup{ if }r\ge c,\\
		\\
		\dfrac{r}{c}&,\textup{ if } 0\le r< c,\\
		\\
		0&,\textup{ if } r<0.
	\end{cases}
\end{equation*}

Given $c>0$ and given $\Psi\in\mathcal{V}(K)$, define $\gls{beta2c}:\bm{V}(\Omega^{(2)})\to\bm{V}(\Omega^{(2)})$ by:
\begin{equation*}
	\bm{\beta}_{2,c}(\bm{v};\Psi):=b_{2,c}((\mathbb{I}_{\Omega^{(2)}}+\bm{v})\cdot\bm{e}_3-\Psi)\bm{e}_3.
\end{equation*}

Since $b_{2,c}$ is Lipschitz continuous in $\bm{L}^2(\Omega^{(2)})$ with Lipschitz constant equal to one and monotone, then $\bm{\beta}_{2,c}$ is Lipschitz continuous with Lipschitz constant equal to one and monotone.
Besides, it results that $0\le b_{2,c}'\le 1$ independently of $c$. We are now ready to state the regularised version of Problem~\ref{problem2rho}.
\begin{customprob}{$\mathcal{P}_{\rho,c}(\Omega^{(2)};\Psi)$}\label{problem2rhoc}
	Find $\bm{u}_{\rho,c}^{(2)}(\Psi)=(u_{\rho,c,i}^{(2)}(\Psi))\in\bm{V}(\Omega^{(2)})$ that satisfies the variational equations
	\begin{equation*}
		\int_{\Omega^{(2)}} A^{ijk\ell,2} e_{k\|\ell}(\bm{u}_{\rho,c}^{(2)}(\Psi)) e_{i\|j}(\bm{v}) \,\mathrm{d}x^{(2)}+\dfrac{1}{\rho}\int_{\Omega^{(2)}}\bm{\beta}_{2,c}(\bm{u}_{\rho,c}^{(2)}(\Psi);\Psi)\cdot\bm{v}\,\mathrm{d}x^{(2)} = \int_{\Omega^{(2)}}\bm{f}^{(2)}\cdot\bm{v}\,\mathrm{d}x^{(2)},
	\end{equation*}
	for all $\bm{v}=(v_i) \in \bm{V}(\Omega^{(2)})$.
	\bqed
\end{customprob}

Denote by $B_{\alpha,c}\in\mathcal{C}^2(\mathbb{R})$ the antiderivative of $b_{\alpha,c}$ satisfying $B_{\alpha,c}(0)=0$. More, specifically, we have
\begin{equation*}
	B_{1,c}(r):=
	\begin{cases}
		\dfrac{r^2}{2}+\dfrac{c}{2}r+\dfrac{c^2}{6}&,r\le -c,\\
		\\
		-\dfrac{r^3}{6c}&,-c<r\le 0,\\
		\\
		0&,r>0,
	\end{cases}
\end{equation*}
and
\begin{equation*}
	B_{2,c}(r):=
	\begin{cases}
		\dfrac{r^2}{2}-\dfrac{c}{2}r+\dfrac{c^2}{6}&,r\ge c,\\
		\\
		\dfrac{r^3}{6c}&,0<r\le c,\\
		\\
		0&,r<0.
	\end{cases}
\end{equation*}

It is straightforward to observe that $B_{\alpha,c}\ge0$, $B_{\alpha,c}''=b_{\alpha,c}'\ge0$, and that:
\begin{equation}
	\label{conv:smooth-penalty}
	\begin{aligned}
		B_{1,c}(r)\to\dfrac{1}{2}(\min\{r,0\})^2,\quad\textup{ uniformly in }\mathbb{R}\textup{ as }c\to0,\\[6pt]
		B_{2,c}(r)\to\dfrac{1}{2}(\max\{r,0\})^2,\quad\textup{ uniformly in }\mathbb{R}\textup{ as }c\to0.
	\end{aligned}
\end{equation}

Let us consider the regularised version of Problem~\ref{cprho}, that we denote by Problem $\mathscr{C}_{\rho,c}(K)$.
\begin{customprob}{$\mathscr{C}_{\rho,c}(K)$}\label{cprhoc}
	Find $\Psi\in\mathcal{V}(K)$ that minimises the reduced cost functional
	\begin{equation*}
		\begin{aligned}
			&\gls{Jrhoc}(\Psi):=\dfrac{\beta}{2}\int_{K}|\Delta\Psi|^2\,\mathrm{d}x
			+\dfrac{1}{2}\int_{\Omega^{(1)}}A^{ijk\ell,1}e_{k\|\ell}(\bm{u}_{\rho,c}^{(1)}(\Psi))e_{i\|j}(\bm{u}_{\rho,c}^{(1)}(\Psi))\,\mathrm{d}x\\
			&\quad+\dfrac{1}{\rho}\int_{\Omega^{(1)}}B_{1,c}((\mathbb{I}_{\Omega^{(1)}}+\bm{u}_{\rho,c}^{(1)}(\Psi))\cdot\bm{e}_3-\Psi)\,\mathrm{d}x-\int_{\Omega^{(1)}}\bm{f}^{(1)}\cdot\bm{u}_{\rho,c}^{(1)}(\Psi)\,\mathrm{d}x\\
			&\quad+\dfrac{1}{2}\int_{\Omega^{(2)}}A^{ijk\ell,2}e_{k\|\ell}(\bm{u}_{\rho,c}^{(2)}(\Psi))e_{i\|j}(\bm{u}_{\rho,c}^{(2)}(\Psi))\,\mathrm{d}x+\dfrac{1}{\rho}\int_{\Omega^{(2)}}B_{2,c}((\mathbb{I}_{\Omega^{(2)}}+\bm{u}_{\rho,c}^{(2)}(\Psi))\cdot\bm{e}_3-\Psi)\,\mathrm{d}x\\
			&\quad-\int_{\Omega^{(2)}}\bm{f}^{(2)}\cdot\bm{u}_{\rho,c}^{(2)}(\Psi)\,\mathrm{d}x,
		\end{aligned}
	\end{equation*}
	subject to $\bm{u}_{\rho,c}^{(1)}(\Psi)$ being the solution of Problem~\ref{problem1rhoc} and to $\bm{u}_{\rho,c}^{(2)}(\Psi)$ being the solution of Problem~\ref{problem2rhoc}.
	\bqed
\end{customprob}

In the same fashion as Theorem~\ref{ex-un-op-con-rho}, we are able to show (Theorem~\ref{th:optimal-1-c}) that Problem~\ref{cprhoc} admits one and only one solution, and that minimisers of Problem~\ref{cprhoc} converge, in a suitable sense, to minimisers of Problem~\ref{cprho} as $c\to0$.
\begin{theorem}
	\label{th:optimal-1-c}
	Let $\rho>0$ and $\bm{f}^{(\alpha)}\in\bm{L}^2(\Omega^{(\alpha)})$ be given. Then, the following properties hold:
	\begin{itemize}
		\item[$(a)$] For each $c>0$, there exists one and only one solution $\Psi_{\rho,c}\in\mathcal{V}(K)$ for Problem~\ref{cprhoc};
		\item[$(b)$] Up to passing to a subsequence, it results that $\Psi_{\rho,c}\rightharpoonup\Psi_\rho$ in $H^2(K)$ as $c\to0$, where $\Psi_\rho$ is the unique solution of Problem~\ref{cprho}.
	\end{itemize}
\end{theorem}
\begin{proof}
	The proof is divided into three steps, numbered $(i)$--$(iii)$.
	
	$(i)$ \emph{Proof of $(a)$.} Consider the \emph{full functional} $\widehat{\mathcal{J}}_{\rho,c}:\mathcal{V}(K)\times\bm{V}(\Omega^{(1)})\times\bm{V}(\Omega^{(2)})\to\mathbb{R}$ defined by:
	\begin{equation*}
		\begin{aligned}
			&\widehat{\mathcal{J}}_{\rho,c}(\Psi,\bm{v}^{(1)},\bm{v}^{(2)}):=\dfrac{\beta}{2}\int_K|\Delta\Psi|^2\,\mathrm{d}x+\dfrac{1}{2}\int_{\Omega^{(1)}}A^{ijk\ell,1}e_{k\|\ell}(\bm{v}^{(1)})e_{i\|j}(\bm{v}^{(1)})\,\mathrm{d}x\\
			&\quad+\dfrac{1}{\rho}\int_{\Omega^{(1)}}B_{1,c}((\mathbb{I}_{\Omega^{(1)}}+\bm{v}^{(1)})\cdot\bm{e}_3-\Psi)\,\mathrm{d}x-\int_{\Omega^{(1)}}\bm{f}^{(1)}\cdot\bm{v}^{(1)}\,\mathrm{d}x\\
			&\quad+\dfrac{1}{2}\int_{\Omega^{(2)}}A^{ijk\ell,2}e_{k\|\ell}(\bm{v}^{(2)})e_{i\|j}(\bm{v}^{(2)})\,\mathrm{d}x+\dfrac{1}{\rho}\int_{\Omega^{(2)}}B_{2,c}((\mathbb{I}_{\Omega^{(2)}}+\bm{v}^{(2)})\cdot\bm{e}_3-\Psi)\,\mathrm{d}x\\
			&\quad-\int_{\Omega^{(2)}}\bm{f}^{(2)}\cdot\bm{v}^{(2)}\,\mathrm{d}x.
		\end{aligned}
	\end{equation*}
	
	The strict convexity and coerciveness of the functional $\widehat{\mathcal{J}}_{\rho,c}$ ensures that $\widehat{\mathcal{J}}_{\rho,c}$ admits a unique minimiser, that we denote by $(\Psi_{\rho,c},\bm{u}^{(1)}_{\rho,c},\bm{u}^{(2)}_{\rho,c})$.
	
	For each fixed $\Psi\in\mathcal{V}(K)$, the full functional $\widehat{\mathcal{J}}_{\rho,c}(\Psi,\cdot,\cdot)$ splits as a sum of two \emph{independent parts}, one depending only on $\bm{v}^{(1)}$ and the other depending only on $\bm{v}^{(2)}$. Minimising each part separately gives exactly Problem~\ref{problem1rhoc} and Problem~\ref{problem2rhoc}, since the penalty satisfies $B_{\alpha,c}'=b_{\alpha,c}$ in $\mathbb{R}$. Therefore, it results $\bm{u}^{(\alpha)}_{\rho,c}=\bm{u}^{(\alpha)}_{\rho,c}(\Psi_{\rho,c})$ and the control $\Psi_{\rho,c}$ thus minimises $\mathcal{J}_{\rho,c}$.
	
	$(ii)$ \emph{Let $\{\Psi_{\rho,c}\}_{c>0}$ be a sequence of minimisers for Problem~\ref{cprhoc}. Let $\bm{u}^{(\alpha)}_{\rho,c}(\Psi_{\rho,c})$ be the unique solution of either Problem~\ref{problem1rhoc} or Problem~\ref{problem2rhoc}. Then, up to passing to a subsequence, it results that $\{\bm{u}^{(\alpha)}_{\rho,c_n}(\Psi_{\rho,c_n})\}_{n=1}^\infty$ weakly converges in $\bm{V}(\Omega^{(\alpha)})$ to the unique solution of either Problem~\ref{problem1rho} or Problem~\ref{problem2rho}}. Arguing exactly as in~\eqref{eq:seb} in item $(iii)$ of the proof of Theorem~\ref{th:optimal-1}, with the comparison element $\bm{v}=\bm{0}$, which is admissible since~\eqref{eq:psi} gives $x_3\ge0$ in $\overline{\Omega^{(1)}}$ and $x_3\le0$ in $\overline{\Omega^{(2)}}$, we obtain that $\mathcal{J}_{\rho,c}(0)$ is bounded independently of $c$. Since $\Psi_{\rho,c}$ minimises $\mathcal{J}_{\rho,c}$, it results $\mathcal{J}_{\rho,c}(\Psi_{\rho,c})\le\mathcal{J}_{\rho,c}(0)\le C$, for some $C>0$ independent of $c$ but possibly dependent on $\rho$.
	Therefore, up to passing to a subsequence $\{c_n\}_{n=1}^\infty$ of positive numbers such that $c_n \to 0$ as $n\to\infty$, we obtain:
	\begin{equation}
		\label{psi:conv}
		\Psi_{\rho,c_n}\rightharpoonup\Psi_\rho,\quad\textup{ in }H^2(K)\textup{ as }n\to\infty.
	\end{equation}
	
	Proceeding as in the derivation of~\eqref{eq:bdd:sol-c}, and using that~\eqref{psi:conv} renders the competitor $\bm{r}_{\rho,c}^{(1)}:=\{(\mathbb{I}_{\Omega^{(1)}}-\Psi_{\rho,c}\bm{e}_3)\cdot\bm{e}_3\}^{-}\bm{e}_3$ bounded in $\bm{V}(\Omega^{(1)})$ independently of $c$ (and similarly for
	$\alpha=2$), the energy comparison gives the uniform bound that $\{\bm{u}^{(\alpha)}_{\rho,c}(\Psi_{\rho,c})\}_{c>0}$ is bounded in $\bm{V}(\Omega^{(\alpha)})$ independently of $c$. Therefore, we can extract a weakly convergent subsequence, still denoted $\{\bm{u}^{(\alpha)}_{\rho,c_n}(\Psi_{\rho,c_n})\}_{n=1}^\infty$, such that $\bm{u}^{(\alpha)}_{\rho,c_n}(\Psi_{\rho,c_n})\rightharpoonup\bm{u}^{(\alpha)}_{\rho}$ in $\bm{V}(\Omega^{(\alpha)})$ as $n\to\infty$. Exploiting the continuity of the operator $\mathcal{A}$ (cf.~\eqref{A}) and item~$(ii)$ in Theorem~\ref{th:exun:reg}, we obtain that the weak limit coincides with the unique solution of either Problem~\ref{problem1rho} or Problem~\ref{problem2rho}.
	
	$(iii)$ \emph{Completion of the proof.} Following the same strategy as in item $(iii)$ of the proof of Theorem~\ref{th:optimal-1}, with the convergence $\bm{u}_{\rho,c}^{(\alpha)}\to\bm{u}_\rho^{(\alpha)}$ in $\bm{V}(\Omega^{(\alpha)})$ as $c\to 0$ of Theorem~\ref{th:exun:reg} in place of Lemma~\ref{lem:conv}, and using~\eqref{conv:smooth-penalty} gives
	\begin{equation*}
		\mathcal{J}_\rho(\Psi_\rho)\le\liminf_{c\to0}\mathcal{J}_{\rho,c}(\Psi_{\rho,c})\le\lim_{c\to0}\mathcal{J}_{\rho,c}(\Phi)=\mathcal{J}_\rho(\Phi),
	\end{equation*}
	for all $\Phi\in\mathcal{V}(K)$, thus completing the proof.
\end{proof}

The purpose of the forthcoming calculations is to establish the differentiability of the reduced cost functional $\mathcal{J}_{\rho,c}$. To begin with, we show that the mapping
\begin{equation*}
	\Psi\in\mathcal{V}(K)\mapsto\bm{u}^{(1)}_{\rho,c}(\Psi)
\end{equation*}
if G\^{a}teaux differentiable.

\begin{lemma}
	\label{lem:Gateaux}
	Let $\rho>0$, $c>0$ and $\Psi,\Phi\in\mathcal{V}(K)$ be given, and set $\tilde{\Psi}:=\Phi-\Psi$. Since $\mathcal{V}(K)$ is convex, we have $(\Psi+\varepsilon\tilde{\Psi})\in\mathcal{V}(K)$ for every $\varepsilon\in[0,1]$. Then the mapping
	\begin{equation*}
		\Psi\in\mathcal{V}(K)\mapsto\bm{u}^{(\alpha)}_{\rho,c}(\Psi)
	\end{equation*}
	is weakly G\^{a}teaux differentiable in $\mathcal{V}(K)$, and the action of weak derivative at $\Psi$ along the direction $\tilde{\Psi}\in\mathcal{V}(K)$ is denoted by $\bm{U}_{\rho,c}^{(\alpha)}[\tilde{\Psi}]$.
\end{lemma}
\begin{proof}
	Let $\varepsilon>0$ and let $\bm{u}^{(\alpha)}_{\rho,c}(\Psi)$ and $\bm{u}^{(\alpha)}_{\rho,c}(\Psi+\varepsilon\tilde{\Psi})$ be the unique states corresponding to $\Psi$ and $(\Psi+\varepsilon\tilde{\Psi})$, respectively. It suffices to consider the case $\alpha=1$, being the proof for the case $\alpha=2$ analogous. Subtracting the variational equations of Problem~\ref{problem1rhoc} associated with $\Psi$ tested at $(\bm{u}^{(1)}_{\rho,c}(\Psi+\varepsilon\tilde{\Psi})-\bm{u}^{(1)}_{\rho,c}(\Psi))$ from the variational equations of Problem~\ref{problem1rhoc} associated with $(\Psi+\varepsilon\tilde{\Psi})$ tested at $(\bm{u}^{(1)}_{\rho,c}(\Psi+\varepsilon\tilde{\Psi})-\bm{u}^{(1)}_{\rho,c}(\Psi))$, we obtain:
	\begin{equation}
		\label{eq:diff:step:1}
		\begin{aligned}
			&\int_{\Omega^{(1)}}A^{ijk\ell,1}e_{k\|\ell}(\bm{u}^{(1)}_{\rho,c}(\Psi+\varepsilon\tilde{\Psi})-\bm{u}^{(1)}_{\rho,c}(\Psi))e_{i\|j}(\bm{u}^{(1)}_{\rho,c}(\Psi+\varepsilon\tilde{\Psi})-\bm{u}^{(1)}_{\rho,c}(\Psi))\,\mathrm{d}x^{(1)}\\
			&\quad+\dfrac{1}{\rho}\int_{\Omega^{(1)}}\left(\bm{\beta}_{1,c}(\bm{u}^{(1)}_{\rho,c}(\Psi+\varepsilon\tilde{\Psi});\Psi+\varepsilon\tilde{\Psi})-\bm{\beta}_{1,c}(\bm{u}^{(1)}_{\rho,c}(\Psi);\Psi)\right)\cdot(\bm{u}^{(1)}_{\rho,c}(\Psi+\varepsilon\tilde{\Psi})-\bm{u}^{(1)}_{\rho,c}(\Psi))\,\mathrm{d}x^{(1)}=0.
		\end{aligned}
	\end{equation}
	
	The monotonicity of $b_{1,c}$ gives
	\begin{equation}
		\label{eq:diff:step:2-1}
		\int_{\Omega^{(1)}}\left(\bm{\beta}_{1,c}(\bm{u}^{(1)}_{\rho,c}(\Psi+\varepsilon\tilde{\Psi});\Psi+\varepsilon\tilde{\Psi})-\bm{\beta}_{1,c}(\bm{u}^{(1)}_{\rho,c}(\Psi);\Psi+\varepsilon\tilde{\Psi})\right)\cdot(\bm{u}^{(1)}_{\rho,c}(\Psi+\varepsilon\tilde{\Psi})-\bm{u}^{(1)}_{\rho,c}(\Psi))\,\mathrm{d}x^{(1)}\ge 0,
	\end{equation}
	while the Lipschitz continuity of $b_{1,c}$ with Lipschitz constant equal to one gives:
	\begin{equation}
		\label{eq:diff:step:2}
		\begin{aligned}
			&|\bm{\beta}_{1,c}(\bm{u}^{(1)}_{\rho,c}(\Psi);\Psi+\varepsilon\tilde{\Psi})-\bm{\beta}_{1,c}(\bm{u}^{(1)}_{\rho,c}(\Psi);\Psi)|\\
			&=|b_{1,c}((\mathbb{I}_{\Omega^{(1)}}+\bm{u}^{(1)}_{\rho,c}(\Psi))\cdot\bm{e}_3-(\Psi+\varepsilon\tilde{\Psi}))-b_{1,c}((\mathbb{I}_{\Omega^{(1)}}+\bm{u}^{(1)}_{\rho,c}(\Psi))\cdot\bm{e}_3-\Psi)|\le\varepsilon|\tilde{\Psi}|.
		\end{aligned}
	\end{equation}

	Dividing~\eqref{eq:diff:step:1} by $\varepsilon^2$ and exploiting~\eqref{eq:diff:step:2-1} and~\eqref{eq:diff:step:2} gives:
	\begin{equation}
		\label{eq:diff:step:3}
		\begin{aligned}
			&-\dfrac{1}{\rho}\int_{\Omega^{(1)}}\dfrac{\bm{\beta}_{1,c}(\bm{u}^{(1)}_{\rho,c}(\Psi+\varepsilon\tilde{\Psi});\Psi+\varepsilon\tilde{\Psi})-\bm{\beta}_{1,c}(\bm{u}^{(1)}_{\rho,c}(\Psi);\Psi)}{\varepsilon}\cdot\dfrac{\bm{u}^{(1)}_{\rho,c}(\Psi+\varepsilon\tilde{\Psi})-\bm{u}^{(1)}_{\rho,c}(\Psi)}{\varepsilon}\,\mathrm{d}x^{(1)}\\
			&\le\dfrac{\|\tilde{\Psi}\|_{L^2(K)}}{\rho}\left\|\dfrac{\bm{u}^{(1)}_{\rho,c}(\Psi+\varepsilon\tilde{\Psi})-\bm{u}^{(1)}_{\rho,c}(\Psi)}{\varepsilon}\right\|_{\bm{V}(\Omega^{(1)})}.
		\end{aligned}
	\end{equation}
	
	Combining~\eqref{eq:diff:step:3} with~\eqref{eq:diff:step:1}, the uniform positive-definiteness of the fourth order three-dimensional elasticity tensor $\{A^{ijk\ell,1}\}$, Korn's inequality (Theorem~\ref{KornCart}) and Young's inequality~\cite{Young1912} gives
	\begin{equation*}
		\left\|\dfrac{\bm{u}^{(1)}_{\rho,c}(\Psi+\varepsilon\tilde{\Psi})-\bm{u}^{(1)}_{\rho,c}(\Psi)}{\varepsilon}\right\|_{\bm{V}(\Omega^{(1)})}\le\dfrac{C_e^{(1)}|c_0^{(1)}|^2}{\rho}\|\tilde{\Psi}\|_{L^2(K)},
	\end{equation*}
	thus showing that the incremental ratio for $\bm{u}^{(1)}_{\rho,c}$ at $\Psi$ in the direction $\tilde{\Psi}$ is independent of $c$ and $\varepsilon$. Therefore, up to passing to a subsequence $\{\varepsilon_n\}_{n=1}^\infty$ of positive numbers such that $\varepsilon_n \to 0$ as $n\to\infty$, we obtain that
	\begin{equation*}
		\dfrac{\bm{u}^{(1)}_{\rho,c}(\Psi+\varepsilon_n\tilde{\Psi})-\bm{u}^{(1)}_{\rho,c}(\Psi)}{\varepsilon_n}\rightharpoonup\bm{U}_{\rho,c}^{(1)}[\tilde{\Psi}],\quad\textup{ in }\bm{V}(\Omega^{(1)})\textup{ as }n\to\infty.
	\end{equation*}
\end{proof}

We are finally in position to establish the G\^{a}teaux differentiability of the reduced cost functional $\mathcal{J}_{\rho,c}$ and compute its G\^{a}teaux derivative explicitly.
\begin{theorem}
	\label{th:Gateaux:Jrhoc}
	Let $\rho>0$, $c>0$ and $\Psi,\Phi\in\mathcal{V}(K)$ be given, and set $\tilde{\Psi}:=\Phi-\Psi$. Then, the mapping
	\begin{equation*}
		\Psi\in\mathcal{V}(K)\mapsto\mathcal{J}_{\rho,c}(\Psi)
	\end{equation*}
	is G\^{a}teaux differentiable in $\mathcal{V}(K)$, and the action of the directional derivative $\mathcal{J}_{\rho,c}'(\Psi)$ along the direction $\tilde{\Psi}\in\mathcal{V}(K)$ is given by
	\begin{equation*}
		\mathcal{J}_{\rho,c}'(\Psi)[\tilde{\Psi}]=\beta\int_K\Delta\Psi \Delta\tilde{\Psi}\,\mathrm{d}x-\dfrac{1}{\rho}\int_{\Omega^{(1)}}b_{1,c}(s_1)\tilde{\Psi}\,\mathrm{d}x-\dfrac{1}{\rho}\int_{\Omega^{(2)}}b_{2,c}(s_2)\tilde{\Psi}\,\mathrm{d}x,
	\end{equation*}
	where $s_\alpha:=(\mathbb{I}_{\Omega^{(\alpha)}}+\bm{u}_{\rho,c}^{(\alpha)}(\Psi))\cdot\bm{e}_3-\Psi$.
\end{theorem}
\begin{proof}
	Observing that the full cost functional $\widehat{\mathcal{J}}_{\rho,c}$ is of class $\mathcal{C}^1(\mathcal{V}(K)\times\bm{V}(\Omega^{(1)})\times\bm{V}(\Omega^{(2)}))$ and recalling that the control-to-state mapping $\Psi\in\mathcal{V}(K)\mapsto\bm{u}_{\rho,c}^{(\alpha)}(\Psi)$ is G\^{a}teaux differentiable (Lemma~\ref{lem:Gateaux}), an application of the mixed chain rule (cf., e.g., Theorem~3.2.12 in~\cite{DM13} or page~50 in~\cite{HPUU18}) gives:
	\begin{equation}
		\label{chain-rule}
		\begin{aligned}
			&\mathcal{J}_{\rho,c}'(\Psi)[\tilde{\Psi}]=\partial_\Psi\widehat{\mathcal{J}}_{\rho,c}\left(\Psi,\bm{u}_{\rho,c}^{(1)}(\Psi),\bm{u}_{\rho,c}^{(2)}(\Psi)\right)[\tilde{\Psi}]+\partial_{\bm{v}^{(1)}}\widehat{\mathcal{J}}_{\rho,c}\left(\Psi,\bm{u}_{\rho,c}^{(1)}(\Psi),\bm{u}_{\rho,c}^{(2)}(\Psi)\right)\left(\bm{U}_{\rho,c}^{(1)}[\tilde{\Psi}]\right)\\
			&\quad+\partial_{\bm{v}^{(2)}}\widehat{\mathcal{J}}_{\rho,c}\left(\Psi,\bm{u}_{\rho,c}^{(1)}(\Psi),\bm{u}_{\rho,c}^{(2)}(\Psi)\right)\left(\bm{U}_{\rho,c}^{(2)}[\tilde{\Psi}]\right).
		\end{aligned}
	\end{equation}
	
	For each fixed $\Psi\in\mathcal{V}(K)$, the full functional $\widehat{\mathcal{J}}_{\rho,c}(\Psi,\cdot,\cdot)$ splits as a sum of two independent parts, one depending only on $\bm{v}^{(1)}$ and the other depending only on $\bm{v}^{(2)}$. Since $B_{\alpha,c}'=b_{\alpha,c}$ in $\mathbb{R}$, the Euler--Lagrange equations of each independent part of $\widehat{\mathcal{J}}_{\rho,c}(\Psi,\cdot,\cdot)$ are exactly Problem~\ref{problem1rhoc}, if $\alpha=1$, and Problem~\ref{problem2rhoc}, if $\alpha=2$. Hence the minimiser $\bm{u}_{\rho,c}^{(\alpha)}(\Psi)$ satisfies:
	\begin{equation*}
		\begin{aligned}
			\partial_{\bm{v}^{(1)}}\widehat{\mathcal{J}}_{\rho,c}\left(\Psi,\bm{u}_{\rho,c}^{(1)}(\Psi),\bm{u}_{\rho,c}^{(2)}(\Psi)\right)&=\bm{0}\in\bm{V}(\Omega^{(1)})^\ast\textup{ with }\Psi\textup{ and }\bm{u}_{\rho,c}^{(2)}(\Psi)\textup{ fixed},\\
			\partial_{\bm{v}^{(2)}}\widehat{\mathcal{J}}_{\rho,c}\left(\Psi,\bm{u}_{\rho,c}^{(1)}(\Psi),\bm{u}_{\rho,c}^{(2)}(\Psi)\right)&=\bm{0}\in\bm{V}(\Omega^{(2)})^\ast\textup{ with }\Psi\textup{ and }\bm{u}_{\rho,c}^{(1)}(\Psi)\textup{ fixed}.
		\end{aligned}
	\end{equation*}
	
	Therefore, the expression for $\mathcal{J}_{\rho,c}'(\Psi)[\tilde{\Psi}]$ in~\eqref{chain-rule} reduces to:
	\begin{equation*}
		\mathcal{J}_{\rho,c}'(\Psi)[\tilde{\Psi}]=\partial_\Psi\widehat{\mathcal{J}}_{\rho,c}\left(\Psi,\bm{u}_{\rho,c}^{(1)}(\Psi),\bm{u}_{\rho,c}^{(2)}(\Psi)\right)[\tilde{\Psi}].
	\end{equation*}
	
	Observe that
	\begin{equation*}
		\partial_\Psi\left(\dfrac{\beta}{2}\|\Delta\Psi\|^2_{L^2(K)}\right)[\tilde\Psi]=\beta\int_{K}\Delta\Psi\Delta\tilde\Psi\,\mathrm{d}x,
	\end{equation*}
	and that, since $B_{\alpha,c}'=b_{\alpha,c}$ in $\mathbb{R}$, and $s_\alpha=(\mathbb{I}_{\Omega^{(\alpha)}}+\bm{u}_{\rho,c}^{(\alpha)}(\Psi))\cdot\bm{e}_3-\Psi$, it results:
	\begin{equation*}
		\partial_\Psi B_{\alpha,c}(s_\alpha)[\tilde{\Psi}]=-b_{\alpha,c}(s_\alpha)\tilde{\Psi}.
	\end{equation*}
	
	
	Therefore, we obtain the sought formula
	\begin{equation*}
		\mathcal{J}_{\rho,c}'(\Psi)[\tilde{\Psi}]=\beta\int_K\Delta\Psi \Delta\tilde{\Psi}\,\mathrm{d}x-\dfrac{1}{\rho}\int_{\Omega^{(1)}}b_{1,c}(s_1)\tilde{\Psi}\,\mathrm{d}x-\dfrac{1}{\rho}\int_{\Omega^{(2)}}b_{2,c}(s_2)\tilde{\Psi}\,\mathrm{d}x,
	\end{equation*}
	and the proof is complete.
\end{proof}

Since $\Psi\in\mathcal{V}(K)$ allows us to write $\Psi=\psi\circ\pi_\kappa$ for some $\psi\in H^2_0(\kappa)$, an application of Theorem~\ref{th:Gateaux:Jrhoc} implies that the action of the G\^{a}teaux derivative of the reduced cost functional $\mathcal{J}_{\rho,c}$ on a direction $\tilde{\Psi}=\tilde{\psi}\circ\pi_\kappa$ takes the form
\begin{equation*}
		\mathcal{J}_{\rho,c}'(\Psi)[\tilde{\Psi}]=2a\beta\int_{\kappa}\Delta\psi\Delta\tilde{\psi}\,\mathrm{d}\xi-\int_{\kappa}\left(\Lambda^{(1)}_{\rho,c}(\xi)+\Lambda^{(2)}_{\rho,c}(\xi)\right)\tilde{\psi}(\xi)\,\mathrm{d}\xi,
\end{equation*}
where, for $\xi=(\xi_1,\xi_2)\in\kappa$ and $\alpha\in\{1,2\}$, the \emph{fibre-averaged multipliers} $\Lambda^{(\alpha)}_{\rho,c}\in L^1(\kappa)$ are defined by
\begin{equation}
	\label{lambda-fibres}
	\begin{aligned}
		\Lambda^{(1)}_{\rho,c}(\xi)&:=
		\begin{cases}
			\dfrac{1}{\rho}\displaystyle\int_{F^{(1)}(\xi)}b_{1,c}\left(s_1(\xi,t)\right)\,\mathrm{d}t&,\quad\textup{ for a.a. } \xi\in\pi_\kappa(\Omega^{(1)}),\\[8pt]
			0&,\quad\textup{ for a.a. } \xi\in\kappa\setminus\pi_\kappa(\Omega^{(1)}),
		\end{cases}\\[14pt]
		\Lambda^{(2)}_{\rho,c}(\xi)&:=
		\begin{cases}
			\dfrac{1}{\rho}\displaystyle\int_{F^{(2)}(\xi)}b_{2,c}\left(s_2(\xi,t)\right)\,\mathrm{d}t&,\quad\textup{ for a.a. } \xi\in\pi_\kappa(\Omega^{(2)}),\\[8pt]
			0&,\quad\textup{ for a.a. } \xi\in\kappa\setminus\pi_\kappa(\Omega^{(2)}),
		\end{cases}
	\end{aligned}
\end{equation}
where
\begin{equation*}
	s_\alpha(\xi,t):=\left(\mathbb{I}_{\Omega^{(\alpha)}}(\xi,t)+\bm{u}^{(\alpha)}_{\rho,c}(\Psi)(\xi,t)\right)\cdot\bm{e}_3-\psi(\xi),\quad\textup{ for a.a. }(\xi,t)\in\Omega^{(\alpha)}.
\end{equation*}

Note that the fibre $F^{(\alpha)}(\cdot)$ is used in two different ways. On the one hand, the common overlap $\omega$ it encodes the two-body non-penetration constraint. On the other hand, on the full projection $\pi_\kappa(\Omega^{(\alpha)})$ it encodes the single-body penalty in $\mathcal{J}_{\rho,c}$. Consequently, the \emph{formal} distributional derivative of $\mathcal{J}_{\rho,c}$ reads:
\begin{equation*}
	\mathcal{J}_{\rho,c}'(\Psi)=2a\beta\Delta^2\psi-\Lambda^{(1)}_{\rho,c}-\Lambda^{(2)}_{\rho,c},\quad\textup{ in }\mathcal{D}'(\kappa).
\end{equation*}

\section{Numerical Verification}\label{Sec:5}
In the final section of this paper, we exploit the G\^{a}teaux differentiability result established in Theorem~\ref{th:Gateaux:Jrhoc} and the uniqueness of the optimal control for the cost functional $\mathcal{J}_{\rho,c}$ established in Theorem~\ref{th:optimal-1-c} to formulate a first order necessary condition for an element $\Psi_{\rho,c}$ to minimise the cost functional $\mathcal{J}_{\rho,c}$ as well as an algorithm allowing to implement this minimisation process numerically. The algorithm we plan on implementing is the following.

\begin{algorithm}[H]
	\caption{Optimal Interface Method}\label{alg:1}
	\begin{algorithmic}[1]
		\Require $\bm{f}^{(\alpha)}$, $\beta$, $\lambda$, $\mu$, $\varepsilon$, $k_{\max}$, initial guess $\psi^{(0)}\in H^2_0(\kappa)$
		\State Set $\Psi^{(0)}:=\psi^{(0)}\circ\pi_\kappa$, $k:=0$
		\While{$k<k_{\max}$ \textbf{and} $\|\mathcal{J}_{\rho,c}'(\Psi^{(k)})\|>\varepsilon$}
		\State Solve Problem~\ref{problem1rhoc} with $\Psi=\Psi^{(k)}$ to get $\bm{u}^{(1)}_{\rho,c}(\Psi^{(k)})$
		\State Solve Problem~\ref{problem2rhoc} with $\Psi=\Psi^{(k)}$ to get $\bm{u}^{(2)}_{\rho,c}(\Psi^{(k)})$
		\State Compute the fibre-averaged multipliers $\Lambda^{(1)}_{\rho,c},\Lambda^{(2)}_{\rho,c}$ via \eqref{lambda-fibres}
		\State Compute the descent direction
		$G^{(k)}:=2a\beta\,\Delta^2\psi^{(k)}-\Lambda^{(1)}_{\rho,c}-\Lambda^{(2)}_{\rho,c}$ (in $\mathcal{D}'(\kappa)$)
		\State Apply pre-conditioning and determine the step-size $\tau_k>0$
		\State Update $\psi^{(k+1)}:=\psi^{(k)}-\tau_k\,G^{(k)}$
		\State Set $\Psi^{(k+1)}:=\psi^{(k+1)}\circ\pi_\kappa$, $k:=k+1$
		\EndWhile
		\State \Return $\Psi^{(k)}$ (approximate optimal interface)
	\end{algorithmic}
\end{algorithm}

We test the algorithm on a benchmark problem posed over a simplified geometrical framework, and for which an explicit solution is known. More specifically, we consider a scalar problem, for which the displacements of the linearly elastic membranes are both defined over the same two-dimensional Lipschitz domain. Note that this simplification in the dimension of the problem allows us to bypass the challenges introduced by the fibre-based argument discussed in the previous sections.

Let $\omega=B(0;1)$ be the unit ball in $\mathbb{R}^2$, let $f^{(1)}, f^{(2)}\in L^2(\omega)$ denote the applied body forces, and let $\theta_1,\theta_2\in\mathcal{C}^0(\overline{\omega})$ denote the references configurations of the two linearly elastic bodies under consideration, which are assumed to be such that $\theta_1(y)>0>\theta_2(y)$ for all $y=(y_\alpha)\in\overline{\omega}$. More precisely, we take $\theta_1\equiv 0.45$ and $\theta_2\equiv-0.05$. The set
\begin{equation*}
	\vec{U}:=\{(v_1,v_2)\in H^1_0(\omega)\times H^1_0(\omega);\theta_1(y)+v_1(y)\ge\theta_2(y)+v_2(y),\textup{ for a.a. }y\in\omega\}
\end{equation*}
is non-empty, closed and convex. We consider the following variational problem: \emph{Find $(u^{(1)},u^{(2)})\in\vec{U}$ satisfying the variational inequalities
\begin{equation*}
	\int_{\omega}\nabla u^{(1)}\cdot\nabla(v-u^{(1)}) \,\mathrm{d}y+\int_{\omega}\nabla u^{(2)}\cdot\nabla(v-u^{(2)}) \,\mathrm{d}y\ge \int_{\omega} f^{(1)} (v_1-u^{(1)}) \,\mathrm{d}y+\int_{\omega} f^{(2)}(v_2-u^{(2)})\,\mathrm{d}y,
\end{equation*}
for all $(v_1,v_2)\in\vec{U}$.}

It can be immediately observed that this variational problem admits a unique solution $(u^{(1)},u^{(2)})\in\vec{U}$ thanks to the Lions--Stampacchia theorem (cf., e.g., Theorem~4.8-1 in~\cite{Cia25}). For each $\psi\in H^2_0(\omega)$, define the sets:
\begin{equation*}
	\begin{aligned}
		K^{(1)}(\omega;\psi)&:=\{v_1\in H^1_0(\omega);\theta_1(y)+v_1(y)\ge\psi(y),\textup{ for a.a. }y\in\omega\},\\
		K^{(2)}(\omega;\psi)&:=\{v_2\in H^1_0(\omega);\theta_2(y)+v_2(y)\le\psi(y),\textup{ for a.a. }y\in\omega\}.
	\end{aligned}
\end{equation*}

Given any $\beta>0$, the corresponding optimal control approximation takes the following form: \emph{Find $\psi_\beta\in H^2_0(\omega)$ that minimises the reduced cost functional
	\begin{equation*}
		\mathcal{J}_\beta(\psi):=\dfrac{\beta}{2}\int_{\omega}|\Delta\psi|^2\,\mathrm{d}y
		+\dfrac{1}{2}\int_{\omega}|\nabla u^{(1)}(\psi)|^2\,\mathrm{d}y-\int_{\omega}f^{(1)} u^{(1)}(\psi)\,\mathrm{d}y
		+\dfrac{1}{2}\int_{\omega}|\nabla u^{(2)}(\psi)|^2\,\mathrm{d}y-\int_{\omega}f^{(2)} u^{(2)}(\psi)\,\mathrm{d}y,
	\end{equation*}
	where $u^{(1)}(\psi)$ and $u^{(2)}(\psi)$ are, respectively, the solutions of the variational inequalities below:
	\begin{equation*}
		\begin{cases}
			\displaystyle\int_{\omega}\nabla u^{(1)}(\psi)\cdot\nabla(v_1-u^{(1)}(\psi))\,\mathrm{d}y\ge\int_{\omega}f^{(1)}(v_1-u^{(1)}(\psi))\,\mathrm{d}y,\quad\textup{ for all }v_1\in K^{(1)}(\omega;\psi),\\[12pt]
			\displaystyle\int_{\omega}\nabla u^{(2)}(\psi)\cdot\nabla(v_2-u^{(2)}(\psi))\,\mathrm{d}y\ge\int_{\omega}f^{(2)}(v_2-u^{(2)}(\psi))\,\mathrm{d}y,\quad\textup{ for all }v_2\in K^{(2)}(\omega;\psi).
		\end{cases}
	\end{equation*}
}

Writing $r=|y|$, $y\in\overline{\omega}$, we consider the \emph{forged solutions}
\begin{equation*}
	\begin{aligned}
		u^{(1)}_{\textup{f}}&:=-(1-r^2)^2,\\
		u^{(2)}_{\textup{f}}&:=\dfrac{1}{2}(r^6-4r^4+4r^2-1),
	\end{aligned}
\end{equation*}
and we observe that these functions are smooth, radial and vanish along $\partial\omega$. Define $f^{(\alpha)}:=-\Delta u^{(\alpha)}_{\textup{f}}$. For these instances of $f^{(\alpha)}$, we will exploit Algorithm~\ref{alg:1} to compute the numerical solution of the penalised version of the optimal control problem considered in this section, and we expect that the optimal states converge to the corresponding forged solutions at least in $L^2(\omega)$ as $\beta\to 0$ (Theorem~\ref{gamma-4})\footnote{Since this paper mainly concerns the mathematical analysis of the model presented in Section~\ref{Sec:3}, the numerical results reported here are meant to corroborate the theoretical conclusions drawn beforehand. The analysis of Algorithm~\ref{alg:1} with respect to $\rho$, $c$ and the mesh size will be deferred to future work.}.

In what follows, we denote by $n_{\textup{grid}}$ the mesh size associated with a triangulation of $\omega$, and we compare how the mesh size influences the convergence of the error between the numerical solution and the forged solution as $\beta\to 0$. In the experiments presented below, we set $\rho=c=\varepsilon=10^{-5}$.

\begin{table}[H]
	\centering
	\begin{tabular*}{\textwidth}{@{\extracolsep{\fill}}ccccc@{}}
		\toprule
		$\beta$ & $\|u^{(1)}_{\rho,c}(\psi_\beta)-u^{(1)}_{\textup{f}}\|_{L^2(\omega)}$ & $\|u^{(2)}_{\rho,c}(\psi_\beta)-u^{(2)}_{\textup{f}}\|_{L^2(\omega)}$ & $|u^{(1)}_{\rho,c}(\psi_\beta)- u^{(1)}_{\textup{f}}|_{H^1_0(\omega)}$ & $|u^{(2)}_{\rho,c}(\psi_\beta)-u^{(2)}_{\textup{f}}|_{H^1_0(\omega)}$\\
		\midrule
		$10^{-1}$ & $2.2879003\times10^{-1}$ & $8.3176311\times10^{-2}$ & $6.3666993\times10^{-1}$ & $3.3172513\times10^{-1}$\\
		$10^{-2}$ & $6.9447810\times10^{-2}$ & $7.0115142\times10^{-2}$ & $2.4279249\times10^{-1}$ & $2.8474955\times10^{-1}$ \\
		$10^{-3}$ & $1.2671657\times10^{-2}$ & $1.7389758\times10^{-2}$ & $1.6206976\times10^{-1}$ & $1.6670129\times10^{-1}$ \\
		$10^{-4}$ & $4.5508943\times10^{-3}$ & $4.4057268\times10^{-3}$ & $1.5982179\times10^{-1}$ & $1.5561489\times10^{-1}$ \\
		$10^{-5}$ & $4.2848861\times10^{-3}$ & $3.5920631\times10^{-3}$ & $1.5980818\times10^{-1}$ & $1.5545775\times10^{-1}$  \\
		$10^{-6}$ & $4.2093786\times10^{-3}$ & $3.5374228\times10^{-3}$ & $1.5980318\times10^{-1}$ & $1.5545171\times10^{-1}$  \\
		$10^{-7}$ & $4.2007427\times10^{-3}$ & $3.5329861\times10^{-3}$ & $1.5980282\times10^{-1}$ & $1.5545140\times10^{-1}$  \\
		$10^{-8}$ & $4.2000934\times10^{-3}$ & $3.5325002\times10^{-3}$ & $1.5980281\times10^{-1}$ & $1.5545139\times10^{-1}$  \\
		\bottomrule
	\end{tabular*}
	\caption{Numerical error of the difference between the forged solution and the corresponding approximate optimal state for $n_{\textup{grid}}=32$.}
	\label{tab:32}
\end{table}

\begin{table}[H]
	\centering
	\begin{tabular*}{\textwidth}{@{\extracolsep{\fill}}ccccc@{}}
		\toprule
		$\beta$ & $\|u^{(1)}_{\rho,c}(\psi_\beta)-u^{(1)}_{\textup{f}}\|_{L^2(\omega)}$ & $\|u^{(2)}_{\rho,c}(\psi_\beta)-u^{(2)}_{\textup{f}}\|_{L^2(\omega)}$ & $|u^{(1)}_{\rho,c}(\psi_\beta)- u^{(1)}_{\textup{f}}|_{H^1_0(\omega)}$ & $|u^{(2)}_{\rho,c}(\psi_\beta)-u^{(2)}_{\textup{f}}|_{H^1_0(\omega)}$\\
		\midrule
		$10^{-1}$ & $2.2879292\times10^{-1}$ & $8.4337893\times10^{-2}$ & $6.2539602\times10^{-1}$ & $3.0822101\times10^{-1}$ \\
		$10^{-2}$ & $6.8579292\times10^{-2}$ & $7.2299090\times10^{-2}$ & $2.0152070\times10^{-1}$ & $2.6061392\times10^{-1}$ \\
		$10^{-3}$ & $1.1053797\times10^{-2}$ & $1.8245401\times10^{-2}$ & $8.4935700\times10^{-2}$ & $1.0302181\times10^{-1}$ \\
		$10^{-4}$ & $1.8281952\times10^{-3}$ & $2.4925103\times10^{-3}$ & $8.0346315\times10^{-2}$ & $7.8735798\times10^{-2}$ \\
		$10^{-5}$ & $1.1031036\times10^{-3}$ & $9.9649740\times10^{-4}$ & $8.0293946\times10^{-2}$ & $7.8286064\times10^{-2}$  \\
		$10^{-6}$ & $1.0749194\times10^{-3}$ & $9.1500472\times10^{-4}$ & $8.0293555\times10^{-2}$ & $7.8281440\times10^{-2}$  \\
		$10^{-7}$ & $1.0650380\times10^{-3}$ & $9.0972981\times10^{-4}$ & $8.0293216\times10^{-2}$ & $7.8281064\times10^{-2}$  \\
		$10^{-8}$ & $1.0631561\times10^{-3}$ & $9.0948436\times10^{-4}$ & $8.0293168\times10^{-2}$ & $7.8281015\times10^{-2}$  \\
		\bottomrule
	\end{tabular*}
	\caption{Numerical error of the difference between the forged solution and the corresponding approximate optimal state for $n_{\textup{grid}}=64$.}
	\label{tab:64}
\end{table}

\begin{table}[H]
	\centering
	\begin{tabular*}{\textwidth}{@{\extracolsep{\fill}}ccccc@{}}
		\toprule
		$\beta$ & $\|u^{(1)}_{\rho,c}(\psi_\beta)-u^{(1)}_{\textup{f}}\|_{L^2(\omega)}$ & $\|u^{(2)}_{\rho,c}(\psi_\beta)-u^{(2)}_{\textup{f}}\|_{L^2(\omega)}$ & $|u^{(1)}_{\rho,c}(\psi_\beta)- u^{(1)}_{\textup{f}}|_{H^1_0(\omega)}$ & $|u^{(2)}_{\rho,c}(\psi_\beta)-u^{(2)}_{\textup{f}}|_{H^1_0(\omega)}$\\
		\midrule
		$10^{-1}$ & $2.2885023\times10^{-1}$ & $8.4611590\times10^{-2}$ & $6.2261974\times10^{-1}$ & $3.0212080\times10^{-1}$ \\
		$10^{-2}$ & $6.8423478\times10^{-2}$ & $7.2970996\times10^{-2}$ & $1.8998098\times10^{-1}$ & $2.5452983\times10^{-1}$  \\
		$10^{-3}$ & $1.0709802\times10^{-2}$ & $1.8610396\times10^{-2}$ & $4.9660771\times10^{-2}$ & $7.9756653\times10^{-2}$  \\
		$10^{-4}$ & $1.2926172\times10^{-3}$ & $2.3139519\times10^{-3}$ & $4.1210507\times10^{-2}$ & $4.0599695\times10^{-2}$  \\
		$10^{-5}$ & $3.2338818\times10^{-4}$ & $3.5649444\times10^{-4}$ & $4.1101572\times10^{-2}$ & $3.9623122\times10^{-2}$  \\
		$10^{-6}$ & $2.7670471\times10^{-4}$ & $2.4120055\times10^{-4}$ & $4.1100729\times10^{-2}$ & $3.9612750\times10^{-2}$  \\
		$10^{-7}$ & $2.7355179\times10^{-4}$ & $2.3409076\times10^{-4}$ & $4.1100715\times10^{-2}$ & $3.9612641\times10^{-2}$  \\
		$10^{-8}$ & $2.7241174\times10^{-4}$ & $2.3368685\times10^{-4}$ & $4.1100696\times10^{-2}$ & $3.9612621\times10^{-2}$  \\
		\bottomrule
	\end{tabular*}
	\caption{Numerical error of the difference between the forged solution and the corresponding approximate optimal state for $n_{\textup{grid}}=128$.}
	\label{tab:128}
\end{table}

\begin{table}[H]
	\centering
	\begin{tabular*}{\textwidth}{@{\extracolsep{\fill}}ccccc@{}}
		\toprule
		$\beta$ & $\|u^{(1)}_{\rho,c}(\psi_\beta)-u^{(1)}_{\textup{f}}\|_{L^2(\omega)}$ & $\|u^{(2)}_{\rho,c}(\psi_\beta)-u^{(2)}_{\textup{f}}\|_{L^2(\omega)}$ & $|u^{(1)}_{\rho,c}(\psi_\beta)- u^{(1)}_{\textup{f}}|_{H^1_0(\omega)}$ & $|u^{(2)}_{\rho,c}(\psi_\beta)-u^{(2)}_{\textup{f}}|_{H^1_0(\omega)}$\\
		\midrule
		$10^{-1}$ & $2.2887484\times10^{-1}$ & $8.4687986\times10^{-2}$ & $6.2192963\times10^{-1}$ & $3.0055647\times10^{-1}$\\
		$10^{-2}$ & $6.8394845\times10^{-2}$ & $7.3151230\times10^{-2}$ & $1.8693995\times10^{-1}$ & $2.5303897\times10^{-1}$ \\
		$10^{-3}$ & $1.0628073\times10^{-2}$ & $1.8720793\times10^{-2}$ & $3.5411197\times10^{-2}$ & $7.2792762\times10^{-2}$ \\
		$10^{-4}$ & $1.1843813\times10^{-3}$ & $2.3039894\times10^{-3}$ & $2.1987746\times10^{-2}$ & $2.2269527\times10^{-2}$ \\
		$10^{-5}$ & $1.6259245\times10^{-4}$ & $2.5030946\times10^{-4}$ & $2.1781872\times10^{-2}$ & $2.0376064\times10^{-2}$ \\
		$10^{-6}$ & $8.0422846\times10^{-5}$ & $7.1711064\times10^{-5}$ & $2.1779776\times10^{-2}$ & $2.0354703\times10^{-2}$ \\
		$10^{-7}$ & $7.7083643\times10^{-5}$ & $6.3980109\times10^{-5}$ & $2.1779761\times10^{-2}$ & $2.0354493\times10^{-2}$ \\
		$10^{-8}$ & $7.6757277\times10^{-5}$ & $6.3432147\times10^{-5}$ & $2.1779761\times10^{-2}$ & $2.0354491\times10^{-2}$ \\
		\bottomrule
	\end{tabular*}
	\caption{Numerical error of the difference between the forged solution and the corresponding approximate optimal state for $n_{\textup{grid}}=256$.}
	\label{tab:256}
\end{table}

\begin{table}[H]
	\centering
	\begin{tabular*}{\textwidth}{@{\extracolsep{\fill}}ccccc@{}}
		\toprule
		$\beta$ & $\|u^{(1)}_{\rho,c}(\psi_\beta)-u^{(1)}_{\textup{f}}\|_{L^2(\omega)}$ & $\|u^{(2)}_{\rho,c}(\psi_\beta)-u^{(2)}_{\textup{f}}\|_{L^2(\omega)}$ & $|u^{(1)}_{\rho,c}(\psi_\beta)- u^{(1)}_{\textup{f}}|_{H^1_0(\omega)}$ & $|u^{(2)}_{\rho,c}(\psi_\beta)-u^{(2)}_{\textup{f}}|_{H^1_0(\omega)}$\\
		\midrule
		$10^{-1}$  & $2.2888010\times10^{-1}$ & $8.4707138\times10^{-2}$ & $6.2173915\times10^{-1}$ & $3.0014621\times10^{-1}$\\
		$10^{-2}$  & $6.8389776\times10^{-2}$ & $7.3199761\times10^{-2}$ & $1.8611523\times10^{-1}$ & $2.5265152\times10^{-1}$ \\
		$10^{-3}$  & $1.0608231\times10^{-2}$ & $1.8751846\times10^{-2}$ & $3.0445633\times10^{-2}$ & $7.0889145\times10^{-2}$ \\
		$10^{-4}$  & $1.1591950\times10^{-3}$ & $2.3045929\times10^{-3}$ & $1.2475877\times10^{-2}$ & $1.4052427\times10^{-2}$ \\
		$10^{-5}$  & $1.2819689\times10^{-4}$ & $2.3840044\times10^{-4}$ & $1.2108736\times10^{-2}$ & $1.0767397\times10^{-2}$ \\
		$10^{-6}$  & $3.4519929\times10^{-5}$ & $3.2947693\times10^{-5}$ & $1.2104930\times10^{-2}$ & $1.0726492\times10^{-2}$ \\
		$10^{-7}$  & $2.9832570\times10^{-5}$ & $2.1986427\times10^{-5}$ & $1.2104893\times10^{-2}$ & $1.0726079\times10^{-2}$ \\
		$10^{-8}$  & $2.9590907\times10^{-5}$ & $2.1607611\times10^{-5}$ & $1.2104893\times10^{-2}$ & $1.0726075\times10^{-2}$ \\
		\bottomrule
	\end{tabular*}
	\caption{Numerical error of the difference between the forged solution and the corresponding approximate optimal state for $n_{\textup{grid}}=512$.}
	\label{tab:512}
\end{table}

For each $n_{\textup{grid}}\in\{32, 64, 128, 256, 512\}$, we observe that the errors decrease monotonically as $\beta\to 0$, and that a plateau is reached. From the qualitative point of view, the results in Tables~\ref{tab:32}--\ref{tab:512} are corroborated by Figure~\ref{bench:states:ctrl} below, which track the changes in the optimal control and optimal states as $\beta\to 0$.

\begin{figure}[H]
	\centering
	\begin{subfigure}[t]{0.38\linewidth}
		\includegraphics[width=0.95\linewidth]{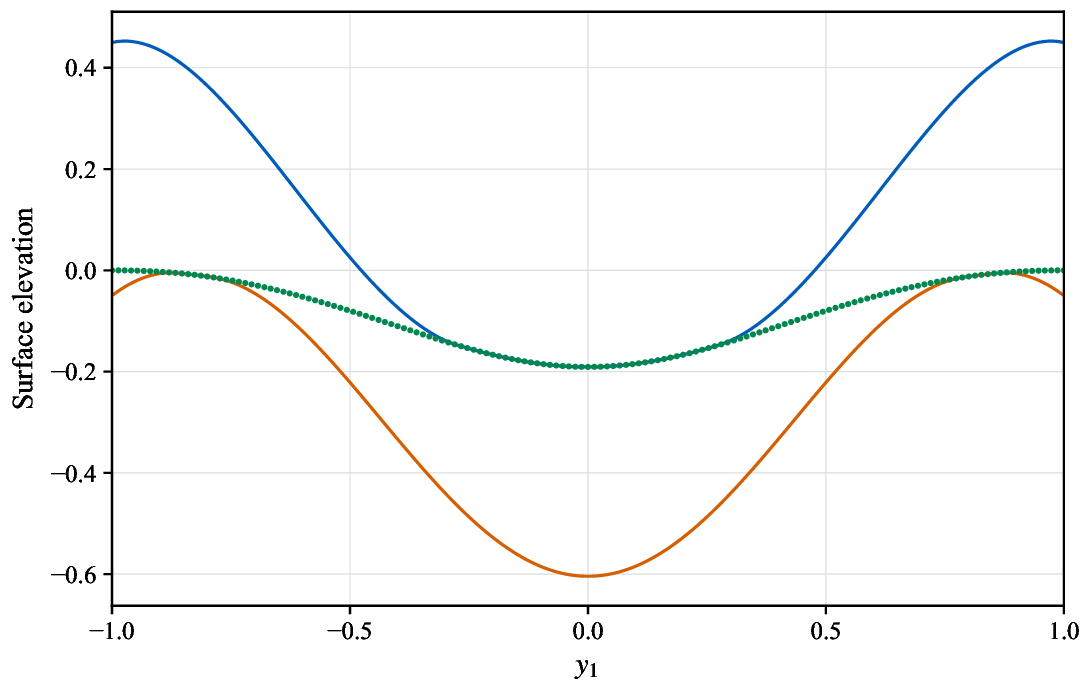}
		\subcaption{$\beta=10^{-1}$}
	\end{subfigure}%
	\hspace{2.0cm}
	\begin{subfigure}[t]{0.38\linewidth}
		\includegraphics[width=\linewidth]{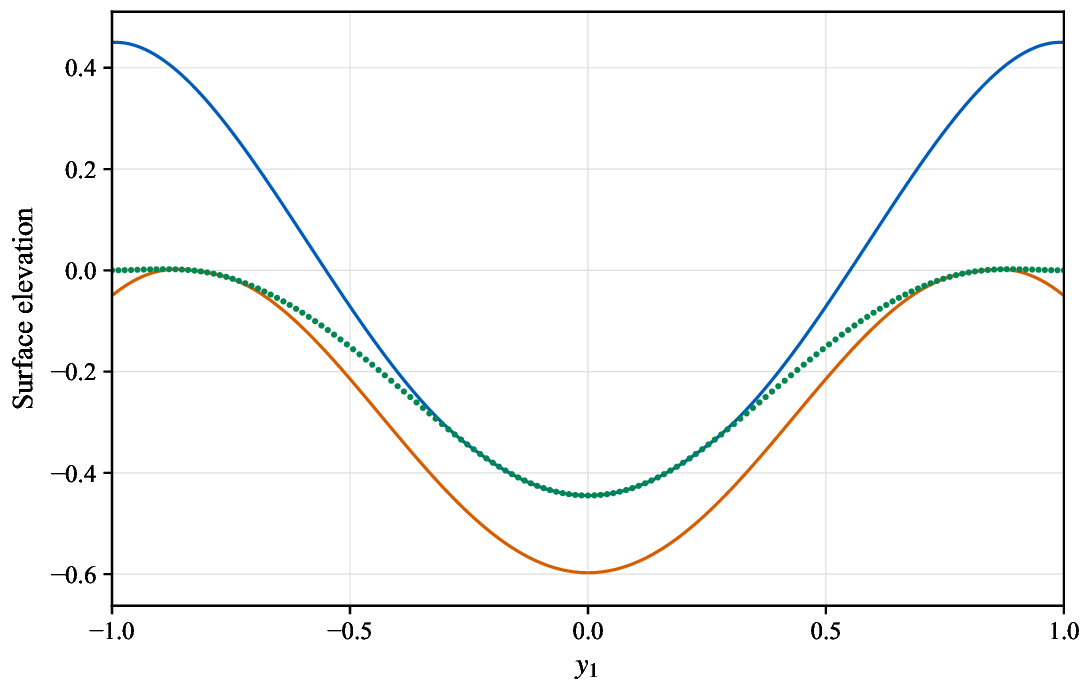}
		\subcaption{$\beta=10^{-2}$}
	\end{subfigure}%
\end{figure}
\begin{figure}[H]
	\ContinuedFloat
	\centering
	\begin{subfigure}[t]{0.38\linewidth}
		\includegraphics[width=0.95\linewidth]{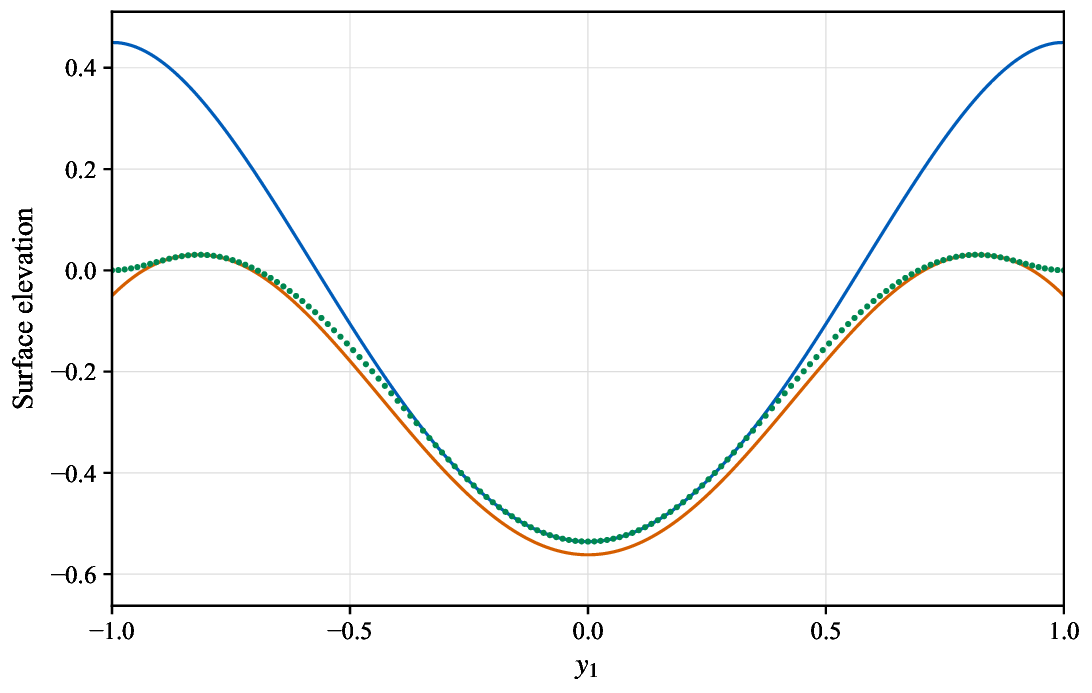}
		\subcaption{$\beta=10^{-3}$}
	\end{subfigure}%
	\hspace{2.0cm}
	\begin{subfigure}[t]{0.38\linewidth}
		\includegraphics[width=\linewidth]{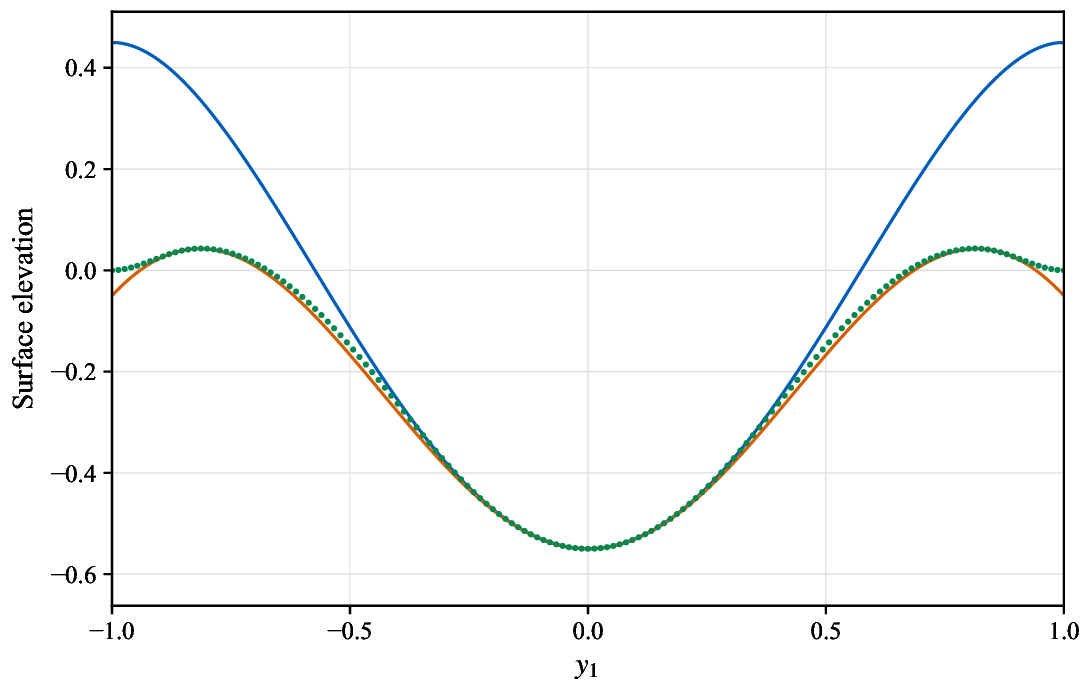}
		\subcaption{$\beta=10^{-8}$}
	\end{subfigure}%
	\caption{Cross sections corresponding to $y_2=0$ of the optimal control $\psi_\beta$ and optimal states $u^{(\alpha)}_{\rho,c}(\psi_\beta)$ for $n_{\mathrm{grid}}=512$ and decreasing values of $\beta$. The green line represents the optimal control $\psi_\beta$, that plays the role of optimal interface, the blue line represents the optimal state $u^{(1)}(\psi_\beta)$ and the orange line represents the optimal state $u^{(2)}(\psi_\beta)$.}
	\label{bench:states:ctrl}
\end{figure}

Figure~\ref{bench:fig:plateaus} below shows that, for $\beta=10^{-8}$, the values $\|u^{(\alpha)}_{\rho,c}(\psi_\beta)-u^{(\alpha)}_{\textup{f}}\|_{L^2(\omega)}$ and $|u^{(\alpha)}_{\rho,c}(\psi_\beta)- u^{(\alpha)}_{\textup{f}}|_{H^1_0(\omega)}$ in Tables~\ref{tab:32}--\ref{tab:512} decrease as $n_{\textup{grid}}$ increases.
\begin{figure}[H]
	\centering
	\begin{subfigure}[t]{0.38\linewidth}
		\centering
		\includegraphics[width=1.0\linewidth]{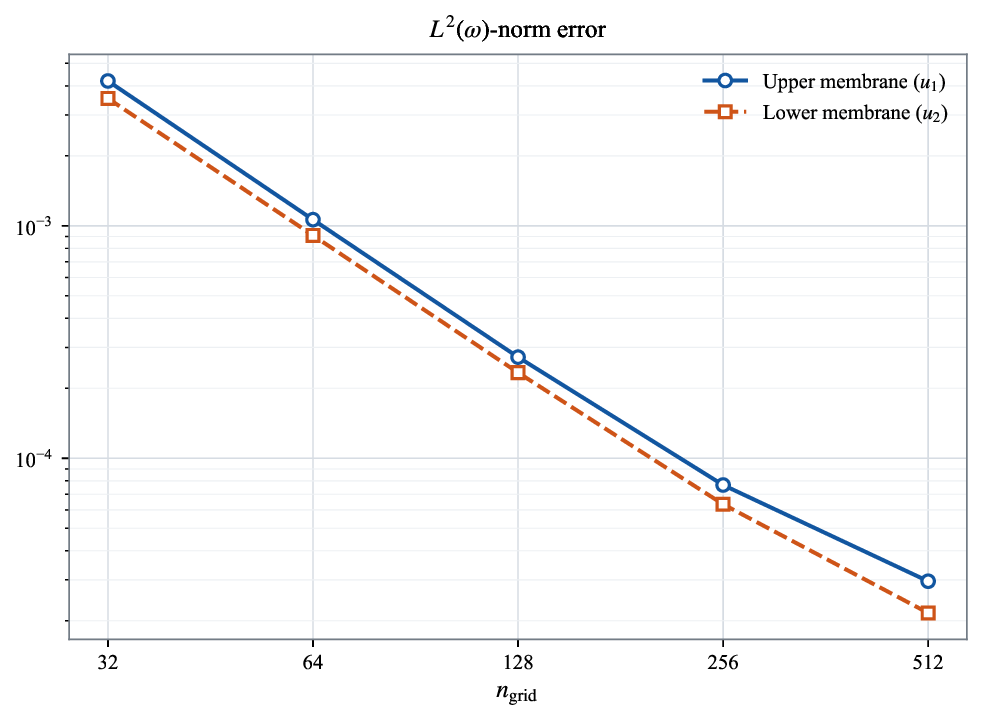}
		\subcaption{$\|u^{(\alpha)}_{\rho,c}(\psi_\beta)-u^{(\alpha)}_{\textup{f}}\|_{L^2(\omega)}$}
	\end{subfigure}%
	\hspace{2.0cm}
	\begin{subfigure}[t]{0.38\linewidth}
		\centering
		\includegraphics[width=1.0\linewidth]{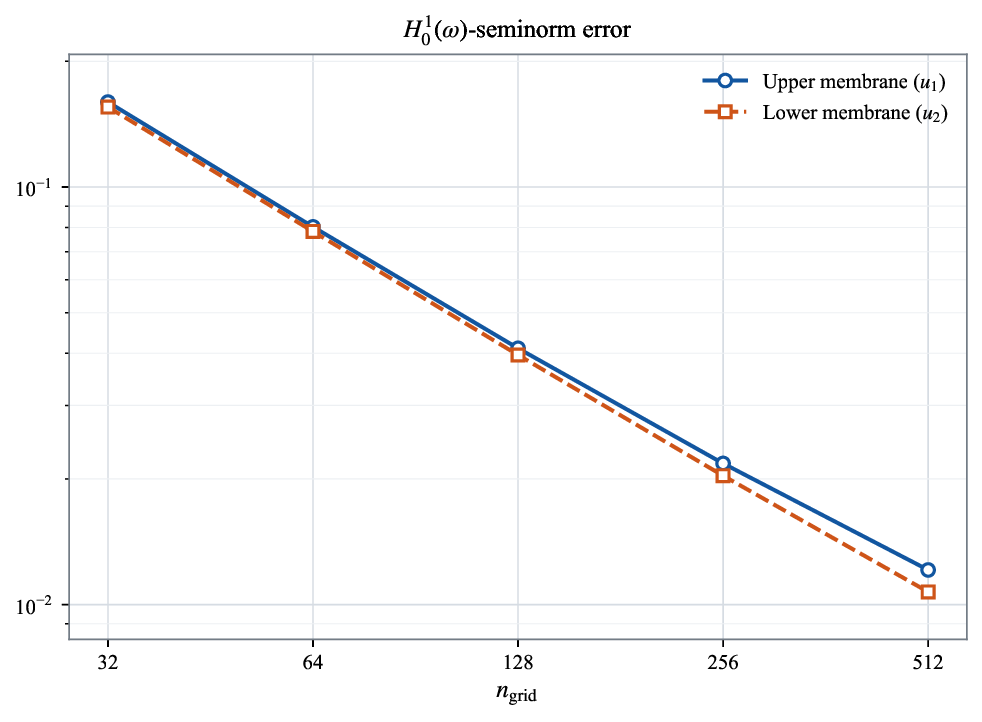}
		\subcaption{$|u^{(\alpha)}_{\rho,c}(\psi_\beta)-u^{(\alpha)}_{\textup{f}}|_{H^1_0(\omega)}$}
	\end{subfigure}%
	\caption{Errors between the forged solution and the approximate optimal states corresponding to $\beta=10^{-8}$.}
	\label{bench:fig:plateaus}
\end{figure}

\addtocontents{toc}{\protect\setcounter{tocdepth}{0}}
\section*{Conclusions and Commentary}
We proposed a new modelling paradigm for two-body unilateral contact in which the separating interface is the primary variable of an optimal control problem, rather than the active set of a coupled variational inequality. By means of a $\Gamma$-convergence argument, we showed that the coincidence set is approximated by a family of surfaces that are the graphs of functions with $H^2$ regularity. Numerical experiment corroborate the theoretical findings. The analysis of the numerical scheme, in particular error estimates with respect to the penalisation and regularisation parameters and to the mesh size, is deferred to future work.

\addtocontents{toc}{\protect\setcounter{tocdepth}{2}}
\addtocontents{toc}{\protect\setcounter{tocdepth}{0}}
\section*{Declarations}

\subsection*{Data Availability Statement}
Data sharing not applicable to this article as no datasets were generated or analysed during the current study.

\subsection*{AI Usage}
Generative Predictive Text chatbots were used to \emph{solely} check the grammar.

\subsection*{Funding Statement}

The research of Z.L. and P.P. is currently being supported by the following agencies:
\begin{itemize}
	\item National Natural Science Foundation of China (NSFC), grant number W2533011;
	\item ZJ Talent Program of the Guangdong Province, grant number 2024QN11X057; 
	\item Peacock Talent Program - Type~C of the Shenzhen Municipality;
	\item University Development Fund of The Chinese University of Hong Kong, Shenzhen.
\end{itemize}

\subsection*{Conflict of Interest Disclosure}
All authors certify that they have no affiliations with or involvement in any organization or entity with any competing interests in the subject matter or materials discussed in this manuscript.

\subsection*{Ethics Approval Statement}
Not applicable.

\subsection*{Patient Consent Statement}

Not applicable.

\subsection*{Permission to Reproduce Material from other Sources}
None of the presented material was drawn from other sources.

\subsection*{Clinical Trial Registration}

Not applicable.

\subsection*{Authors' Contribution}
All authors equally contributed to the realisation of this manuscript.

\addtocontents{toc}{\protect\setcounter{tocdepth}{1}}

\bibliographystyle{abbrvnat} 
\bibliography{references}

\addresseshere

\newpage
\printglossary[type=mathsymbols, title={List of Mathematical Symbols}]
\printglossary[type=physsymbols, title={List of Physics Symbols}]

\end{document}